\documentclass[11pt]{article}
\usepackage{times,fancyhdr}
\usepackage{amsmath}
\usepackage{latexsym}
\usepackage{amssymb}
\usepackage{amsfonts}
\usepackage{amsthm}
\title{DIFERENTIAL GALOIS THEORY \\AND MECHANICS}
\author{ J.-F. POMMARET  \\ CERMICS, Ecole des Ponts Paris-Tech, France  \\(http://cermics.enpc.fr/$\sim$pommaret/home.html)}
\date{  }
\begin{document}
\maketitle

\thispagestyle{empty}

\noindent
{\bf ABSTRACT} \\

The {\it classical Galois theory} deals with certain finite algebraic extensions and establishes a bijective order reversing correspondence between the intermediate fields and the subgroups of a group of permutations called the {\it Galois group} of the extension. It has been the dream of many mathematicians at the end of the nineteenth century to generalize these results to systems of algebraic partial differential (PD) equations and the corresponding finitely generated differential extensions, in order to be able to add the word {\it differential} in front of any classical statement. The achievement of the Picard-Vessiot theory by E. Kolchin between 1950 and 1970 is now well known.  \\

The purpose of this paper is to sketch the general theory for such differential extensions and algebraic pseudogroups by means of new methods 
mixing {\it differential algebra}, {\it differential geometry} and {\it algebraic geometry}. As already discovered by E. Vessiot in 1904 through 
the use of {\it automorphic systems}, a concept never acknowledged, the main point is to notice that the Galois theory (old and new) is mainly 
a study of {\it principal homogeneous spaces} (PHS) for algebraic groups or pseudogroups. Hence, all the formal theory of PD equations developped by D.C. Spencer around 1970 must be used together with modern algebraic geometry, in particular tensor products of rings and fields.  \\

However, the combination of these new tools is not sufficient and we have to create the analogue for Lie pseudogroups of the so-called {\it invariant derivations} introduced by A. Bialynicki-Birula after 1960 in the study of algebraic groups and {\it fields with derivations}. We shall finally prove the usefulness of the resulting {\it differential Galois theory} through striking applications to mechanics, revisiting {\it shell theory}, {\it chain theory}, the {\it Frenet-Serret formulas} and the integration of {\it Hamilton-Jacobi equations}.\\

\vspace{1cm}

\noindent
{\bf FOREWORD}\\

The {\it classical Galois theory} deals with certain finite algebraic extensions and establishes a bijective order reversing correspondence between the intermediate fields and the subgroups of a group of permutations called the {\it Galois group} of the extension.  \\

It has been the dream of many mathematicians at the end of the nineteenth century to generalize these results to systems of linear or algebraic 
ordinary differential (OD) or partial differential (PD) equations and the corresponding finitely generated differential extensions, in order to be able to add the word {\it differential} in front of any classical statement. Among the tentatives, we may quote the {\it Picard-Vessiot theory} dealing with differential extensions having finite transcendence degree, where the Galois group is an {\it algebraic group} that can be considered as a linear algebraic group of matrices. We may also quote the {\it Drach-Vessiot theory} dealing with differential extensions having an infinite transcendence degree but a finite differential transcendence degree where the Galois group is an {\it algebraic Lie pseudogroup}. The achievement of the Picard-Vessiot theory by E. Kolchin and coworkers between $1950$ and $1970$ is now well known.  \\

The purpose of this chapter is to sketch the general theory for arbitrary partial differential extensions and algebraic Lie pseudogroups by means of new methods mixing {\it differential algebra}, {\it differential geometry} and {\it algebraic geometry}. As already discovered by Vessiot in $1904$ through the use of {\it automorphic systems}, a concept still neither known nor acknowledged, the main point is to notice that the Galois theory (old and new) is mainly a study of {\it principal homogeneous spaces} (PHS) for algebraic groups or pseudogroups. Hence all the modern formal theory of OD or PD equations (D.C. Spencer and coworkers around 1970) must be used together with the modern algebraic geometry missing in the work of Kolchin, in particular tensor products of rings and fields.  \\

However, as will be shown by means of explicit counterexamples, the combination of these new tools is not sufficient and we have to create the analogue for Lie pseudogroups of the so-called {\it invariant derivations} introduced by A. Bialynicki-Birula after 1960 in the study of algebraic groups and {\it fields with derivations}.  \\

After recalling the mathematical foundations of the resulting {\it differential Galois theory}, our main purpose will be to prove its usefulness through striking applications to mechanics, revisiting in particular {\it shell theory}, {\it chain theory}, the {\it Frenet-Serret formulas} and the integration of {\it Hamilton-Jacobi equations}.\\

\newpage
\noindent{\bf 1) INTRODUCTION}  \\

Evariste Galois died on may 311832, at the age of $21$ in a duel. Though he introduced the word " {\it group} " in mathematics for the firs time in $1830$, his work has only been known fifteen years later. Then group theory attracted more and more people, slowly passing from the field of pure algebra to the field of differential algebra, with applications ranging from the domain of pure geometry to the domain of differential geometry. \\

A major step ahead has been achieved by Sophus Lie in 1880 with the introduction of {\it Lie groups of transformations} and, in $1890$, with the understanding that these groups of transformations were in fact only examples of a wider class, now called {\it Lie pseudogroups of transformations}, that is groups of transformations solutions of a system of OD or PD equations, in general non-linear and of rather high order. Let us illustrate this point of view with a few examples that will be used later on in a quite different setting. The group $y=ax+b$ of affine transformations of the real line can be considered after differentiating twice as the set of solutions of the second order OD equation $y_{xx}=0$ with standard notations, on the invertibility condition $y_x=a\neq 0$. However, it is not so evident that the group of projective transformations $y=(ax+b)/(cx+d)$ is made by the invertible solutions $y_x\neq 0$ of the third order schwarzian OD equation $(y_{xxx}/y_x) - \frac{3}{2}(y_{xx}/y_x)^2=0$ and it is a pure chance that an explicit integration can be exhibited by any student after some work. Hence the concet of {\it parameter} which is crucial in the first approach has no longer any meaning in the second approach. Going on this way in ${\mathbb{R}}^2$, volume preserving transformations are defined by the jacobian condition $y^1_1y^2_2-y^1_2y^2_2=1$ while complex transformations are defined by the linear Cauchy-Riemann PD equations $y^1_1-y^2_2=0, y^1_2 + y^2_1=0$ and one may introduce similarly contact transformations of ${\mathbb{R}}^3$ or symplectic transformations of ${\mathbb{R}}^4$ and so on, with no longer any parameter involved.  \\

With more details, the idea is to consider a group $G$ as a manifold of dimension $p$ with local coordinates $a=(a^1, ... ,a^p)$ indexed by greek indices $\rho, \sigma, \tau$, with a {\it composition law} $G\times G \rightarrow G:(a,b)\rightarrow ab$, an {\it invers law} $GÊ\rightarrow G: a \rightarrow a^{-1}$ and an identity $e \in G$ such that $a(bc)=(ab)c, ae=ea=a, aa^{-1}=a^{-1}a=e, \forall a,b,c \in G$. The group $G$ is said to {\it act} on a manifold $X$ with local coordinates $x=(x^1, ..., x^n)$ if there is an {\it action map} $X\times G \rightarrow X$ or rather its {\it graph} $X\times G \rightarrow X \times X:(x,a) \rightarrow (x,y=ax=f(x,a))$. The point $x$ is called the {\it source} of the transformation while the point 
$y$ is called the {\it target}. The action is said to be {\it free} if its graph is injective and {\it transitive} if its graph is surjective. Moreover, $X$ is said to be a {\it principal homogeneous space} (PHS) for $G$ if the graph is an isomorphism. The set $Gx=\{ax \mid a\in G \}$ is called the 
{\it orbit} of $x$ under $G$ and the action is said to be {\it effective} if $ax=x, \forall x \in X \Rightarrow a=e$. \\

Let $T$ be the tangent vector bundle of vector fields on $X$, $T^*$ be the cotangent vector bundle of 1-forms on $X$ and ${\wedge}^sT^*$ be the vector bundle of s-forms on $X$ with usual bases $\{dx^I=dx^{i_1}\wedge ... \wedge dx^{i_s}\}$ where we have set $I=(i_1< ... <i_s)$. Also, let $S_qT^*$ be the vector bundle of symmetric q-covariant tensors. Moreover, if  $\xi,\eta\in T$ are two vector fields on $X$, we may define their {\it bracket} $[\xi,\eta]\in T$ by the local formula $([\xi,\eta])^i(x)={\xi}^r(x){\partial}_r{\eta}^i(x)-{\eta}^s(x){\partial}_s{\xi}^i(x)$ leading to the {\it Jacobi identity} $[\xi,[\eta,\zeta]]+[\eta,[\zeta,\xi]]+[\zeta,[\xi,\eta]]=0, \forall \xi,\eta,\zeta \in T$. We have also the useful formula $[T(f)(\xi),T(f)(\eta)]=T(f)([\xi,\eta])$ where $T(f):T(X)\rightarrow T(Y)$ is the tangent mapping of a map $f:X\rightarrow Y$. Finally, we may introduce the {\it exterior derivative} $d:{\wedge}^rT^*\rightarrow {\wedge}^{r+1}T^*:\omega={\omega}_Idx^I \rightarrow d\omega={\partial}_i{\omega}_Idx^i\wedge dx^I$ with $I=\{i_1< ... <i_r\}$ and we have $d^2=d\circ d\equiv 0$ in the {\it Poincar\'{e} sequence}:\\
\[  {\wedge}^0T^* \stackrel{d}{\longrightarrow} {\wedge}^1T^* \stackrel{d}{\longrightarrow} {\wedge}^2T^* \stackrel{d}{\longrightarrow} ... \stackrel{d}{\longrightarrow} {\wedge}^nT^* \longrightarrow 0  \]
Defining the {\it algebraic bracket} $j_q([\xi,\eta])=\{j_{q+1}(\xi),j_{q+1}(\eta)\}$ with a slight abuse of language and recalling the {\it Spencer operator} $D:J_{q+1}(E)\rightarrow T^*\otimes J_q(E): f_{q+1} \rightarrow j_1(f_q)-f_{q+1}$ with $(Df_{q+1})^k_{\mu ,i}=
{\partial}_if^k_{\mu} - f^k_{\mu +1_i}$ for any vector bundle $E$, we may set $E=T$ and define a {\it differential bracket} of Lie algebra on $J_q(T)$ by the formula:  \\
\[     [{\xi}_q, {\eta}_q]= \{ {\xi}_{q+1},{\eta}_{q+1}\} + i(\xi)D{\eta}_{q+1} - i(\eta)D{\xi}_{q+1}, \,\,\,  \forall {\xi}_q,{\eta}_q\in J_q(T)  \]
which does not depend on the lifts and where $i()$ is the standard {\it interior multiplication} of a $1$-form by a vector [19,22,23].  \\

We now recall two results of Lie that will be of constant use in this chapter:   \\

\noindent
$\bullet$ {\it First fundamental theorem}: \\
The orbits $x=f(x_0, a)$ satisfy $\partial x^i/\partial a^{\sigma}={\theta}^i_{\rho}(x){\omega}^{\rho}_{\sigma} (a)$ with $det(\omega)\neq 0$. The vector fields ${\theta}_{\rho}={\theta}^i_{\rho}(x){\partial}_i$ are called {\it infinitesimal generators} of the action and are linearly independent over the constants when the action is effective.  \\

\noindent 
$\bullet$ {\it Second fundamental theorem}: \\
If $\{{\theta}_1, ... ,{\theta}_p\}$ is a basis of the infinitesimal generators of the effective action of a Lie group $G$ on $X$, then $[{\theta}_{\rho},{\theta}_{\sigma}]=c^{\tau}_{\rho\sigma}{\theta}_{\tau}$ where the $c^{\tau}_{\rho\sigma}$ are the {\it structures constants} of the {\it Lie algebra} ${\cal{G}}=T_e(G)$.  \\ 

Coming back to the work of Lie and Vessiot on what is now called {\it Lie pseudogroup} while {\it Lie groups} correspond to actions of groups, we denote by $aut(X)$ the pseudogroup of all local diffeomorphisms of $X$ and define the sub-fibered manifold ${\Pi}_q(X,X)\subset J_q(X\times X)$ by the condition $det(y^k_i)\neq 0$. \\

\noindent
{\bf DEFINITION 1.1}: We may say that $\Gamma \subset aut(X)$ is a Lie pseudogroup of transformations defined by a system of OD or PD equations ${\cal{R}}_q\subset {\Pi}_q(X,X)$ with $n$ independent variables $x$ and the same number of unknowns $y$ if, whenever $y=f(x)$ and $z=h(y)$ are two local invertible transformations solutions of this system that can be composed in $aut(X)$, then $z=h\circ f (x)$ and $x=f^{-1}(y)=g(y)$ are again solutions. \\

As seen on the previous examples, such a definition is totally useless in actual practice unless one can provide an explicit form for the generic solutions of the defining system of ODE or PDE which may be highly nonlinear. Let us now introduce {\it two} manifolds, namely the {\it source manifold} $X$ with $dim(X)=n$ and the {\it target manifold} $Y$ with $dim(Y)=m$ and local coordinates $y=(y^1, ... ,y^m)$.  Hence, considering a target transformation $\bar{y}=g(y)\in \Gamma$, a map $f:X\rightarrow Y:x\rightarrow y=f(x)$ and using the chain rule for derivatives, we obtain successively ${\bar{y}}_x=\frac{\partial g}{\partial y} y_x, {\bar{y}}_{xx}=\frac{\partial g}{\partial y} y_{xx}+ \frac{{\partial}^2g}{\partial y \partial y}y_xy_x$ in a symbolic way, and so on. Accordingly, we may look for {\it differential invariants}, namely functions $\Phi (y,y_x,y_{xx}, ...)$ preserved by the action of $\Gamma$ on the target, that is such that $\Phi (\bar{y},{\bar{y}}_x, {\bar{y}}_{xx},...)=\Phi (y,y_x,y_{xx}, ...)$. We may also look for infinitesimal transformations of the target $y \rightarrow \bar{y}=y + t {\eta}(y) + ...$ where $t$ is a small constant parameter and extend it to ${\bar{y}}_x=y_x+t \frac{\partial \eta}{\partial y}y_x, {\bar{y}}_{xx}=y_{xx}+t (\frac{\partial \eta}{\partial y}y_{xx}+\frac{{\partial}^2\eta}{\partial y\partial y}y_xy_x),... $. If $\mu=({\mu}_1,...,{\mu}_n)$ is a multi-index with length $\mid \mu \mid={\mu}_1 + ... + {\mu}_n$ and $\mu + 1_i=({\mu}_1, ..., {\mu}_{i-1},{\mu}_i+1, {\mu}_{i+1},...,{\mu}_n)$, we may set $y_q=\{y^k_{\mu} \mid 0\leq k\leq m,0\leq \mid \mu \mid \leq q,y_0=y\}$ and introduce the following {\it formal derivative} on functions of $(x,y_q)$ in order to get functions of $(x,y_{q+1})$:  \\
\[  \frac{d}{dx^i}=\frac{\partial }{\partial x^i} + y^k_{\mu +1_i} \frac{\partial }{\partial y^k_{\mu}}  \]
We may thus define the $q$-prolongation ${\rho}_q(\eta)$ of the target infinitesimal transformation $\eta={\eta}^k(y)\frac{\partial}{\partial y^k}\in T(Y)$ by the formula ${\rho}_q(\eta)=d_{\mu}{\eta}^k\frac{\partial}{\partial y^k_{\mu}}$. Now, if $\Theta\subset T$ denotes the set of infinitesimal transformations of $\Gamma \subset aut(X)$, one can prove that $[R_q,R_q] \subset R_q$ when $R_q=id_q^{-1}V({\cal{R}}_q)$ if we set $id_q=j_q(id)$ for the identity transformation $y=x$. It follows that such a condition can be checked by means of computer algebra, contrary to the condition $[\Theta,\Theta]\subset \Theta$. In the formula for ${\rho}_q(\eta)$, we may replace the derivatives of $\eta$ with respect to $y$ by a section ${\eta}_q\in R_q(Y)$ and denote by $\sharp({\eta})_q$ the vertical vector in $V({\cal{R}}_q)$ obtained. One can then prove that we have the important formula:  \\
\[   [\sharp ({\eta}_q),\sharp ({\zeta}_q)]=\sharp ([{\eta}_q,{\zeta}_q]     \]
Applying the Frob\'{e}nius theorem, on the resulting distribution of vertical vector fields on ${\cal{R}}_q$, we may obtain a fundamental set $\{ {\Phi}^{\tau}(y_q)$ of funtionnaly invariant {\it differential invariants}.

Setting now $Y=X$ and $y=x$, we may summarize the previous results as follows:  \\

\noindent
$\bullet$ {\it First fundamental result of Vessiot}: \\
Any source transformation commutes with any target transformation and exchanges therefore among them te differential invariants of a fundamental set. Patching coordinates, we may therefore obtain a {\it natural bundle} ${\cal{F}}$ over $X$ of order $q$, also called {\it bundle of geometric objects} of order $q$, both with a section $\omega$ in such a way that $\Gamma=\{f\in aut(X)\mid {j_q(f)}^{-1}(\omega)=\omega \}$. \\

\noindent
{\bf EXAMPLE 1.2}: With $n=2, m=2$, let us introduce the manifolds $X$ with local coordinate $(x^1,x^2)$ and let $Y$ be a copy of $X$ with local coordinates $(y^1,y^2)$. We may consider the {\it algebraic Lie pseudogroup} $\Gamma \subset aut(X)$ of (local, invertible) transformations of $X$ preserving the $1$-form $\alpha= x^2dx^1$ and thus also the $2$-form $\beta=dx^1\wedge dx^2$, that is to say made up by transformations $y=f(x)$ solutions of the Pfaffian system  $ y^2dy^1=x^2dx^1$ and thus $dy^1\wedge dy^2=dx^1\wedge dx^2$. Equivalently, we have to look for the invertible solutions of the algebraic first order involutive system ${\cal{R}}_1\subset {\Pi}_1(X,X)$ defined  by the first order involutive system of algebraic PD equations in Lie form:\\ 
\[   {\Phi}^1\equiv y^2y^1_1=x^2, \hspace{5mm} {\Phi}^2\equiv y^2y^1_2= 0 \hspace{5mm}\Rightarrow  \hspace{5mm}{\Phi}^3\equiv \frac{\partial (y^1,y^2)}{\partial (x^1,x^2)}=y^1_1y^2_2-y^1_2y^2_1=1  \]
We let the reader check that the corresponding natural bundle over $X$ is ${\cal{F}}=T^*{\times }_X{\wedge}^2T^*$ with section $\omega = (\alpha,\beta)$ and we notice that $d_1{\Phi}^1-d_2{\Phi}^1+{\Phi}^3=0$, that is $d\alpha + \beta=0$. By chance one can obtain the generic solution $y^1=f(x^1), \hspace{2mm} y^2=x^2/(\partial f(x^1)/\partial x^1)$ where $f(x^1)$ is an arbitrary (invertible) function of one variable.
The linearized system over the target $Y$ is:  \\
\[  y^2\frac{\partial {\eta}^1}{\partial y^1}+{\eta}^2=0, \frac{\partial {\eta}^1}{\partial y^2}=0, \frac{\partial {\eta}^1}{\partial y^1}+ \frac{\partial {\eta}^2}{\partial y^2}=0  \]
 Now, with $n=1,m=2$, introducing a manifold $X$ of dimension $n=1$ and a {\it different} manifold $Y$ of dimension $2$ with a map $y=f(x)$ while considering the corresponding transformations of the jets $(y^1,y^2,y^1_x,y^2_x, ...)$, we obtain the distribution generated by:  \\
 \[   \{ {\theta}_1= \frac{\partial}{\partial y^1}, \,\, {\theta}_2= y^2\frac{\partial}{\partial y^2}-y^1_x\frac{\partial}{\partial y^1_x}+y^2_x\frac{\partial}{\partial y^2_x}, \,\, {\theta}_3= y^1_x\frac{\partial}{\partial y^2_x}  \}  \]
which has the only generating differential invariant $\Phi\equiv {\bar{y}}^2  {\bar{y}}^1_x=y^2y^1_x$ because its generic rank is $3$ and it is easy to check the commutation relations $[{\theta}_1,{\theta}_2]=0, [{\theta}_1,{\theta}_3]=0, [{\theta}_2,{\theta}_3]=-2{\theta}_3$. The corresponding natural bundle with local coordinates $(x,u)$ is $T^*=T^*(X)$ because it has the transition rules $(\bar{x}=\varphi(x), \bar{u}=u/\frac{\partial (\varphi(x)}{\partial x})$. \\

The next proposition is important but its proof is out of the scope of this book [22]:  \\

\noindent
{\bf PROPOSITION 1.3}: For any function $\Phi \in {\Pi}_q(X,X)$ we have:  \\
\[  \sharp({\eta}_{q+1}) d_i\Phi=d_i(\sharp({\eta}_q).\Phi) - y^k_i\sharp(D{\eta}_{q+1}(\frac{\partial}{\partial y^k})).\Phi    \]

If $\Phi$ is a differential invariant at order $q$ and ${\eta}_{q+1}\in R_{q+1}$, then $D{\eta}_{q+1}\in T^*\otimes R_q$ over the target and the right member vanishes, that is $d_i\Phi$ is a differential invariant at order $q+1$. Consider now a maximal number of formal linear combinations of the $d_i{\Phi}^{\tau}$ that do not contain jets of strict order $q$. We can always suppose that they begin with a leading term equal to $1$ and we may apply $\sharp(R_{q+1})$ in order to find a contradiction unless the other coefficients of the combinations are killed by $\sharp(R_q)$ and are thus only functions of the $\Phi$, a result leading to identities of the symbolic form:  \\
\[                  I(j_1(\Phi))\equiv A(\Phi) d\Phi +B(\Phi)=0   \]
Taking now the reciprocal images with respect to the corresponding sections when $Y$ is a copy of $X$, we obtain the vector bundles:  \\
\[            T=id^{-1}(V(X\times Y), \,\,\,  R_q=id_q)^{-1}(V({\cal{R}}_q), \,\,\, F={\omega}^{-1}(V({\cal{F}}))      \]
and an operator:  \\
\[  {\cal{D}}: T \Rightarrow F: \xi \rightarrow {\cal{L}}(\xi)\omega= \frac{d}{dt}j_q(exp(t\xi))^{-1}(\omega)\mid_{t=0}=\Omega   \]
where ${\cal{L}}(\xi)$ is the {\it Lie derivative} of a geometric object with respect to a vector field $\xi$. We notice that ${\cal{D}}$ is a {\it Lie operator} in the sense that ${\cal{D}}\xi=0, {\cal{D}}\eta=0 \Rightarrow {\cal{D}}[\xi,\eta]=0$.\\
Finally, introducing the vector bundle $J_q(T)$ of $q$-jets of $T$ with sections ${\xi}_q$ over $\xi$ transforming like $j_q(\xi)$, we have $F=J_q(T)/R_q$ and we may introduce its {\it Medolaghi form} [19,22]:  \\
\[    {\Omega}^{\tau}\equiv - L^{\tau\mu}_k(\omega (x)){\xi}^k_{\mu} + {\xi}^r{\partial}_r{\omega}^{\tau}(x)=0   \]
with $0< \mid \mu \mid \leq q$, a result showing that the coefficients of ${\cal{D}}$ only depend on $\omega$ but the last, exacily like the standard Lie derivative of tensors. In particular, if $\Gamma$ contains the translations $\xi=cst$, then $\omega$ {\it must} be locally constant. \\

\noindent
$\bullet$ {\it Second fundamental result of Vessiot}:  \\
The map ${\pi}^{q+1}_q:R_{q+1} \rightarrow R_q$ is surjective for a {\it general section} $\omega$ of ${\cal{F}}$ if and only if it is indeed surjective for the {\it special section} and if the {\it Vessiot structure equations}:  \\
\[      I(j_1(\omega))=c(\omega)   \]
are satisfied, where the section $\omega \rightarrow c(\omega)$ is invariant under any diffeomorphism and only depends on a certain number of constants called {\it structure constants} satisfying purely algebraic equations $J(c)=0 $ called {\it Jacobi conditions} by analogy with the case of Lie groups, even though there is no Lie algebra structure behind.  \\

\noindent
{\bf EXAMPLE 1.4}: In the preceding example we have $d\alpha=-\beta$ when $\alpha=x^2dx^1,\beta=dx^1\wedge dx^2$ and thus one structure constant $c=-1$ only. However, we may choose $\alpha=dx^1, \beta=dx^1\wedge dx^2$ a choice leading to $c=0$ and to the new (non-isomorphic) pseudogroup $y^1=x^1+a, y^2=x^Ž+g(x^1)$. Moreover, we have proved in [26] that, in the case of a Riemann structure, that is when ${\cal{F}}=S_2T^*$ and $\omega$ is such that $det(\omega)\neq 0$, there are indeed {\it two} structure constants but $J(c)$ is linear in such a way that they must be equal and there is finally the {\it only} constant of the constant Riemann curvature condition. We do not believe that such a result or even situation is known. \\

We have explained the deep contribution of Vessiot to the formal theory of Lie pseudogroups, made as early as in 1903 [19,22,26] but still neither known nor acknowledged today, and we are in position to sketch the other important contribution of Vessiot to the differential Galois theory, made as early as in 1904 [29] but still neither known nor acknowledged today.\\

First of all, a group of permutation can be represented as a group of matrices with entries equal to $0$ or $1$, having a single $1$ in each row or column. For example, the group ${\mathfrak{S}}_2$ of permutations in $2$ variables is $\{(12)\rightarrow (12), (12) \rightarrow (21)\}$, the first matrix is the identity $2\times 2$ matrix while the second is the anti-diagonal $2\times 2$ matrix. Then it is well known that any algebraic group can be realized by a linear algebraic group of matrices. Also, an algebraic pseudogroup must be defined by a system of algebraic OD or PD equations, the above permutation group being defined by the algebraic equations $( y^1 +y^2=x^1  + x^2, y^1y^2=x^1x^2)$. However, the pseudogroup defined by the condition that the Jacobian matrix should be of the form:  \\
\[  \left ( \begin{array}{ccc}
    1 & 0 & 0  \\
    A & 1  &  0  \\
    0 & 0 & e^A
    \end{array}  \right )  \]                                                                                                                                                                                                                                                                                                                                                                                                                                                                                                                                                                                                                                                                                                                                                                                                                                                                                                                                                                                                                                                                                                                                                                                                                                                                                                                                                                                                                                                                                                                                                                                                                                                                                                                                                                              
is {\it not} an algebraic pseudogroup because one of the defining OD equations is easily seen to be $ y^3_3 - \exp(y^2_1)=0$. We have the inclusions:  \\
\[     finite \,\,\,  groups  \,\,\, \subset \,\,\, algebraics  \,\;                 groups   \, \, \,  \subset  \,\,\,  algebraic  \,\,\, pseudogroups  \,\,\,    \]

Now, if $K\subset L$ are fields, then $L$ can be considered as a vector space over $K$ and we shall set $\mid L/K\mid=dim_K(L)$. Also,if the extension $L/K$ is not an algebraic extension, that is if $L$ contains elements which are not algebraic over $K$, that is which are not roots of a polynomial with coefficients in $K$, then the maximum number of algebraically independent such elements is well defined, only depends on $L/K$ and is called the {\it transcendence degree} $trd(L/K)$ of $L/K$. More generally, a {\it differential integral domain} $A$ is a ring without any divisor of zero and with $n$ commuting derivations ${\partial}_1,...,{\partial}_n$ such that ${\partial}_i(a+b)={\partial}_ia + {\partial}_i b, {\partial}_i(ab)=({\partial}_ia)b + a {\partial}_ib, \forall a,b\in A, \forall i=1,...,n$. Such a definition can easily be extended in order to define a {\it differential field} by introducing the field $K=Q(A)$ of quotients of $A$ while setting ${\partial}_i(a/b)=({\partial}_ia)b-a{\partial}_ib)/b^2, \forall a,b\in A$ with $b\neq 0$. Similarly, if $L/K$ is a {\it differential extension}, that is if $K\subset L$ and the derivations of $K$ are induced by those of $L$, then the maximum number of elements of $L$ which are not differentially algebraic, that is which are not solutions of a differential polynomial PD equation with coefficients in $K$, is well defined and is called the {\it differential transcendence degree} 
$diff\, trd (L/K)$. The first basic idea of Vessiot has been first to establish a kind of "{\it classification} " of the {\it differential Galois theory}, namely:\\

\noindent
$\bullet$ CLASSICAL GALOIS THEORY : $\mid L/K \mid < \infty $.  \\
systems of algebraic equations $\leftrightarrow$  finite groups.  \\

\noindent
$\bullet$ PICARD-VESIOT THEORY: $trd(L/K)< \infty, diff\, trd(L/K)=0$.  \\
Systems of algebraic OD or PD equations $\leftrightarrow$ algebraic groups.  \\

\noindent
$\bullet$ DRACH-VESSIOT THEORY: $trd(L/K)=\infty, diff\,trd(L/K)<\infty$.  \\
Systems of algebraic PD equations $\leftrightarrow$ algebraic pseudogroups.  \\

The next crucial step is to establish a specific link between the systems and the groups and this will be the heart of this chapter. In order to sketch the underlying idea in the most elementary way, let us review for a few lines the classical Galois theory and see how one could add the word 
"{\it differential} " in front of the concepts. If $L/K$ is an algebraic extension, we denote by $iso(L/K)$ the set (care) of isomorphisms 
$\varphi:L\rightarrow M$ of $L$ into another field $M$ containing $K$ and such that $\varphi(a)=a, \forall a\in K$ and $\varphi(a+b)=\varphi(a) + \varphi(b), \varphi(ab)=\varphi(a)\varphi(b), ,\forall a,b\in L$. We denote by $aut(L/K)$ the group of automorphisms of $L$ fixing $K$ and by $inv(\Gamma)$ the subfield of $L$ fixed by a group $\Gamma \subset aut(L/K)$ with $\mid \Gamma \mid$ elements. Classical Galois theory deals with {\it Galois extensions} and we recall the following three equivalent definitions that can be found in any textbook [1,28]:  \\

\newpage
\noindent
{\bf DEFINITION 1.5}: \\
\noindent
1) $L/K$ is a Galois extension if $iso(L/K)=aut(L/K)=\Gamma$ with $\mid \Gamma \mid = \mid L/K\mid$.  \\
2  $L/K$ is a Galois extension if $inv(\Gamma)=K$ with $\Gamma = aut(L/K)$.  \\
3) $L/K$ is a Galois extension if $L$ is the {\it splitting field} of an irreducible polynomial with coefficients in $K$, that is $L$ is obtained by ajoining to $K$ all the roots of such a polynomial which can be thereore decomposed into linear factors over $L$.  \\

Once these definitions are assumed, the two main results of the classical Galois theory useful for applications are the following {\it fundamental theorem} and its corollary:  \\

\noindent
{\bf THEOREM 1.6}: When $L/K$ is a Galois extension, there is a bijective order reversing {\it Galois correspondence} between intermediate fields $K\subset K' \subset L$ and subgroups $id \subset{\Gamma}'\subset \Gamma=\Gamma(L/K)$ given by: \\
\[  K' \longrightarrow {\Gamma}'=aut(L/K'), \hspace{2cm}  {\Gamma}' \longrightarrow K'=in({\Gamma}')  \]

\noindent
{\bf COROLLARY 1.7}: Let $L/K$ be a Galois extension and $M$ be an arbitrary extension of $K$. If $L$ and $M$ are contained in a bigger field $N$ and we denote by $(L,M)$ the {\it composite field} of $L$ and $M$ in $N$, that is the smallest subfield of $N$ containing both $L$ and $M$, then $(L,M)/M$ is a Galois extension and there is an isomorphism $\Gamma((L,M)/M)\simeq \Gamma(L/(L\cap M))$. Moreover, $L$ and $M$ are linearly disjoint over $L\cap M$ in $(L,M)$.\\

\noindent
{\bf REMARK 1.8}: The additional well known result saying that $K'/K$ is again a Galois extension if and only if ${\Gamma}'\lhd \Gamma$, that is ${\Gamma}'$ is a {\it normal subgroup} of $\Gamma$ will not be considered here because the study of the {\it normalizer} of a Lie pseudogroup $\Gamma$ in $aut(X)$ is one of the most difficult problem to be found in the formal theory of Lie pseudogroups [19,26].\\ 

The first and second previous properties cannot be extended because there is no reason at all that the trnsformations of the Lie group/pseudogroup of invariance of the system do respect {\it any kind} of extension. As for the last definition, a system of OD or PD equations has in general an infinite number of solutions that cannot even be explicitly described or added in general. In particular, we should like to strongly react against the abstract fashion a few people are using "{\it universal extension}, a kind of huge reserve into which one ould put all the solutions of all systems of algebraic OD or PD equations , exactly like people use to do in in algebra with the field of complex numbers through the so-called {\it fundamental theorem of algebra} [14,15,27]. \\

Before providing the striking answer given by Vessiot, let us examine the auxiliary  though preliminary problem of how to define a {\it differential extension} by relating it only to the formal theory of systems of OD or PD equations. If $K$ is a differential field as above and $(y^1,...,y^m)$ are indeterminates over $K$, we transform the polynomial ring $K\{y\}={lim}_{q\rightarrow \infty}K[y_q]$ into a differential ring by introducing as usual the {\it formal derivations} $d_i={\partial}_i+y^k_{\mu+1_i}\partial/\partial y^k_{\mu}$ and we shall set $K<y>= Q(K\{y\})$ as field of quotients.  \\

\noindent
{\bf DEFINITION 1.9}: We say that $\mathfrak{a}\subset K\{y\}$ is a {\it differential ideal} if it is stable by the $d_i$, that is if  $d_ia\in\mathfrak{a}, \forall a \in \mathfrak{a}, \forall i=1,...,n$. We shall also introduce the {\it radical} ideal $rad(\mathfrak{a})=\{a\in A\mid \exists r,a^r\in \mathfrak{a}\}\supseteq \mathfrak{a}$ and say that $\mathfrak{a}$ is a {\it perfect} (or {\it radical}) differential ideal if $rad(\mathfrak{a})=\mathfrak{a}$. We say that $\mathfrak{p}\subset K\{y\}$ is a {\it prime} differential ideal if it is a prime ideal and a differential ideal. If $S$ is any subset of a differential ring $A$, we shall denote by $\{S\}$ the differential ideal generated by $S$ and introduce the (non-differential) ideal 
${\rho}_r(S)=\{ d_{\nu}a \mid a\in S, 0 \leq \mid\mu\mid \leq r\}$ in $A$.  \\

\noindent
{\bf LEMMA 1.10}: If $\mathfrak{a}\subset A$ is differential ideal, then $rad(\mathfrak{a}) $ is a differential ideal containing 
$\mathfrak{a}$.  \\

\noindent
{\it Proof}: If $d$ is one of the derivations, we have $ a^{r-1}da=\frac{1}{r}da^r \in \{a^r\}$ and thus:  \\
\[ (r-1) a^{r-2}(da)^2 + a^{r-1}d^2a \in \{a^r\}\Rightarrow a^{r-2}(da)^3\in \{a^r\},... \Rightarrow (da)^{2r-1} \in \{a^r\}   \]
\hspace*{12cm}   Q.E.D.   \\

We shall say that a differential extension $L=Q(K\{y\}/\mathfrak{p})$ is a {\it finitely generated} differential extension of $K$ and we may define the {\it evaluation epimorphism} $K\{y\} \rightarrow K\{\eta \}\subset L$ with kernel $\mathfrak{p}$ where $\eta$ or $\bar{y}$ is the residual image of $y$ modulo $\mathfrak{p}$. In particular, the following Lemma will be used in the next important Theorem:  \\

\noindent
{\bf LEMMA 1.11}: If $\mathfrak{p}$ is a prime differential ideal of $K\{y\}$, then, for $q$ sufficiently large, there is a polynomial $P\in K[y_q]$ such that $P\notin {\mathfrak{p}}_q$ and :   \\
\[          P{\mathfrak{p}}_{q+r} \subset rad ({\rho}_r({\mathfrak{p}}_q)) \subset {\mathfrak{p}}_{q+r}, \hspace {1cm}  \forall r\geq 0  \] 

\noindent
{\bf THEOREM  1.12}: ({\it Primality test}) Let ${\mathfrak{p}}_q\subset K[y_q]$ and ${\mathfrak{p}}_{q+1}\subset K[y_{q+1}]$ be prime ideals  such that ${\mathfrak{p}}_{q+1}={\rho}_1({\mathfrak{p}}_q)$ and ${\mathfrak{p}}_{q+1}\cap K[y_q]={\mathfrak{p}}_q$. If the symbol $g_q$ of the algebraic variety ${\cal{R}}_q$ defined by ${\mathfrak{p}}_q$ is $2$-acyclic and if its first prolongation $g_{q+1}$ is a vector bundle over ${\cal{R}}_q$, then $\mathfrak{p}={\rho}_{\infty}({\mathfrak{p}}_q)$ is a prime differential ideal with $\mathfrak{p} \cap K[y_{q+r}]={\rho}_r({\mathfrak{p}}_q), \forall r\geq 0 $.  \\

\noindent
{\bf EXAMPLE 1.13}: With $n=2$ and $K=\mathbb{Q}$, let us consider the two differential polynomials $P_1\equiv y_{22}-\frac{1}{3}(y_{11})^3, P_2\equiv y_{12}- \frac{1}{2}(y_{11})^2$ in $K\{y\}$. We have $d_2P_2-D_1P_1-y_{11}d_1P_2\equiv 0$ and the differential ideal 
$\mathfrak{p}=\{P_1,P_2\}\subset K\{y\}$ is prime as we have indeed:   \\
\[  K\{y\}/ \mathfrak{p} \simeq K[y,y_1,y_2, y_{11}, y_{111},...]      \]

\noindent
{\bf COROLLARY  1.14}: Every perfect differential ideal of $\{y\}$ can be expressed in a unique way as the non-redundant intersection of a finite number of prime differential ideals.  \\

\noindent
{\bf COROLLARY 1.15}: ({\it Differential basis}) If $\mathfrak{r}$ is a perfect differential ideal of $K\{y\}$, then we have $\mathfrak{r}=rad({\rho}_{\infty} ({\mathfrak{r}}_q))$ for $q$ sufficiently large.  \\

\noindent
{\bf PROPOSITION  1.16}: If $\zeta$ is differentially algebraic over $K<\eta>$ and $\eta$ is differentially algebraic over $K$, then $\zeta$ is differentially algebraic over $K$. Setting $\xi=\zeta - \eta$, it follows that, if $L/K$ is a differential extension and $\xi,\eta \in L$ are both differentially algebraic over $K$, then $\xi + \eta$, $\xi\eta$ and $d_i\xi$ are differentially algebraic over $K$.  \\

If $L=Q(K\{y\}/\mathfrak{p})$, $M=Q(K\{z\}/\mathfrak{q})$ and $N=Q(K\{y,z\}/\mathfrak{r})$ are such that $\mathfrak{p}=\mathfrak{r}\cap K\{y\}$ and $\mathfrak{q}=\mathfrak{r}\cap K\{z\}$, we have the two towers $K\subset L\subset N$ and $K\subset M\subset N$ of differential extensions and we may therefore define the new tower $K \subseteq L\cap M \subseteq <L,M> \subseteq N$. However, if only $L/K$ and $M/K$ are known and we look for such an $N$ containing both $L$ and $M$, we may use the universal property of tensor products an deduce the existence of a differential morphism $L{\otimes}_KM\rightarrow N$ by setting $d(a\otimes b)=(d_La) \otimes b+a \otimes (d_Mb)$ whenever $d_L\mid K=d_M\mid K=\partial$. The construction of an abstract {\it composite} differential field amounts therefore to look for a prime differential ideal in $L{\otimes}_K M$ which is a direct sum of integral domains [20,30].  \\

\noindent
{\bf DEFINITION  1.17}: A differential extension $L$ of a differential field $K$ is said to be {\it differentially algebraic} over $K$ if every element of $L$ is differentially algebraic over $K$. The set of such elements is an intermediate differential field $K' \subseteq L$, called the {\it differential algebraic closure} of $K$ in $L$. If $L/K$ is a differential extension, one can always find a maximal subset $S$ of elements of $L$ that are differentially transcendental over $K$ and such that $L$ is differentially algebraic over $K<S>$. Such a set is called a {\it differential transcedence basis} and the number of elements of $S$ is called the {\it differential transcendence degree} of $L/K$.  \\  

\noindent
{\bf THEOREM  1.18}: The number of elements in a differential basis of $L/K$ does not depent on the generators of $L/K$ and his value 
is $difftrd(L/K)=\alpha$. Moreover, if $K\subset L \subset M$ are differential fields, then $difftrd(M/K)=difftrd(M/L) + difftrd(L/K)$.  \\

\noindent
{\bf THEOREM 1.19}: If $L/K$ is a finitely generated differential extension, then any intermediate differential field $K'$ between $K$ and $L$ is also finitely generated over $K$.  \\

We shall now slightly transform the third property of Definition $1.4$ into:   \\

\noindent
4) $L/K$ is a Galois extension if  $L{\otimes}_KL\simeq  L \oplus ...  \oplus L$ with $\mid L/K \mid$ terms.  \\                                                                                                                                                                         

However, we have $L\oplus ... \oplus L =L{\otimes}_{\mathbb{Q}} (\mathbb{Q} \oplus ... \oplus \mathbb{Q} )$ with $\mid L/K \mid$ terms in the direct sum of fields equal to $\mathbb{Q}$ which is isomorphic to $\mathbb{Q}[\Gamma]$ because $\Gamma$ is a group of permutations which splits entirely over $\mathbb{Q}$ as we have already seen by exhibiting square invertible matrices with coefficients equal to $0$ or $1$ only. As a byproduct, we obtain the isomorphism:  \\
\[        L{\otimes}_{\mathbb{Q}} L\simeq L {\otimes }_k  k[\Gamma]   \]
where each member is a direc sum of fields [4,20,30].\\

Let finally $k\subset K \subset L$ be fields and consider an irreducible algebraic set $X$ or {\it variety} defined over $K$ by a prime ideal $\mathfrak{p}\subset K[y]$. We may denote as usual by $K[X]=K[y]/\mathfrak{p}$ the ring of polynomial functions on $X$ {\it which is an integral domain} and introduce its field of quotients $L=K(X)=Q(K[y]/\mathfrak{p})$. We shall say that $X$ is the {\it model variety} of the extension. Similarly, if $G$ is an algebraic group defined on $k$, we denote by $k[G]$ the ring of polynomial functions on $G$. Accordingly, if $X$ is a PHS for $G$, we have $X\times X\simeq X\times G$ and obtain therefore the 
{\it fundamental isomorphism} also called {\it Hopf duality}:  \\
\[    Q(L{\otimes}_KL) \simeq Q(L{\otimes}_kk[\Gamma])    \]
where both members are direct sums of fields. It follows that classical Galois theory is a theory of algebraic PHS and we thus obtain the key idea of Vessiot obtained as early as in 1904 in a clever paper where each chapter is studying the previous classification:  \\

\noindent
DIFFERENTIAL GALOIS THEORY IS A THEORY OF DIFFERENTIAL ALGEBRAIC PHS FOR ALGEBRAIC PSEUDOGROUPS.  \\

Of course the usual definition of of a PHS saying that, if $y=f(x)$ and $\bar{y}=\bar{f}(x)$ are two solutions of the defining system of algebraic OD or PD equations, then there exists one and only one transformation $\bar{y}=g(y)\in \Gamma$ such that $\bar{f}=g\circ f$ is {\it totally useless} in actual practice, thoug it can be checked sometimes. Nevertheless, we may state with Vessiot:  \\

\noindent 
{\bf DEFINITION 1.20}: A system of OD or PD equations having such an above property is called an {\it automorphic system} for $\Gamma$.ÊThe corresponding differential extension is called a {\it differential automorphic extension}. \\

\noindent
{\bf REMARK 1.21}: Contrary to the situation existing in he classical Galois theory, the irreducible components of a PHS may not be themselves PHS for subgroups. We notice that the equation $y^4+1=0$ defines a PHS for the group $\{ \bar{y}=\epsilon y\mid {\epsilon}^4=1 \}$ made up by the four roots of unity $(\epsilon=1,i,-1,-i)$. However, we have the identity $y^4+1\equiv (y^2-2y+1)(y^2+2y+1)$.  \\

\noindent
{\bf EXAMPLE 1.22}: Coming back to example $1.2$ with $n=1,m=2$, the defining system $y^2y^1_x=\omega \in K$ is an automorphic system for the algebraic pseudogroup $\Gamma=\{{\bar{y}}^1=g(y^1), {\bar{y}}^2= y^2/(\partial g/ \partial y^1)\}$. Indeed, giving $y^1=f^1((x)$, we get $y^2=\omega / {\partial}_xf^1$ provided the derivative is non-zero. Hence we can get $x=h(x^1)$ by the implicit function theorem and set ${\bar{y}}^1= \bar{f}(h(y^1))=g(y^1)$. However, the new system $y^2y^1_x - y^1y^2_x=\omega \in K$ cannot be integrated in an explicit way and it does not seem evident to prove that it is an automorphic system for the algebraic pseudogroup preserving the $1$-form $y^2dy^1-y^1dy^2$ and the $2$-form $dy^1\wedge dy^2$. There is no reason to tansform also "$x$" which does not appear explicitly in $K$. Finally, if we set: \\
\[ K=\mathbb{Q}<y^2y^1_x> \subset K'= \mathbb{Q}<y^2y^1_x, y^2_x> \subset L=\mathbb{Q}<y^1,y^2> \]
we get at once the subpseudogroup ${\Gamma}'=\{{\bar{y}}^1=y^1 + a, {\bar{y}}^2=y^2\}$ preserving $K'$ with $a=cst$ and cannot obtain a Galois correspondence because it leaves invariant the intermediate differential field $K"=\mathbb{Q}<y^2y^1_x, y^2>$ which strictly contains $K'$. On the contrary, if we choose now $K'=\mathbb{Q}<y^2y^1_x, y^2_x/y^2>$, we get ${\Gamma}'=\{{\bar{y}}^1=ay^1+b, {\bar{y}}^2=(1/a)y^2$. Hence there are only two possible ways to escape from such a contradiction: one is to say that the Galois correspondence does not exist in the differential Galois theoy while the other is to say that not all intermediate differential fields can be chosen.  \\

\noindent
{\bf EXAMPLE 1.23}: ({\it Picard-Vessiot}) As another reason not to believe in the Picard-Vessiot theory of Kolchin and others, let us prove that it cannot even allow to study the simplest second order OD equation $y_{xx}=0$. For this, following Vessiot, let us copy twice this equation in order to obtain the automorphic system $P_1\equiv y^1_{xx}=0, P_2 \equiv y^2_{xx}=0$ for the action of the linear group $GL(2)$ with ${\bar{y}}^1=ay^1+by^2, {\bar{y}}^2=cy^1+dy^2$ where $a,b,c,d=cst$. With $K=\mathbb{Q}$, we have $L=Q(K\{y^1\}/\{P_1,P_2\})\simeq K(y^1,y^2, y^1_x,y^2_x)$ which is a differential field with $d_xy^k=y^k_x,d_xy^k_x=0$ providing therefore a differential automorphic extension $L/K$ with $trd(L/K)=4$. Choosing the intermediate ifferential field $K'=K(y^1_x)$, we obtain easily the subgroup ${\Gamma}'$ of $GL(2)$ defined by ${\bar{y}}^1=y^1, {\bar{y}}^2= cy^1+dy^2$ which preserves $K"=K(y^1)$ strictly containing $K'$. Finally, choosing $K'=K(y^1y^2_x - y^2y^1_x) $, we get ${\Gamma}'=SL(2)$ and thus $K"=K'$. \\

Let us consider now a quite different problem brought by the above examples. Indeed, for certain groups or subgroups, we have found Lie groups with constant parameters, namely $a=cst$ in the first and $a,b,c,d=cst$ in the second. such a result has no formalmeaning because no explicit integration can be achieved in general. Hence, setting $C=cst(K)=\{a\in K\mid {\partial}_ia=0, \forall i=1,...,n\}$ for the {\it subfield of constants} of $K$, the question is now where to find these constant parameters. The next result will prove the confusion that has been done on the concept of constants. When the PHS $X$ for a Lie group $G$ is $G$ tself, it has been first discovered by A. Bialynicki-Birula in $1960$ and presented in two very difficult papers [3,4].\\

\noindent
{\bf THEOREM 1.24}: ({\it PHS revisited}) The group parameters are constant on $X \times X$.  \\

\noindent
{\it Proof}: Let $X$ be PHS for $G$ and consider the {\it graph isomorphism} $X\times G \simeq X \times X$ under the action of $G$ on $X$ already considered. We have $dim(X)=dim(G) \Rightarrow n=p$ and thus:  \\
\[    det     \left ( \begin{array}{cc} 
                        1 & 0 \\
                                  \frac{\partial f}{\partial x}   & \frac{\partial f}{\partial a}     
                                     \end{array}   \right)
= det(\frac{\partial f}{\partial a}) \neq 0  \]                                                                                                                                                                                                                                                                                                                                                                                                                                                                                                                                                                                                                                                                                                                                                                                                                                                                                                                                                                                                                                                                                                                                                                                                                                                                                                                                                                  
Using this isomorphism, we may exhibit $p$ functions $a=\varphi (x,y)$ when $dim(G)=p$ in such a way to have the $dim(X)=n$ identities $y\equiv f(x,\varphi(x,y))$. Let $\delta={\xi}^i(x)\frac{\partial}{\partial x^i}\in \Delta$ be a transformation commuting with all the infinitesimal generators $\theta \in \Theta$ of the action and extend $\delta$ to $X\times X$ or rather $X\times Y$ where $Y$ is a copy of $X$ by setting anew 
$ \delta={\xi}^i(x)\frac{\partial}{\partial x^i} + {\xi}^k(y)\frac{\partial}{\partial y^k}$. Applying to the above identity, we obtain the formula:  \\
\[     {\xi}^k(y)={\xi}^i(x)\frac{\partial f^k}{\partial x^i}(x,a) + (\delta.a^{\tau})\frac{\partial f^k}{\partial a^{\tau}}(x,a)  \]
whenever $a=\varphi(x,y)$. Then $\delta$ commutes with the action, that is $[\delta,\theta]=0$ and thus $\xi (y) = \xi(x)\frac{\partial f}{\partial x}$. But we have $n=p$ and thus $  det(\frac{\partial f}{\partial a})\neq 0  \Rightarrow \delta . a=0  $. Finally, as $[\Delta, \Theta]=0$ and $\Phi$ is an ivariant of thaction, we get $\delta ( \theta \Phi) - \theta (\delta \Phi)=0 \Rightarrow \theta (\delta \Phi)=0$ and thus $\delta \Phi$ {\it must} also be invariant by $\Theta$.  \\
\hspace*{12cm}    Q.E.D.    \\

We discover that the group parameters must be "{\it constants} " but not necessarily killed by the ${\partial}_i$ as in the standard sense of the word.  \\

\noindent
{\bf COROLLARY 1.25}: $k[G]\subset cst(Q(L {\otimes}_K L) \Rightarrow Q(L {\otimes}_K L)\simeq Q(L{\otimes}_k cst(Q(L {\otimes}_K L)$. where the first $L$ in the right member is identified with $L\otimes 1$. \\ 

\noindent
{\bf EXAMPLE 1.26}: Let us consider the affine group of the real line and set $z=b^1y+b, y=a^1x + a^2\Rightarrow z=(b^1a^1)x + (b^1a^2+b^2)$ in order to get the group law $b=(b^1,b^2), a=(a^1,a^2) \Rightarrow \bar{a}=ba \Rightarrow ({\bar{a}}^1=b^1a^1,{\bar{a}}^2=b^1a^2+b^2)$. With $K=\mathbb{Q}\subset K'=\mathbb{Q}(a^2/a^1)\subset \mathbb{Q}(a^1,a^2)=L$ we obtain $G'$ by setting $b^2=0$ and get again $K'$. However, with the new $K'=\mathbb(Q)(a^1a^2)$, we get for the new $G'$ the conditions $(b^1)^2=1, b^2=0$ and the invariant field 
$K"=\mathbb{Q}('a^1)^2,a^1a^2, (a^2)^2)$ which is strictly bigger than $K'$. In this case, we may use the reciprocal commuting left and right invariant disributions:  \\ 
\[  \Theta=\{ {\theta}_1=a^1\frac{\partial}{\partial a^1} +a^2 \frac{\partial}{\partial a^2},\,\,{\theta}_2=\frac{\partial}{\partial a^2}\}  \]

\[\Delta=\{{\delta}_1=a^1\frac{\partial}{\partial a^1}, \,\, {\delta}_2=a^1\frac{\partial}{\partial a^2} \}  \]
and we check indeed that ${\delta}^2(a^1a^2)=(a^1)^2 \notin K'$.  \\

\noindent
{\bf EXAMPLE 1.27}: ({\it Picard-Vessiot revisited}) Prolonging the action of $GL(2)$ up to the first order jets, we get:  \\

\[    \left( \begin{array}{cc}  
a & b \\
 c & d  
 \end{array} \right)                   
=
    \left( \begin{array}{cc}  
     {\bar{y}}^1 & {\bar{y}}^1_x\\  {\bar{y}}^2 & {\bar{y}}^2_x  \end{array}  \right)  
   \left(  \begin{array}{cc}
   y^1 & y^1_x  \\
   y^2 & y^2_x
   \end{array} \right)^{-1}  \]                                                                                                                                                                   
that we can write $A=\bar{M}M^{-1}$. Now, the two reciprocal distributions are:  \\
{ \small 
   \[     \Theta= \{ {\theta}_1=y^1 \frac{\partial}{\partial y^1 } + y^1_x \frac{\partial}{\partial y^1_x },
           {\theta}_2     =     y^1 \frac{\partial}{\partial y^2 } + y^1_x \frac{\partial}{\partial y^2_x },
           {\theta}_3  =          y^2 \frac{\partial}{\partial y^1 } + y^2_x \frac{\partial}{\partial y^1_x },
          {\theta}_4   =       y^2 \frac{\partial}{\partial y^2 } + y^2_x \frac{\partial}{\partial y^2_x }   \}      \]
\[  \Delta=\{ {\delta}_1=y^1 \frac{\partial}{\partial y^1 } + y^2\frac{\partial}{\partial y^2 },
                  {\delta}_2     =       y^1 \frac{\partial}{\partial y^1_x } + y^2 \frac{\partial}{\partial y^2_x },
           {\delta}_3=          y^1_x \frac{\partial}{\partial y^1 } + y^2_x \frac{\partial}{\partial y^2 },
          {\delta}_4=          y^1_x \frac{\partial}{\partial y^1_x } + y^2_x \frac{\partial}{\partial y^2_x }   \}      \]   }
\noindent
We check that $\bar{\delta }M=A\delta M$ and thus $\delta A=0$ whenever $det(M)\neq 0$, that is the well known non-zero wronskian condition. We also check that each $\delta$ stabilizes $K"$ but {\it not} $K'$. We finally notice that the action is generically free because the rank of $\Theta$ is equal to $4$ whenever the {\it wronskian condition} $y^1y^2_x-y^2y^1_x \neq 0$ is satisfied.  \\

\noindent
{\bf EXAMPLE  1.28} Revisiting Examples 1.2 and 1.22, we obtain the reciprocal distribution at order $1$:    \\
\[ \Delta= \{{\delta}_1= y^1_x \frac{\partial}{\partial y^1_x} + y^2_x \frac{\partial}{\partial y^2_x}, \,\, {\delta}_2= y^2 \frac{\partial }{\partial y^2_x}  \}    \]
Taking into account that we have a tensor product over $K$, we may thus use the identification ${\bar{y}}^2{\bar{y}}^1_x=y^2y^1_x$ in the formulas and check that each $\delta$ separately kills:   \\
\[  \frac{\partial {\bar{y}}^1}{\partial y^1}= \frac{{\bar{y}}^1_x}{y^1_x}=\frac{y^2}{{\bar{y}}^2}, \,\, \frac{\partial{\bar{y}}^1}{\partial y^2}=0, \frac{\partial {\bar{y}}^2}{\partial y^1}= \frac{{\bar{y}}^1_x{\bar{y}}^2_x-y^1_xy^2_x}{y^1_x{\bar{y}}^1_x}= \frac{y^2{\bar{y}}^2_x - {\bar{y}}^2y^2_x}{y^2y^1_x}
          , \,\,  \frac{\partial {\bar{y}}^ 2}{\partial y^2}= \frac{y^1_x}{{\bar{y}}^1_x}=\frac{{\bar{y}}^2}{y^2}  \]

Now, in order to apply the differential Galois theory to mechanics, we need to answer to the following important question:  \\

\noindent
{\bf QUESTION 1.29}: When a given linear or nonlinear system of algebraic OD or PD equations is given, how to check in a formal way, that is to say without introducing solutions, that it is an automorphic system for an algebraic pseudogroup and can thus be considered as a model differential algebraic variety for a differential automorphic extension ? .  \\

Among the systems of OD/PD equations with $m=n$, we shall distinguish the ones defining Lie pseudogroups. For this, let us introduce the open sub-bundle ${\Pi}_q={\Pi}_q(X,X) \subset J_q(X\times X)$ defined by the condition $det(y^k_i)\neq 0$ when ${\cal{E}}=X\times Y$ and $Y$ is a copie of $X$. We introduce the {\it source projection} ${\alpha}_q:{\Pi}_q \rightarrow X:(x,y_q) \rightarrow x$ and the {\it target projection} ${\beta}_q:{\Pi}_q \rightarrow Y:(x,y_q) \rightarrow y$ while identifying a map $f:X \rightarrow Y$ with its graph $f \rightarrow X\times X\times X$. In particular, we denote by $id:X\rightarrow X\times X \times X:x\rightarrow'x,x)$ is called the {\it identity map}.  \\

\noindent
{\bf DEFINITION 1.30}: A fibered manifold ${\cal{R}}_q\subset {\Pi}_q$ is called a system of {\it finite Lie equations} or a {\it Lie groupoid} of oder $q$ if we have an induced {\it source projection} ${:alpha}_q:{\cal{R}}_q \rightarrow X$ an induced {\it target projection} ${\beta}_q: {\cal{R}}_q \rightarrow X$, an induced {\it identity} $j_q(id)=id_q$ as a section $X \rightarrow   {\cal{R}}_q$, an induced {\it composition} ${\gamma}_q:{\cal{R}}_q {\times}_X {\cal{R}}_q$ where the fibered product is taken with respect to the target projection on the left and to the  source projection on the right, both with an induced {\it inversion} ${\iota}_q{\cal{R}}_q \rightarrow {\cal{R}}_q$ exchanging source and target. The set $\gamma\subset aut(X)$ of (local) solutions of ${\cal{R}}_q$ is called a {\it Lie pseudogroup of order q}. We shall suppose that ${\cal{R}}_q$ is {\it transitive}that is the projectin $({\alpha}_q, {\beta}_q):{\cal{R}}_q \rightarrow X \times X$ is surjective.   \\

With evident notations, we set formally in a pointwise way from left to right: \\
\[ ((x,y,\frac{\partial y}{\partial x}, ...), (y,z, \frac{\partial z}{\partial y},...)) \rightarrow  (x,z,\frac{\partial z}{\partial y}.\frac{\partial y}{\partial x}, ...)       \]
\[  (x,y,\frac{\partial y}{\partial x}, ...)      \rightarrow   (y,x,(\frac{\partial y}{\partial x})^{-1},...)\]
We set $ j_q(g \circ f )=j_q(g) \circ j_q(f)$ whenever the composition is defined and $j_q(f)^{-1}=j_q(f^{-1})$.  \\

\noindent
{\bf DEFINITION  1.31}: Setting $y=x+t{\xi}(x)+ ...$ and liearizing, the pull-back $R_q=id_q^{-1}(V({\cal{R}}_q))\subset J_q(T)$ defines a system of {\it infinitesimal Lie equations} and, as we already saw, we have $[R_q,R_q] \subset R_q \Rightarrow [\Theta, \Theta ] \subset \Theta$.                                                                                                                                                                                                                                                                                                                                                                                                                                                                                                                                       

When $X$ with local coordinates $(x^1,...,x^n)$ and $Y$ wit local coordinates $(y^1,...,y^m)$ and $m$ is no longer equal to $n$, we may use the preceding results, in particular the bracket on sections of $R_q(Y)$ when $\Gamma \subset aut(Y)$ is a Lie pseudogroup of transformations of $Y$ in order to find a fundamental set of generating differential invariants ${\Phi}^{\tau}(y_q)$ at order $q$. Then we know that the $d_i{\Phi}^{\tau}$ are again differential invariants at order $q+1$, though, as we shall see in  many examples of mechanics, new differential invariants may be added at higher order. Accordingly, we just need a criterion in order to stop the adjonction procedure when $q$ is large enough.    \\

\noindent
{\bf REMARK 1.32}: In the case of an algebraic pseudogroup, the use of the Frob\'{e}nius theorem may not provide rational differential invariants. However, when $m=n$, it is possible to exhibit rational differential invariants by using {\it formal translation} techniques first introduced by J. Drach (See [6,9] and [20], p 467 for details).      \\

Using the composition of jets when $m\neq n$ exactly as we did when $m=n$, we have an action morphism with graph:  \\
\[   J_q(X \times Y){\times}_Y {\Pi}_q(Y,Y) \stackrel{(graph)}\longrightarrow  J_q(X\times Y){\times}_X J_q(X \times Y)  \]
and a system ${\cal{A}}_q\subset J_q(X \times Y)$ will be said to be {\it invariant} by the action of the Lie groupoid ${\cal{R}}_q\subset {\Pi}_q(Y,Y)$ if we have the following restricted action morphism with induced graph:   \\
\[     {\cal{A}}_q {\times}_Y {\cal{R}}_q \longrightarrow  {\cal{A}}_q {\times}_X {\cal{A}}_q    \]

As we shall see, the following definition and the two corresponding criteria will become crucial for the applications to mechanics.\\

\noindent
{\bf DEFINITION 1.33}: ${\cal{A}}_q$ is said to be a PHS for ${\cal{R}}_q$ if the above morphism is an isomorphism. Then ${\cal{A}}_q$ is said to be an {\it automorphic system} for ${\cal{R}}_q$ if  ${\cal{A}}_{q+r}={\rho}_r({\cal{A}}_q)$ is a PHS for ${\cal{R}}_{q+r}={\rho}_r({\cal{R}}_q), \forall r\geq 0$.   \\

Of course, if ${\cal{A}}_q$ is an automorphic system for ${\cal{R}}_q$, then, {\it necessarily}, ${\cal{A}}_q$ is a PHS for ${\cal{R}}_q$ 
{\it and} ${\cal{A}}_{q+1}$ is a PHS for ${\cal{R}}_{q+1}$ but such a double condition may not be sufficient, a result proving the importance of the following theorem (See [20], p 330 for the technical proof):  \\

\noindent
{\bf THEOREM 1.34}: ({\it First criterion for automorphic systems}) If an involutive system ${\cal{A}}_q \subset J_q(X\times Y)$ is a PHS for a Lie groupoid ${\cal{R}}_q \subset {\Pi}_q(Y,Y)$ and if ${\cal{A}}_{q+1}={\rho}_1({\cal{A}}_1) \subset J_q(X\times Y)$ is a PHS for the Lie groupoid ${\cal{R}}_{q+1}={\rho}_1({\cal{R}}_q) \subset {\Pi}_{q+1}(Y,Y)$, then ${\cal{R}}_q$ is an involutive system with the same non-zero characters and 
${\cal{A}}_q$ is an automorphic system for ${\cal{R}}_q$.  \\

Similarly, we have (See [20], p 339 for the technical proof):   \\

\noindent
{\bf THEOREM  1.35}: ({\it Second criterion for automorphic systems}) If ${\cal{R}}_q \subset {\Pi}_q(Y,Y)$ is an involutive sytem of finite Lie equations such that the action of ${\cal{R}}_q$ on $J_q(X\times Y)$ is generically free, then the action of ${\cal{R}}_{q+r}={\rho}_r({\cal{R}}_q)$ on $J_{q+r}(X \times Y)$ is generically free and all the differential invariants are generated by a fundamental set of order $q+1$ (care to the order).\\
 
\noindent
{\bf EXAMPLE  1.36}: With $K=\mathbb{Q}, n=1, m=2, q=2$, let $\Gamma$ be the Lie pseudogroup determined by the Lie group $G$ with $3+2=5$ parameters, defined by $\bar{y}=Ay+B$ with $det(A)=1$. We have the involutive defining system of finite Lie equations: \\
\[      \frac{\partial({\bar{y}}^1,{\bar{y}}^2)}{\partial (y^1,y^2)}=1, \frac{{\partial}^2\bar{y}}{\partial y\partial y}=0  \]
with linearized system ${\eta}^1_1+{\eta}^2_2=0, {\eta}^k_{rs}=0, \forall k,r,s=1,2$ having a zero symbol and thus trivially involutive. The only generating differential invariant at order $2$ is $\Phi\equiv y^1_xy^2_{xx}-y^2_xy^1_{xx}$. It follows that $d_x\Phi\equiv y^1_xy^2_{xxx}-y^2_xy^1_{xxx}$ is a differential invariant at order $3$ but we have also $\Psi\equiv  y^1_{xx}y^2_{xxx}-y^2_{xx}y^1_{xxx}$. The involutive distribution $\Theta$ involved has $5$ infinitesimal generators at order $2$, namely:  \\
\[  {\theta}_1= \frac{\partial}{\partial y^1}, \,\,{\theta}_2= \frac{\partial}{\partial y^2}, 
\,\,  {\theta}_3=y^1_x\frac{\partial}{\partial y^1_x}+y^1_{xx}\frac{\partial}{\partial y¬1_{xx}} - y^2_x\frac{\partial}{\partial y^2_x} - y^2_{xx}\frac{\partial}{\partial y^2_{xx}},\]
\[   {\theta}^4= y^2_x\frac{\partial}{\partial y^1_x}+y^2_{xx}\frac{\partial}{\partial y^1_{xx}}, 
\,\,  {\theta}_5= y^1_x\frac{\partial}{\partial y^2_x}+ y^1_{xx}\frac{\partial}{\partial y^2_{xx}}   \]
We notice that the symbol of ${\cal{A}}_3$ vanishes if and only if $y^1_xy^2_{xx}-y^2_xy^1_{xx}\neq 0$ and, in this case, we have thus a generically free action, in agrement with the second criterion. Finally, we notice that ${\cal{A}}_2$ is a PHS for ${\cal{R}}_2$, both fibered manifolds having a fiber dimension equal to $5$ (namely $(2+2+2)-1$ for the first and $(2+4)-1$ for the second), and that ${\cal{A}}_3$ is thus  an automorphic system with the same fiber dimension at any higher order.   \\

It is finally important to notice that the definition of prime differential ideals was not known at all by Drach and Vessiot because it has only been introduced by Ritt after 1930 (See [20,27] for details) and we shall correct the old definition by saying:   \\

\noindent
{\bf DEFINITION 1.37}: A system of algebraic PD equations is {\it irreducible} if it is defined by a {\it prime} differential ideal.  \\

Contrary to the classical Galois theory, the irreducible components of a PHS for a group may not be again at all PHS for subgroups. It follows that Drach, Vessiot or even Kolchin made a confusion between {\it prime} differential ideals and {\it maximum} differential ideals in order to define the Galois group of a system. Hence, the starting point of any differential Galois theory must be an irreducible automorphic system (See the introduction 
of [20] for more details). However, despite this comment, modern works on the Picard-Vessiot (PV) theory or on the differential Galois theory are based on three {\it conceptual misunderstandings}:  \\

\noindent
1) {\it First misunderstanding}: The Galois group of a PV extension cannot be computed in general. On the other end, the system $y^1_{xx}=0,y^2_{xx}=0$ is indeed an irreducible ({\it prime} because {\it linear}) automorphic system for $GL(2)$. Hence, the Galois group of {\it any} linear OD equation of order $m$ must be $GL(m)$.  \\

\noindent
2) {\it Second misunderstanding}: The OD equation $y_{xx}=0$ cannot be treated by Kolchin  because the corresponding automorphic 
extension $L/K$ defined by $y^1_{xx}=0,y^2_{xx}=0$ is such that $L$ contains differential constants other than the ones of $K$, for example $y^1_x$. Hence, the definition of "{\it constants} " never took into account the modern approach of Bialynicki-Birula.   \\

\noindent
3) {\it Third misunderstanding}: {\it Last but not least}, instead of using the concept of {\it algebraic pseudogroups} like Drach and Vessiot, Kolchin used the concept of {\it differential algebraic group} initiated by Ritt [16].     \\

\newpage

\noindent
{\bf  2) SHELL THEORY REVISITED}  \\

Using the previous notations, we study first the case $n=2,m=3$ when $G$ is the group of rigid motions ($3$ translations + $3$ rotations) acting on the space $Y={\mathbb{R}}^3$ with cartesian coordinates $(y^1,y^2,y^3)$ and the space $X={\mathbb{R}}^2$ with coordinates $(x^1,x^2)$ is the parametrizing manifold for the surface considered but we shall insist on the intrinsic aspect by using {\it jet theory}, both fom the differential geometric and algebraic aspects. Meanwhile, we invite the reader to try to imagine how to extend the results below to arbirary $n$ and $m$.  \\
The starting motivation for such a study has been the clever idea of Pierre Oscar Bonnet in $1867$ ([2]) to look for the following vague and difficult problem:  \\

\hspace*{2mm}              {\it HOW TO DETERMINE A SURFACE IN\, ${\mathbb{R}}^3$ UP TO A RIGID MOTION}\\

Of course, translated into the modern language of the previous sections, it just amounts to construct the corresponding {\it automorphic systems}, that is to exhibit a generating set of differential invariants ${\Phi}^{\tau}$, the corresponding bundle $\cal{F}$ of geometric objects and a section $\omega$ both with its compatibility conditions (CC). Meanwhile, it will be rather striking to notice that the differential Galois theory may be applied  as it works in this framework because we have an algebraic Lie group of transformations, namely $\bar{y}=Ay+B$ where $A$ is an orthogonal $3\times 3$ matrix with $det(A)=1$ (see later on when such a condition is used) and $B$ is a vector. Of course, eliminating the parameters provides at once an {\it algebraic Lie pseudogroup} defined by the $m(m+1)/2$ differential polynomial equations in symbolic form: \\
\[     {\delta}_{uv} \frac{\partial {\bar{y}}^u}{\partial y^k} \frac{\partial {\bar{y}}^l}{\partial y^l}={\delta}_{kl}  \]
and all the second order jets vanish. The novelty of this presentation is that not only {\it all} known results of shell theory will appear for the first time as an effective application of a general theory with no reference to classical geometry and moving frames but other results will be obtained which are not known up to our knowledge, in particular when $m$ and $n$ are arbitrary.   \\

   First of all, we look for the corresponding differential invariants by examining the various prolongations ${\rho}_0(\theta), {\rho}_1(\theta), {\rho}_2(\theta), {\rho}_3(\theta), ...$ of the infinitesimal generators:\\
   \[ {\theta}_1=\frac{\partial}{\partial y^1}, {\theta}_2= \frac{\partial}{\partial y^2}, {\theta}_3=\frac{\partial}{\partial y^3},\]
   \[ {\theta}_4=y^2\frac{\partial}{\partial y^3}-y^3\frac{\partial}{\partial y^2}, {\theta}_5=y^3\frac{\partial}{\partial y^1}- y^1\frac{\partial}{\partial y^3}, {\theta}_6=y^1\frac{\partial}{\partial y^2}-y^2\frac{\partial}{\partial y^1}      \] 
as we need only stop when the symbol of the automorphic system is zero, according to the criteria for automorphic systems.  \\
Using the identity $y^1{\theta}_4+y^2{\theta}_5+y^3{\theta}^6\equiv 0$, the rank of the vertical distribution generated by the ${\rho}_0(\theta)$ on 
$V(X\times Y)=E$ is maximum and equal to $3$, so that there is no differential invariant of order zero and also any differential invariant cannot depend explicitly on $y$. The fiber of $V(J_1(X\times Y))\simeq J_1(V(X\times Y))=J_1(E)$ has dimension $3+6=9$ and there is therefore $9-6=3$ (rational) differential invariants of first order killed by ${\rho}_1(\theta)$, namely we get the Lie form at order $1$: \\
\[   \begin{array}{rcccccl}
{\Omega}_{11} & \equiv & {\sum}_k (y^k_1)^2& = & (y^1_1)^2+(y^2_1)^2+(y^3_1)^2 & = & 
{\omega}_{11}(x)\\
{\Omega}_{12} & \equiv & {\sum}_k y^k_1y^k_2& = & y^1_1y^1_2+y^2_1y^2_2+y^3_1y^3_2 & = & 
{\omega}_{12}(x)\\
{\Omega}_{22} & \equiv & {\sum}_k (y^k_2)^2& = & (y^1_2)^2+(y^2_2)^2+(y^3_2)^2 & = & 
{\omega}_{22}(x)
\end{array}   \]
or simpy:  \\
\[  {\Omega}_{ij}\equiv {\delta}_{kl}y^k_iy^l_j = {\omega}_{ij}(x)    \]
It is easy to check that $\omega=({\omega}_{ij}={\omega}_{ji}) \in S_2T^* $ and $\omega$ is called the {\it first fundamental form} though this name can be rather confusing. Then we know that each $d_i\Omega$ is also a differential invariant and we may set:  \\
\[{\Gamma}_{rij}=\frac{1}{2}(d_i{\Omega}_{rj}+d_j{\Omega}_{ri}-d_r{\Omega}_{ij}) \Rightarrow {\gamma}_{rij}=\frac{1}{2}({\partial}_i{\omega}_{rj}+{\partial}_j{\omega}_{ri}-{\partial}_r{\omega}_{ij})  \]
in order to get the $6$ linearly and thus functionally independent second order differential invariants:  \\
\[   {\Gamma}_{rij}\equiv {\vec{y}}_r.{\vec{y}}_{ij}={\delta}_{kl} y^k_ry^l_{ij} ={\gamma}_{rij}(x)={\gamma}_{rji}(x)  \]
In order to check this result, let us introduce the $2\times 3$ matrix $(y^k_i)$ of strict first order jets and consider the $3$ different $2\times 2$ subdeterminants like $y^1_1y^2_2-y^1_2y^2_1$. Then it is easy to obtain the relation:  \\
\[   det(\omega)\equiv {\omega}_{11}{\omega}_{22} - ({\omega}_{12})^2=
{\sum}_{det} (y^1_1y^2_2-y^1_2y^2_1)^2    \]
and we shall suppose from now on that $det(\omega)\neq 0$. Introducing the action and its prolongations:  \\
\[   {\bar{y}}^u=a^u_ky^k+b^u, {\bar{y}}^u_i=a^u_ky^k_i, {\bar{y}}^u_{ij}=a^u_ky^k_{ij}, ...  \]
we get:  \\
\[   {\bar{y}}^1_1{\bar{y}}^2_2-{\bar{y}}^2_1{\bar{y}}^1_2= ... + (-(a^1_1a^2_3-a^2_1a^1_3))(y^3_1y^1_2-y^1_1y^3_2)+ ...   \]
and obtain therefore with $A^t=transposed(A)$ and thus $A^{-1}=(1/det(A)) cof(A)^t$ by using the matrix of cofactors:   \\
\[      ( {\vec{\bar{y}}}_1 \wedge {\vec{\bar{y}}}_2)= det(A) (A^t)^{-1} ({\vec{y}}_1\wedge {\vec{y}}_2)     \]
But $AA^t=I\Rightarrow A^t=A^{-1}\Rightarrow (A^t)^{-1}=A$ and, if $det(A)=1$ ({\it care}), that is we deal with the connected component $G$ of the identity, then:  \\
\[  ( {\vec{\bar{y}}}_1 \wedge {\vec{\bar{y}}}_2).{\vec{\bar{y}}}_{ij}= ({\vec{y}}_1\wedge {\vec{y}}_2).{\vec{y}}_{ij}   \]  
Accordingly, we may introduce the unexpected $3$ {\it additional} differential invariants:  \\
\[ {\Sigma}_{ij}\equiv (y^2_1y^3_2-y^3_1y^2_2)y^1_{ij}+(y^3_1y^1_2-y^1_1y^3_2)y^2_{ij}+(y^1_1y^2_2-y^2_1y^1_2)y^3_{ij}={\sigma}_{ij}(x)  \]
and we notice that the $3\times 3$ matrix made by the factors of the second order jets $y^k_{ij}$ with $k=1,2,3$ and each $(ij)\in ((11), (12), (22))$ in $({\Sigma}_{ij},{\Gamma}_{1ij}, {\Gamma}_{2ij})$ has determinant equal to $det(\omega)$ and is thus of maximal rank $3$. It follows that the second order symbol defined by the system $(\Gamma=\gamma, \Sigma=\sigma)$ does vanish and the complete system $(\Omega=\omega, \Gamma=\gamma,\Sigma=\sigma)$ is not only automorphic but also involutive if and only if convenient {\it first order} generating compatibility conditions (CC) among $(\omega, \gamma, \sigma)$ are satisfied as a necessary and sufficient condition for applying the Cartan-K\"{a}hler theorem in the analytic case. The reader may compare the present approach with the one of M. Janet in ([13]) and with the one of P.G.Ciarlet in ([7,11]), always getting in mind that, in the formal theory of non-linear systems, the difficulty is to look for the various numbers of CC through techniques of acyclicity and diagram chasing which are in general very far from just doing "crossed derivatives " (Just consider the case $n=3, m= 6$ or the many examples provided in ([19-22,26])).\\

{\it We exhibit them for the first time within the formal framework of automorphic systems, still unknown today}.  \\

First of all, we get at once $\hspace{5mm}  {\partial}_r{\omega}_{ij}= {\gamma}_{irj}+{\gamma}_{jri} $. \\

As for the other CC, they are similar to the ones of a system like $y^k_{ij}=0$ with $k=1,2,3$ and $i,j=1,2$, by setting $y^1_1=1,y^2_2=1$ and the other first order jets equal to zeo, that is to say these CC are induced by the cokernel of the Spencer monomorphism $\delta: S_3T^* \otimes E\rightarrow T^*\otimes S_2T^*\otimes E$ for each $k=1,2,...,m$, that is $m(n^2(n+1)/2-n(n+1)(n+2)/6)=mn(n^2-1)/3$ which is equal to $6$ when $n=2,m=3$. Among these $6$ CC, we have the $2\times 2=4$ CC with symbols:  \\
\[   {\partial}_2{\gamma}_{i11} - {\partial}_1{\gamma}_{i12}, \hspace{3mm} {\partial}_2{\gamma}_{i12} - {\partial}_1{\gamma}_{i22}  \]
coming from the ordinary Riemann tensor and it just remains to compute the $2$ {\it new} CC having symbols:  \\
\[  {\partial}_2{\sigma}_{11} - {\partial}_1{\sigma}_{12}, \hspace{3mm} {\partial}_2{\sigma}_{12} - {\partial}_1{\sigma}_{22}   \]
in such a way that these $6$ quantities do not contain third order jets any longer. Indeed, thanks to differential algebra and the fact that the first prolongation of a nonlinear partial differential equation of order $q$ is quasi-linear in the jets of strict order $q+1$, we are thus led to a problem of pure linear algebra while studying cokernels.\\

We provide details for this {\it tricky} computation but we know from the differential Galois theory that these $6$ CC come from differential invariants of order $2$ and {\it must} therefore be expressed by means of rational functions of $(\omega, \gamma, \sigma)$, a fact not at all evident at first sight as we shall see but explaining why our definion of the second form slightly supersedes the standard one of textbooks as it {\it is the only one} fitting with differential algebra and differential Galois theory while using only rational differential functions in $\mathbb{Q}<y^1,y^2,y^3>=\mathbb{Q}<y>$. \\

Second, we get at once:  \\
\[d_2{\Gamma}_{111}-d_1{\Gamma}_{112}\equiv {\delta}_{kl}( y^k_{12}y^l_{11} -  y^k_{11}y^l_{12})=0 \Rightarrow {\partial}_2{\gamma}_{111} - {\partial}_2{\gamma}_{112}=0  \]
\[d_2{\Gamma}_{212}-d_1{\Gamma}_{222}\equiv {\delta}_{kl} (y^k_{22}y^l_{12} - y^k_{12}y^l_{22})=0 \Rightarrow {\partial}_2{\gamma}_{212} - {\partial}_2{\gamma}_{222}=0  \]
Similarly, we get:  \\
\[   {\partial}_2{\gamma}_{211}-{\partial}_1{\gamma}_{212}={\partial}_1{\gamma}_{122} - {\partial}_2{\gamma}_{112} =  {\vec{y}}_{11}.{\vec{y}}_{22} - {\mid} {\vec{y}}_{12}{\mid}^2  \]
and obtain therefore {\it only three} CC {\it for the only} $\gamma$ by adding {\it one more} CC, namely: \\
\[    {\partial}_2{\gamma}_{211}-{\partial}_1{\gamma}_{212}+{\partial}_2{\gamma}_{112} - {\partial}_1{\gamma}_{122} =0  \]
a result showing that {\it the other} CC {\it must contain} $\sigma$. \\

Meanwhile, we also get (See ([19], p 126-129]):  \\
\[  \frac{1}{2} ( {\partial}_{11}{\omega}_{22} + {\partial}_{22}{\omega}_{11}-2{\partial}_{12}{\omega}_{12})=
{\partial}_1{\gamma}_{212} - {\partial}_2{\gamma}_{211}= {\mid} {\vec{y}}_{12} {\mid}^2- {\vec{y}}_{11}.{\vec{y}}_{22}  \]
Using the formula ${\Sigma}_{ij}\equiv ({\vec{y}}_1\wedge {\vec{y}}_2).{\vec{y}}_{ij}={\sigma}_{ij} $ and various projections in the tangent plane
 to the surface like $({\vec{y}}_2\wedge ({\vec{y}}_1\wedge {\vec{y}}_2)).{\vec{y}}_1=det(\omega)$, we finally obtain:  \\
\[  \begin{array}{rcl}
det(\omega)y^k_{11} & = & {\sigma}_{11}({\vec{y}}_1\wedge {\vec{y}}_2)^k +{\gamma}_{111}
({\vec{y}}_2\wedge ({\vec{y}}_1\wedge {\vec{y}}_2))^k-{\gamma}_{211}({\vec{y}}_1\wedge ({\vec{y}}_1\wedge {\vec{y}}_2))^k  \\
det(\omega)y^k_{22} & = & {\sigma}_{22}({\vec{y}}_1\wedge {\vec{y}}_2)^k +{\gamma}_{122}
({\vec{y}}_2\wedge ({\vec{y}}_1\wedge {\vec{y}}_2))^k-{\gamma}_{222}({\vec{y}}_1\wedge ({\vec{y}}_1\wedge {\vec{y}}_2))^k  \\
det(\omega)y^k_{12} & = & {\sigma}_{12}({\vec{y}}_1\wedge {\vec{y}}_2)^k +{\gamma}_{112}
({\vec{y}}_2\wedge ({\vec{y}}_1\wedge {\vec{y}}_2))^k-{\gamma}_{212}({\vec{y}}_1\wedge ({\vec{y}}_1\wedge {\vec{y}}_2))^k  \\
\end{array}   \]
or simply:  \\
\[  det(\omega)\,{\vec{y}}_{ij}= {\sigma}_{ij} ({\vec{y}}_1\wedge {\vec{y}}_2) + 
{\gamma}_{1ij}({\vec{y}}_2\wedge ({\vec{y}}_1\wedge {\vec{y}}_2))-{\gamma}_{2ij}({\vec{y}}_1\wedge ({\vec{y}}_1\wedge {\vec{y}}_2))  \]
where we recall that:   \\
\[    det(\omega)=({\vec{y}}_2\wedge ({\vec{y}}_1 \wedge {Ê\vec{y}}_2)).{\vec{y}}_1= ({\vec{y}}_1,{\vec{y}}_2,    {\vec{y}}_1 \wedge {Ê\vec{y}}_2 )=
{\mid} {\vec{y}}_1 \wedge {Ê\vec{y}}_2 {\mid }^2   = {\mid }{\vec{y}}_1{\mid}^2{\mid}{\vec{y}}_2{\mid}^2 - {\mid}{\vec{y}}_1.{\vec{y}}_2{\mid}^2  \]
It follows that we have:  \\
\[  \begin{array}{lcl}
det(\omega)\,{\vec{y}}_{11}.{\vec{y}}_{22} & = &{\sigma}_{11}{\sigma}_{22} + {\omega}_{22}{\gamma}_{111}{\gamma}_{122} +{\omega}_{11}{\gamma}_{211}{\gamma}_{222}-{\omega}_{12}({\gamma}_{111}{\gamma}_{222} + {\gamma}_{211}{\gamma}_{122}) \\
 det(\omega)\,{\mid} {\vec{y}}_{12}{\mid}^2 & = & ({\sigma}_{12})^2 + {\omega}_{22}({\gamma}_{112})^2 +{\omega}_{11}({\gamma}_{212})^2- 2{\omega}_{12}{\gamma}_{112}{\gamma}_{212}  
 \end{array}  \]
Substracting the second equation from the first and setting $det(\sigma)={\sigma}_{11}{\sigma}_{22} - ({\sigma}_{12})^2$,  we get an equation of the form:  \\
\[  det(\omega)({\vec{y}}_{11}.{\vec{y}}_{22}-{\mid} {\vec{y}}_{12} {\mid}^2)=det (\sigma) + \omega \gamma\gamma  \]
where the last expression is a rational function of $j_1(\omega)\simeq (\omega,\gamma)$ according to the well known Levi-Civita isomorphism. This is {\it exactly} the " {\it theorema egregium} " of K.F. Gauss in ([10]) and, for this reason, $det(\sigma)/det(\omega)=(j_2(\omega))$ is called the {\it Gauss curvature} or {\it total curvature} and only depend on $\omega$. It is however important to notice that our definition of $\sigma$ does not involve the square root of $det(\omega)$. For this reason {\it it allows to use only rational differential functions} and must therefore be preferred to the standard one existing in the literature which is using the so-called {\it normal vector} $\vec{n}=({\vec{y}}_1\wedge {\vec{y}}_2)/\mid {\vec{y}}_1\wedge {\vec{y}}_2\mid $.\\
As another motivation for such a choice, we notice that $\omega \in S_2T^*$ is usually called {\it first fundamental form} while $\sigma$, which is called {\it second fundamental form}, is also considered as a section of $S_2T^*$ but this is not correct. Indeed, for any change $\bar{x}=\varphi (x)$ of independent variables, we have successively:  \\
\[  \frac{\partial \vec{y}}{\partial x^1}\wedge  \frac{\partial \vec{y}}{\partial x^2} =  (\frac{\partial \vec{y}}{\partial {\bar{x}}^1} \wedge \frac{\partial \vec{y}}{\partial {\bar{x}}^2} ) \frac{\partial ({\bar{x}}^1,{\bar{x}}^2)}{\partial (x^1,x^2)} \]
\[  \frac{{\partial}^2 \vec{y}}{\partial x^i\partial x^j} =  \frac{{\partial}^2 \vec{y}}{\partial {\bar{x}}^r\partial {\bar{x}}^s}    \frac{ \partial {\bar{x}}^r}{\partial x^i} \frac{\partial {\bar{x}}^s}{\partial x^j} + \frac{\partial \vec{y}}{ \partial {\bar{x}}^r} \frac{{\partial}^2 {\bar{x}}^r}{\partial x^i\partial x^j}         \]
with $i,j,r,s=1,2$ and we deduce that $\sigma \in S_2T^* \otimes {\wedge}^2T^*$ is in fact a {\it metric density}. \\

As for the last two CC for $\sigma$, an easy but tedious computation does provide:  \\
\[  \begin{array}{rcr}
det(\omega)({\partial}_2{\sigma}_{12} - {\partial}_1{\sigma}_{22}) & = & ({\gamma}_{122}\,{\omega}_{22} - {\gamma}_{222}\,{\omega}_{12}) \,{\sigma}_{11}   \\
                  &    & + (2{\gamma}_{112}\,{\omega}_{12}-2 {\gamma}_{212}\,{\omega}_{11}+{\gamma}_{211}\,{\omega}_{12}-
                  {\gamma}_{111}\,{\omega}_{22}) \,{\sigma}_{22}   \\
                  &    &  +(2{\gamma}_{212}\,{\omega}_{12}- 2{\gamma}_{112}\,{\omega}_{22}
                  +2{\gamma}_{222}\,{\omega}_{11}- 2{\gamma}_{122}\,{\omega}_{12}) \,{\sigma}_{12}
\end{array}  \]
and the other CC by exchanging $1$ and $2$. These equations have been first found  independently by D. Codazzi ([8]) and G. Mainardi ([18]) in a classical setting, with no reference to differential algebra. It follows that, contrary to $\omega$ which can be given arbitrarily provided that $det(\omega)\neq 0$, $\sigma$ {\it cannot be given arbitrarily}.  \\
 
 In certain cases, a known symmetry of the surface may be taken into account in order to determine explicitly $(\omega,\gamma, \sigma)$ and we provide two particular cases.  \\

$\bullet$  \hspace{1cm}  $y^1=x^1, y^2=x^2, y^3= \frac{1}{6}((x^1)^3 + (x^2)^3)$  \\
\[   {\omega}_{11}=1 + \frac{1}{4}(x^1)^4, {\omega}_{22}= 1 + \frac{1}{4}(x^2)^4, {\omega}_{12}= \frac{1}{4}(x^1x^2)^2  \Rightarrow det(\omega)= 1 +\frac{1}{4}(x^1)^4 + \frac{1}{4}(x^2)^4 \]
\[  {\gamma}_{111}=\frac{1}{2} (x^1)^3, {\gamma}_{112}=0, {\gamma}_{212}=0, {\gamma}_{222}=\frac{1}{2}(x^2)^3, {\gamma}_{211}=\frac{1}{2}x^1(x^2)^2, {\gamma}_{122}=\frac{1}{2}(x^1)^2x^2  \]
\[  {\sigma}_{11}= x^1, {\sigma}_{22}=x^2, {\sigma}_{12}=0  \Rightarrow det(\sigma)=x^1x^2\]
We check easily:  \\
\[ det(\omega)x^1x^2=x^1x^2 +\frac{1}{2}{\omega}_{22}(x^1)^5x^2+\frac{1}{2}{\omega}_{11}x^1(x^2)^5 - \frac{1}{2}{\omega}_{22}(x^1x^2)^3, {\partial}_2{\sigma}_{12}-{\partial}_1{\sigma}_{22}=0  \] 
and all the results must remain unchanged by the permutation $1 \leftrightarrow 2$.  \\

$\bullet$  \hspace{1cm}Again with $n=2,m=3$, let us consider the sphere of radius $R$ centered at the origin of the cartesian frame $(Oy^1y^2y^3)$ and apply a stereographic projection of the northern half sphere (draw a picture) with north pole $N$ on the equatorial plane $(Ox^1x^2)$ with coordinates $(x^1,x^2)$ from the south pole $S$. Using well known similarity of the triangles $SNM$ for a point $M$ on the sphere and $SXO$ where $X$ is the intersection of $SM$ with the equatorial plane, we get the classical formulas:  \\
\[  r^2=(x^1)^2+(x^2)^2, L^2=\mid {\vec{SX}}^2 \mid=R^2 + r^2, {\vec{SM}}^2=D^2=(y^1)^2+(y^2)^2+ (R+y^3)^2, \]
\[ LD=2R^2 \Rightarrow   \frac{y^1}{x^1}=\frac{y^2}{x^2}=1+\frac{y^3}{R}=\frac{2R^2}{L^2}, \frac{y^3}{R}=\frac{R^2-r^2}{L^2}\]
and obtain the parametrization of the sphere:  \\
\[  \frac{y^1}{R}=\frac{2Rx^1}{R^2+(x^1)^2+(x^2)^2},  \frac{y^2}{R}=\frac{2Rx^2}{R^2+(x^1)^2+(x^2)^2}, \frac{y^3}{R}=\frac{R^2-(x^1)^2-(x^2)^2}{R^2+(x^1)^2+(x^2)^2}  \]
leading to:  \\
\[  {\omega}_{11}={\omega}_{22}=4R^4/(R^2+(x^1)^2+(x^2)^2)^2={\phi}, \hspace{3mm}{\omega}_{12}={\omega}_{21}=0\Rightarrow det(\omega)={\phi}^2\]
We obtain at once: \\
\[  {\gamma}_{111}={\gamma}_{212}= - {\gamma}_{122}=\frac{1}{2}{\partial}_1\phi, \hspace{3mm} {\gamma}_{121}={\gamma}_{222}= - {\gamma}_{211}=\frac{1}{2}{\partial}_2\phi   \]
and, from the previous second order CC for $\omega$, we get:
\[  \frac{1}{2}{\phi}^2({\partial}_{11}\phi + {\partial}_{22}\phi)= - det(\sigma) + \frac{1}{2}\phi (({\partial}_1\phi)^2 + ({\partial}_2\phi)^2)  \]
An easy computation then leads to:   \\
\[ det(\sigma)=(16R^7/(R^2 + (x^1)^2 + (x^2)^2)^4)^2=\frac{1}{R^2}{\phi}^4   \]
Finally, though we have indeed:  \\
\[  {\sigma}_{ij}=det \left( \begin{array}{ccc}
                                        y^1_1 & y^1_2& y^1_{ij}  \\
                                        y^2_1 & y^2_2 & y^2_{ij} \\
                                        y^3_1 & y^3_2 & y^3_{ij}
                                        \end{array} \right)  \]
a direct computation is very tedious and we prefer to use the relations ${\omega}_{11}={\omega}_{22}=\phi=4R^4/L^4$ and ${\omega}_{12}=0$ in the formula:  \\
\[  {\phi}^2{\mid} {\vec{y}}_{12}{\mid }^2=({\sigma}_{12})^2+\phi (({\gamma}_{112})^2+({\gamma}_{212})^2)  \]
Substituting and using the relation $L{\partial}_iL=x^i$, we obtain:  \\
\[ y^1_{12}= - \frac{4R^2x^2}{L^4} + \frac{16R^2(x^1)^2x^2}{L^6},
   y^2_{12}= - \frac{4R^2x^1}{L^4} + \frac{16R^2x^1(x^2)^2}{L^6},
   y^3_{12}=\frac{16R^3x^1x^2}{L^6}  \]
and thus:  \\
\[ (y^1_{12})^2+(y^2_{12})^2+(y^3_{12})^2=\frac{16R^4((x^1)^2+(x^2)^2)}{L^8}   \]
\[  ({\gamma}_{112})^2+({\gamma}_{212})^2=\frac{1}{4}({\partial}_1\phi)^2+({\partial}_2\phi)^2)=\frac{64R^8((x^1)^2+(x^2)^2)}{L^{12}}  \]
that is to say ${\sigma}_{12}={\sigma}_{21}=0$. A similar computation left to the reader as an exercise gives ${\sigma}_{11}={\sigma}_{22}= - \frac{1}{R}{\phi}^2$ in agrement with the value of $det(\sigma)$ already obtained. As $\sigma$ is only determined up to the sign, we could also go to the north pole and notice that all the jets with $i=1,j=1$ vanish at $x^1=0,x^2=0$ but $y^1_1=y^2_2=2, y^3_{11}= - 4/R$ in order to obtain:\\
\[  {\sigma}_{11}(0,0)=det \left( \begin{array}{ccc}
                                        2 & 0& 0  \\
                                        0 & 2 & 0 \\
                                        0 & 0& -4/R
                                        \end{array} \right) = - 16/R     \]
as a simple way to know about the right sign.   \\

Let us now treat the general situation for arbirary $m$ and $n$ by using the formal theory of systems of partial differential eqquations and algebraic analysis ([19-24]). First of all, the case $n=2,m=3$ has been fully examined in ([19], p 126-129) and we recall it briefly in a more modern setting.\\

The initial system ${\cal{A}}_1\subset J_1(X\times Y)$ is, by construction, a PHS for the Lie groupoid ${\cal{R}}_1\subset {\Pi}_1(Y,Y)$ already defined as we have a groupoid action with:   \\  
\[     dim_X ({\cal{A}}_1)=dim_Y ({\cal{R}}_1)=6  \]
However, its first prolongation ${\cal{A}}_2={\rho}_1({\cal{A}}_1)\subset J_2(X\times Y)$ over $X$ is not at all a PHS for ${\cal{R}}_2={\rho}_1({\cal{R}}_1)$ over $Y$ because we have now 
$dim_X ({\cal{A}}_2)=9$ but $dim_Y ({\cal{R}}_2)=6$ and it follows that ${\cal{A}}_1$ is {\it not} an automorphic systems. Moreover, there is an additional difficulty coming from the fact that ${\cal{A}}_1$ is not even formally integrable because, introducing the second prolongation ${\cal{A}}_3={\rho}_2({\cal{A}}_1)$, there is a {\it strict inclusion} ${\cal{A}}^{(1)}_2 = {\pi}^3_2({\cal{A}}_3)\subset {\cal{A}}_2$ with $dim_X({\cal{A}}^{(1)}_2)=8$ as we have exhibited an additional second order equation and ${\cal{A}}_3$ is an affine bundle over ${\cal{A}}^{(1)}_2$ modelled on 
$g_3$. Dealing with nonlinear systems, we have the following commutative and exact diagram of affine bundles where the top row is made by the corresponding model vector bundles with $E=V(X\times Y)$:  \\
\[   \begin{array}{rcccccl}
0\longrightarrow & g_3 & \longrightarrow & S_3T^*\otimes E & \longrightarrow &S_2T^*\otimes S_2T^* &       \\
      &    \vdots  &    &   \vdots   &  & \vdots &    \\
      0  \longrightarrow & {\cal{A}}_3 &  \longrightarrow & J_3(X\times Y)) & \longrightarrow & J_2(S_2T^*) & \longrightarrow 0  \\
      & \downarrow &   &  \downarrow  & & \downarrow  &   \\
     0 \longrightarrow &  {\cal{A}}_2  & \longrightarrow  & J_2(X\times Y)  & \longrightarrow & J_1(S_2T^*)  & \longrightarrow 0
 \end{array}     \]
By chasing, the defect of surjectivity by $1$ of the arrow on the right of the top row is equal to the defect of surjectivity by $1$ of the arrow in the left column and $dim(g_3)= dim_X({\cal{A}}_3)-dim_X({\cal{A}}^{(1)}_2)=(30-18)-8=4$, a result not evident at all as it depends on a quite difficult {\it prolongation theorem} ([19,22,23):  \\

\noindent
{\bf THEOREM 2.1}: If $\cal{E}$ is a fibered manifold over $X$ and ${\cal{R}}_q\subset J_q({\cal{E}})$ is a system of order $q$ on $\cal{E}$, then, if the symbol $g_{q+1}$ is a vector bundle over ${\cal{R}}_q$ and $g_q$ is involutive or at least $2$-acyclic, one has the following non-trivial prolongation formula:  \\
\[  \begin{array}{rcl}
{\cal{R}}^{(1)}_{q+r}={\pi}^{q+r+1}_{q+r}({\cal{R}}_{q+r+1})={\pi}^{q+r+1}_{q+r}({\rho}_{r+1}({\cal{R}}_q))  & = &  {\rho}_r({\cal{R}}^{(1)}_q)   \\
          &  =  &  ({\cal{R}}^{(1)}_q)_{+r}\subset J_{q+r}({\cal{E}}) , \forall r\geq 0 
          \end{array}   \]

In the present situation, the symbol $g_2$ of ${\cal{A}}_2$ is defined by the $6$ linear equations:  \\
\[  {\delta}_{kl} y^k_rv^l_{ij}=0   \]
with $2$ equations of class $2$ and $4$ equations of class $1$ providing characters $(2,1)$ as we may always suppose that one determinant, say $y^1_1y^2_2-y^1_2y^2_1$ does not vanish. As the $4$ corresponding CC are easily seen to be satisfied, it follows that $g_2$ is involutive, a result leading to $dim(g_{3+r})=r+4$. \\
The symbol $g^{(1)}_2$ of ${\cal{A}}^{(1)}_2$ is defined by adding the only equation:  \\
\[  {\delta}_{kl}(y^k_{11}v^l_{22}+ y^k_{22}v^l_{11}-2y^k_{12}v^k_{12})=0  \]
Meanwhile, $3$ among the $7$ equations can be solved with respect to $(y^1_{22},y^2_{22},y^3_{22})$ provided that the determinant ${\sigma}_{11}$ does not vanish. This additional equation reduces the second character to zero and, for the same reason as above because the equations of class $1$ are untouched. It follows that $dim(g^{(1)}_2)=2$ and $g^{(1)}_2$ is thus also involutive with characters $(2,0)$ providing therefore $dim(g^{(1)}_{2+r})=2$. We obtain therefore ${\cal{A}}^{(1)}_{2+r}=({\cal{A}}^{(1)}_2)_{+r}$ with fiber dimension equal to:  \\
\[   dim_X({\cal{A}}^{(1)}_{2+r})= dim(g^{(1)}_3) + ... +dim(g^{(1)}_{2+r}) +dim_X ({\cal{A}}^{(1)}_2)= 2r + 8  \]
and ${\cal{A}}^{(1)}_2$ is an involutive system. It follows that: \\
\[   dim_X({\cal{A}}_{3+r})= dim g_{3+r} + dim_X({\cal{A}}_{2+r})=(r+4) + (2r+8)= 3r + 12  \]
while   \\
\[  dim_X(J_{3+r}(X\times Y))=3(r^2+9r+20)/2,   \,  dim(J_{2+r}(S_2T^*))=3(r^2+7r+12)/2   \]
and thus:   \\
\[ dim({\cal{A}}_{3+r})=dim_X(J_{3+r}(X\times Y)) - dim_X(J_{2+r}(S_2T^*))  \]
It follows that $\omega$ can be given arbitrarily and we have the following commutative diagram of affine bundles:  \\
\[   \begin{array}{rcccccccl}
0 \rightarrow  & g_{3+r} & \longrightarrow & S_{3+r}T^*\otimes E & \longrightarrow &S_{2+r}T^*\otimes S_2T^* &  \longrightarrow   & ? &\rightarrow 0     \\
      &    \vdots  &    &   \vdots   &  & \vdots &   & & \\
      0  \rightarrow & {\cal{A}}_{3+r} &  \longrightarrow & J_{3+r}(X\times Y) & \longrightarrow & J_{2+r}(S_2T^*) & \longrightarrow  & 0 &\\
      & \downarrow &   &  \downarrow  & & \downarrow  &  & &  \\
0 \rightarrow &  {\cal{A}}_{2+r}  & \longrightarrow  & J_2(X\times Y)  & \longrightarrow &  J_{1+r}(S_2T^*) & \longrightarrow  &0 &
 \end{array}     \]
where the defect of surjectivity of the central upper arrow in the first row is equal to the defect of surjectivity of the left downarrow that is exactly:  \\
\[  dim_X({\cal{A}}_{2+r}) - dim_X({\cal{A}}^{(1)}_{2+r}= (3r+9) - (2r+8)=r+1\]
coming from the fact that ${\cal{A}}_2$ is {\it not} formally integrable. We have thus provided a modern proof of the fact that it is always possible to embed a riemannian surface into an euclidean space when $n=2,m=3$, a beautiful and tricky result first obtained by M. Janet in 1926 ([13]) and by E. Cartan in 1927 ([5]).  \\

Changing slightly the notations while adding the three additional Lie equations ${\Sigma}_{ij}={\sigma}_{ij}$, we obtain a PHS ${\cal{A}}_2$ for ${\cal{R}}_2$, both systems having a zero symbol and we have the successive inclusions ${\cal{A}}_2 \subset {\pi}^3_2({\rho}_1({\cal{A}}_1)) \subset {\rho}_1({\cal{A}}_1)$ with respective fiber dimensions $ 6 < 8 < 9$. According to the fundamental theorems on automorphic systems, it follows that {\it only} ${\cal{A}}_2$ {\it is indeed an automorphic system} provided that convenient CC are satisfied, being described by the top row of the new following commutative and exact diagram of affine bundles with symbolic notations, that must be compared to the previous ones:  \\
\[   \begin{array}{rcccccccl}
       & 0 & \longrightarrow & S_3T^*\otimes E & \longrightarrow &T^*\otimes (\Omega,\Gamma,\Sigma) &  \longrightarrow   & CC &\rightarrow 0     \\
      &    \vdots  &    &   \vdots   &  & \vdots &   & & \\
      0  \rightarrow & {\cal{A}}_3 &  \longrightarrow & J_3(X\times Y)) & \longrightarrow & J_1(\Omega,\Gamma,\Sigma) &  &  &\\
      & \downarrow &   &  \downarrow  & & \downarrow  &  & &  \\
     0 \rightarrow &  {\cal{A}}_2  & \longrightarrow  & J_2(X\times Y)  & \longrightarrow &  (\Omega,\Gamma,\Sigma) & \longrightarrow  &0 &
 \end{array}     \]
The number of desired CC is determined by the exactness of the top row and is: \\
\[Ê \begin{array}{rcl}
nb(CC)  & =   & dim( T^*\otimes (\Omega,\Gamma,\Sigma)) - dim (S_3T^*\otimes E)\\
             & =  &2 \times dim(\Omega,\Gamma,\Sigma) - 4 \times dim (E)   \\
             & = &  2\times (3+6+3) - 4\times 3= 24-12= 12
             \end{array} \]
             
Among these $12$ CC, we have $6$ CC of the form $\partial \hspace{1mm} \omega=\gamma + \gamma$, then $2$ CC of the form 
$\partial \hspace{1mm}\gamma - \partial \hspace{1mm} \gamma =0$, then $1$ CC of the form $\partial \hspace{1mm}\gamma - \partial \hspace{1mm} \gamma + \partial \hspace{1mm}\gamma - \partial \hspace{1mm} \gamma =0$, then $1$ CC of the form:  \\
\[  det(\omega)({\partial}_2{\gamma}_{211} - {\partial}_1{\gamma}_{212}) - \omega\gamma\gamma = det(\sigma)   \]
providing Gauss theorem and finally $2$ additional CC for $\sigma$ providing the Codazzi-Mainardi equations. \\
{\it It is therefore only now that we do really understand the structure of shell theory and its relation with the theory of automorphic systems}.\\

The generalization of the previous results to arbitrary dimensions is of course quite more difficult and could provide a lot of work for the future. For simplicity, we shall restrict our study to the local isometric embedding problem of a riemann surface of dimension $n$ in an euclidean space of 
dimension $m$. Though surprising it may look like at first sight, our study will highly depend on the following result involving the differential rank of differential modules and {\it double duality} (Compare to a modern proof of the Janet conjecture in ([24],[23], p 539):  \\

\noindent
{\bf THEOREM 2.2}: If $K$ is a differential field and a left differential module $M={ }_DM$ is defined over the ring $D=K[d_1, ... ,d_n]=K[d]$ of differential operators with coefficients in $K$ by a finite presentation $D^p \stackrel{{\cal{D}}}{\longrightarrow} D^m \longrightarrow M \rightarrow 0 $, then ${\cal{D}}$ is formally injective if and only if $rk_D(M)=m-p$. In particular, if $m=p$, then ${\cal{D}}$ is injective if and only if $M$ is a torsion module, that is $rk_D(M)=0$.  \\

\noindent
{\it Proof}: First of all, we recall that $rk_D(M)$ is the dimension of the biggest free differential module that can be contained in $M$ and that the differential rank satisfies the additivity condition $ rk_D(M)=rk_D(M') + rk _D(M")$ for any short exact sequence $0 \rightarrow M' \rightarrow M \rightarrow M" \rightarrow 0$ ([23,24]). A similar property is existing in the non-linear framework for differential extensions $K \subset L \subset M$ where $diff\,rk (M/K)=diff\,rk (L/K) +diff\,rk(M/L)$ and generalizes the well known classical purely algebraic situation ([20,22,25]). As a second comment, this result is valid even if $D$ is non-commutative and if ${\cal{D}}$ is not formally integrable and this will be the case of the situations we shall study. Applying $hom_D(\bullet,D)$ to the presentation, we may define the "{\it right} " (care) differential module $N_D$ by the long exact sequence $0 \leftarrow N_D \longleftarrow D^p \longleftarrow D^m \longleftarrow hom_D(M,D) \leftarrow 0$ of {\it right} differential modules where we use the bimodule structure of $D={ }_DD_D$ where the left action of $P\in D$ on $D$ is defined by $Q \rightarrow PQ$ while the right action is defined by $Q \rightarrow QP$ for any $Q\in D$. Meanwhile, the right action of $P\in D$ on $f \in hom_D(M,D)$ is defined by $ (fP)(m)=f(m)P$ but a left action cannot be defined and we check that:\\
\[ ((fP)Q)(m)=((fP)(m))Q=(f(m)P)Q=f(m)PQ=(f(PQ))(m), \forall m\in M\]
in such a way that:  \\
 \[(fP)(Qm)=f(Qm)P=Qf(m)P=Q(f(m)P)=Q(fP)(m) , \forall P,Q\in D \] 
Of course, when $D$ is commutative, both actions coincide and $hom_D(D,D)\simeq D$ because $1\in D$ and thus 
$f(Q)=f(Q.1)=Q(f(1))=QP$ if we set $f(1)=P$.  
We may also pass from the right differential module $N_D$ to a left differential module $N={ }_DN=hom_K({\wedge}^nT^*,N_D)$ by the {\it side changing procedure} if we notice that $D$ is generated by $K$ and $T=\{a^id_i\}\subset D$ with $T^*=hom_K(T,K)$. In this case, the {\it dual operator} $hom_D({\cal{D}},D)$ becomes the {\it adjoint operator} $ad({\cal{D}})$ used in variational calculus.\\

Finally, it can be proved that the additive property of the differential rank is such that $rk_D(hom_D(M,D))=rk_D(M)$ and $rk_D(ad({\cal{D}}))=rk_D({\cal{D}})$ if we set $rk_D({\cal{D}})=rk_D(im({\cal{D}}))$. Counting the differential ranks, we get:   \\
\[ rk_D(N)-p+m-rk_D(M)=0  \]
and thus discover that $rk_D(N)=0$ if and only if $rk_D(M)=m-p$. As $D$ is an integral domain and $N$ is a torsion module, then $hom_D(N,D)=0$ and we have obtained by biduality the short exact sequence $0 \rightarrow D^p \stackrel{{\cal{D}}}{\longrightarrow } D^m \rightarrow M \rightarrow 0$. For the reader not familiar with homological algebra, if $L=ker({\cal{D}})$ in the long exact sequence $0 \rightarrow L \rightarrow D^p \stackrel{{\cal{D}}}{\rightarrow} D^m \rightarrow M \rightarrow 0 $, then $rk_D(L)=0$ and thus $L$ is a torsion module over the differential integral domain $D$. Hence, any element $x\in L$ is such that there exists {\it at least one} $0\neq P\in D$ such that $Px=0$ and thus $x=0, \forall x\in L$ because $D$ and thus $D^p$ do not admit divisors of zero., that is $L=0$. We can also say that $L$ is a torsion module and obtain a contradiction with the fact that $L$ is contained in the free and thus torsion-free module $D^p$, unless $L=0$. \\
\hspace*{12cm}     Q.E.D.  \\

Using the notations of M. Janet in (13]), we have the non-linear system ${\cal{A}}_1$:  \\
\[   {\Omega}_{i'j'}\equiv {\delta}_{kl}y^k_{i'}y^l_{j'}={\omega}_{i'j'}, \hspace{3mm}
     {\Omega}_{i'n}\equiv {\delta}_{kl}y^k_{i'}y^l_{n}={\omega}_{i'n}, \hspace{3mm}
     {\Omega}_{nn}\equiv {\delta}_{kl}y^k_{n}y^l_{n}={\omega}_{nn}   \]
where $i',j'=1,...,n-1$ and $k,l=1,...,m$. We may, as before, differentiate all these equations once in order to obtain the various 
${\Gamma}_{rij}={\gamma}_{rij}$ with now $i,j,r=1,...,n$, in particular the $n$ differential invariants:  \\
\[  {\Gamma}_{rnn}\equiv {\delta}_{kl}(y^k_{r}y^l_{nn} )={\gamma}_{rnn}, \hspace{4mm} r=1,...,n  \]
and consider the $n(n-1)/2$ new differential invariants of order two (Look for example in ({AIRY}) for $m=n=3$):  \\
\[    d_{nn}{\Omega}_{i'j'} + d_{i'j'}{\Omega}_{nn} - d_{i'n}{\Omega}_{j'n} - d_{j'n}{\Omega}_{i'n}\equiv 
2{\delta}_{kl}(y^k_{i'n}y^l_{j'n}- y^k_{i'j'}y^l_{nn})  \]
Accordingly, the possibility to compute the character of index $n$ of the system ${\cal{A}}^{(1)}_2\subset {\rho}_1({\cal{A}}_1)$ only depends on 
the matrix with $m$ rows and $n+(n(n-1)/2)=n(n+1)/2$ columns needed for solving the equations with respect to the jets $y^1_{nn},...,y^m_{nn}$:\\
\[ (y^1_{nn} ... y^m_{nn}) \left( \begin{array}{cccc}
                                        y^1_1 & ... & y^1_n &  y^1_{i'j'}  \\
                                        ... & ... & ... & ...                 \\
                                        y^m_1& ... & y^m_n & y^m_{i'j'}
                                        \end{array} \right)     \]
We find the results previously obtained for $n=2,m=3$ where the corresponding $3\times 3$ square matrix has a determinant equal to 
${\sigma}_{11}$. Hence, when $m=n(n+1)/2=p$ and generic jets, the rank of the above system is equal to $m$ and {\it the character} $n$ {\it of the symbol must vanish}. It follows from the previous theorem that we cannot have CC for the $\omega$ and we obtain a new proof in a modern setting of the result first found by M. Janet in $1926$.  \\

Coming back to the case $n=2,m=3$, we may introduce the matrix :  \\
\[   M(x)=({\vec{y}}_1,{\vec{y}}_2,{\vec{y}}_1\wedge{\vec{y}}_2)  \]
This matrix is invertible because $det(M)=\mid {\vec{y}}_1\wedge{\vec{y}}_2 {\mid}^2 = det(\omega)\neq 0$. Also, looking at the action of the group of rigid motions, we have:  \\
\[  \bar{y}=Ay+B \Rightarrow {\vec{y}}_i=A{\vec{y}}_i \Rightarrow \bar{M}=AM\Rightarrow A=\bar{M}M^{-1}, B=\bar{y}-Ay  \]
Accordingly, if $y_2=f_2(x)\neq j_2(f)(x)$ and ${\bar{y}}_2={\bar{f}}_2(x)\neq j_2(\bar{f})(x)$ are two sections of $J_2(X\times Y)$, then $M=M(x), \bar{M}=\bar{M}(x)$ leads to $A=A(x), B=B(x)$ with:\\
\[ \bar{f}(x)=A(x)f(x)+B(x), {\bar{f}}_x(x)=A(x)f_x(x), {\bar{f}}_{xx}(x)=A(x)f_{xx}(x)  \]
and the relations:  \\
\[  {\partial}_x\bar{f}- {\bar{f}}_x=A({\partial}_xf-f_x)+{\partial}_xAf+{\partial}_xB,\hspace{3mm}
{\partial}_x\bar{f}_x - {\bar{f}}_{xx}=A({\partial}_xf_x-f_{xx})+{\partial}_xAf_x   \]
If ${\bar{f}}_2\neq j_(\bar{f})$ but $f_2=j_2(f)$ and thus $f$ is one solution of the automorphic system, then we have:  \\
\[   {\partial}_x\bar{f}- {\bar{f}}_x=({\partial}_xA{A}^{-1})\bar{f}+({\partial}_xB-{\partial}_xAA^{-1}B),\hspace{3mm}
{\partial}_x\bar{f}_x - {\bar{f}}_{xx}=({\partial}_xAA^{-1}){\bar{f}}_x   \]
bringing right invariant $1$-forms which are the pull-back of the Maurer-Cartan forms by the gauging procedure $X\rightarrow G:(x)\rightarrow (A(x),B(x))$ induced by the two sections $f_2,{\bar{f}}_2\in {\cal{A}}_2$. A similar left invariant version is obtained by composing the gauging with the inversion map $G\rightarrow G:(A,B) \rightarrow (A^{-1}, - A^{-1}B)$. Sections of jet bundles provide therefore a modern version of the so-called Darboux vectors but this is out of the scope of this paper.  \\

\newpage

\noindent
{\bf 3) CHAIN THEORY REVISITED}  \\

The shape $y=f(x)$ of a thin {\it inelastic} free hanging chain ( a bike chain for example ) in a vertical plan with cartesian coordinates ($x$ horizontal,$y$ vertical) under gravity $g$ is usually called {\it chainette}, {\it catenary} or sometimes {\it funicular}. If $\vec{T}(x)$ is the tension in the chain and $s$ the curvilinear abcissa along the chain, we have $ds=\sqrt{1+({\partial}_xf)^2}$ and the equilibrium condition of a piece of chain with length $ds$ and mass $m$ per unit of length is:\\
\[  \frac{d\vec{T}}{ds} + m\vec{g}=0  \]
As $\vec{T}=\lambda (1,{\partial}_xf)$, we get ${\partial}_x\lambda=0 \Rightarrow \lambda =cst$ and thus ${\partial}_{xx}f/\sqrt{1+({\partial}_xf)^2}=cst=1/a$. One possible solution is $y=a\, ch(x/a)$ and we shall only consider the particular case $y=ch(x)$ obtained when $a=1$. \\
As before, the problem becomes:   \\

\hspace*{5mm}              {\it HOW TO DETERMINE A CURVE IN\, ${\mathbb{R}}^2$ UP TO A RIGID MOTION}\\

Changing slightly the notations with $n=1,m=2$, we may now consider the algebraic group $G$ of rigid motions of the plan $Y={\mathbb{R}}^2$ with cartesian coordinates $y=(y^1,y^2)$ while $X=\mathbb{R}$ has coordinate $(x)$. The finite action is defined by $\bar{y}=Ay+B$ with $AA^t=I\Rightarrow (det(A))^2=1$ and we may also consider, as before, the connected component $G^0\subset G$ of the identity with $det(A)=1$. The three infinitesimal generators of the group action are:\\
\[   \{  {\theta}_1=\frac{\partial}{\partial y^1},{\theta}_2=\frac{\partial}{\partial y^2}, 
{\theta}_3=y^1\frac{\partial}{\partial y^2} - y^2 \frac{\partial}{\partial y^1}  \}  \]
and the prolongation to the jets of order $2$ only depends on ${\theta}_1, {\theta}_2$ and:  \\
\[  {\rho}_2( {\theta}_3)=y^1\frac{\partial}{\partial y^2} - y^2 \frac{\partial}{\partial y^1}+y^1_x\frac{\partial}{\partial y^2_x} - y^2_x \frac{\partial}{\partial y^1_x} + y^1_{xx}\frac{\partial}{\partial y^2_{xx}} - y^2_{xx} \frac{\partial}{\partial y^1_{xx}}     \]
The corresponding PHS ${\cal{A}}_1$ is easily seen to be ${\Omega}^1\equiv (y^1_x)^2 + (y^2_x)^2=\omega (x)$ but is {\it not} an automorphic system. Exactly like in the study of shell, we may therefore consider the automorphic system ${\cal{A}}_2 \subset {\rho}_1({\cal{A}}_1)$:  \\
\[  {\Omega}^1\equiv (y^1_x)^2 + (y^2_x)^2=\omega (x), \hspace{3mm} \Gamma \equiv y^1_xy^1_{xx} + y^2_xy^2_{xx}=\gamma(x), 
\hspace{3mm} \Sigma \equiv y^1_xy^2_{xx}- y^2_xy^1_{xx}=\sigma (x)  \]
The symbol of order $2$ vanishes if and only if $\omega \neq 0$ and we shall assume such a condition. The only CC is $\gamma = \frac{1}{2}{\partial}_x\omega $ but, contrary to shell theory, $\sigma$ is arbitrary. However, our definitions allow to avoid using square roots in the geometrical approach. Indeed, setting as usual $\frac{d\vec{y} }{ds} = \vec{t}, \frac{d\vec{t}}{ds}= \kappa \vec{n}, \frac{d\vec{n}}{ds}= - \kappa \vec{t} $, we have $ds^2=\omega dx^2$ and thus $\mid  \vec{t}\,{\mid}^2= 1 \Rightarrow \vec{t}.\vec{n}=0 \Rightarrow \vec{t} \wedge \vec{n}=\vec{b}\Rightarrow \mid \vec{b}\,{\mid}^2=1$ with $\vec{b}$ fixed perpendicular to the plane considered. It follows that:  \\
\[  \frac{d\vec{y}}{dx}\wedge \frac{d^2\vec{y}}{dx^2}= (\frac{ds}{dx}\vec{t}\,)\wedge (\frac{d}{dx}(\frac{ds}{dx}\vec{t}\,))=
  \omega   \vec{t}\wedge  \frac{d\vec{t}}{dx} = \omega \vec{t}\wedge (\frac{ds}{dx}\frac{d\vec{t}}{ds})=\omega \frac{ds}{dx}\kappa  \vec{b}    \]
and thus $\sigma^2={\omega}^3 {\kappa}^2  $. It is however important to notice that $\bar{y}=Ay+B \Rightarrow {\bar{y}}_x=Ay_x, {\bar{y}}_{xx}=Ay_{xx}, ...$ and thus $\bar{\sigma}=
det(A)\sigma$, that is $\Sigma$ is only invariant by $G^0$ while the new Lie equation $\Upsilon\equiv (y^1_{xx})^2 + (y^2_{xx})^2=\upsilon(x)$ is invariant by $G$. Nevertheless, an elementary computation provides the idenity $\Omega \Upsilon=(\Gamma)^2 + (\Sigma)^2 \Rightarrow \omega\upsilon={\gamma}^2 + {\sigma}^2$ and we have the successive inclusions of differential extensions $k\subset K \subset K_0 \subset L$:  \\
\[    \mathbb{Q} < \mathbb{Q}<\Omega, \Upsilon>\subset \mathbb{Q}<\Omega, \Sigma >\subset \mathbb{Q} < y >   \]
because $\Gamma \in \mathbb{Q}<\Omega>$. Then $K_0=\mathbb{Q}<\Omega,\Sigma>$ is thus purely algebraic over $K=\mathbb{Q}<\Omega,\Upsilon>$ in agrement with the fact that $G/G^0=\{1, - 1\}$ while $L=\mathbb{Q}<y>$ is a regular extension of $K_0=\mathbb{Q}<\Omega,\Sigma>$ because $G$ is defined over $k=\mathbb{Q}$ and $k(G_0)$ is a regular extension of $k$, that is $k$ is algebraically closed in $k(G_0)$. As a byproduct, we obtain now:\\
\[   \upsilon = \mid\frac{d^2\vec{y}}{d^2x}{\mid}^2=\mid (\frac{d^2s}{dx^2}\vec{t}+ (\frac{ds}{dx})^2\frac{d\vec{t}}{ds}{\mid}^2=
(\frac{d^2s}{dx^2})^2 + (\frac{ds}{dx})^4{\kappa}^2=({\gamma}^2/\omega ) + {\omega}^2 {\kappa}^2  \]
because $\frac{d^2s}{dx^2}= \gamma / \sqrt{\omega}$ and, multiplying by $\omega$, we find back the previous formulas. It follows that 
$\omega=1 \Rightarrow \gamma=0 \Rightarrow \sigma=\kappa$ in the case of an inelastic chain.\\
Let us now consider two sections of ${\cal{A}}_2$, namely $f_2=j_2(f)$ with $y^1=x, y^2=ch(x)$ giving $ \omega= ch^2(x), \gamma= sh(x)ch(x), \sigma=ch(x), \upsilon= {ch}^2(x)$ and ${\bar{f}}_2\neq j_2(\bar{f})$ defined by $ {\bar{y}}^1=sh(x), {\bar{y}}^2= 1, {\bar{y}}^1_x=ch(x), {\bar{y}}^2_x= 0, {\bar{y}}^1_{xx}= sh(x), {\bar{y}}^2_{xx}=1$ in order to obtain right invariant Maurer-Cartan forms. Indeed, the gauging is defined by the formulas:  \\
\[ \bar{y}=Ay+B\Rightarrow \bar{f}(x)=A(x)f(x)+B(x), \bar{f}_x(x)=A(x)f_x(x), \bar{f}_{xx}(x)A(x)f_{xx}\]
that is to say:  \\
\[  \forall f_2,\bar{f}_2 \in {\cal{A}}_2 \Rightarrow \exists \, gauging \,(A(x),B(x)) \]
An easy computation provides:  \\
\[ A(x)=\left( \begin{array}{cc}
\hspace{2mm} \frac{1}{ch(x)} & \frac{sh(x)}{ch(x)}\\
 - \frac{sh(x)}{ch(x)} & \frac{1}{ch(x)}
 \end{array}  \right )  , \hspace{5mm} B(x)= \left(  \begin{array}{c}
 - \frac{x}{ch(x)}  \\  \hspace{3mm}\frac{x\,sh(x)}{ch(x)}
 \end{array}  \right)  \]
where of course $A$ is an orthogonal $2\times 2$ matrix with $det(A)=1$. Equivalently, we may introduce the $2\times 2$ matrices $M(x)$ and 
$\bar{M}(x)$ defined by:  \\
\[ M=\left(  \begin{array}{cc}
f^1_x & f^1_{xx}  \\   f^2_x & f^2_{xx}
\end{array} \right) ,\, \bar{M}=\left(  \begin{array}{cc}
{\bar{f}}^1_x & {\bar{f}}^1_{xx}  \\   {\bar{f}}^2_x & {\bar{f}}^2_{xx}
\end{array} \right)  \Rightarrow  \bar{M}(x)=A(x)M(x)   \] 
and the inversion formula $A=\bar{M}M^{-1}$ can be used whenever $\bar{\sigma}=det(\bar{M})=det(M)=\sigma \neq 0 \Rightarrow det(A)=1$, leading to $B=\bar{f} -Af$.The components of the Spencer operator are now $Df_2=0, \, D{\bar{f}}_2$ with:  \\
\[ \frac{\partial \bar{f}}{\partial x} - {\bar{f}}_x= (\frac{\partial A}{\partial x}A^{-1})\bar{f}+(\frac{\partial B}{\partial x} - \frac{\partial A}{\partial x}A^{-1}B), \,  \frac{\partial {\bar{f}}_x}{\partial x} - {\bar{f}}_{xx}= (\frac{\partial A}{\partial x}A^{-1}){\bar{f}}_x  \]
We finally obtain:  \\
\[  \frac{\partial A}{\partial x}A^{-1}= \left( \begin{array}{cc} 0 & \frac{1}{ch(x)} \\ - \frac{1}{ch(x)} & 0  \end{array} \right), \,
\frac{\partial B}{\partial x} - \frac{\partial A}{\partial x}A^{-1}B= \left(  \begin{array}{c} - \frac{1}{ch(x)} \\  
\frac{sh(x)}{ch(x)}\end{array} \right) \]
where the first matrix is skew-symmetric as expected but there is no Maurer-Cartan equation because $dim(X)=1$.   \\
Conversely, if these $1$-forms are known and provide $(A,B)$, then any other solution $(\bar{A},\bar{B})$ is such that:  \\
\[   \frac{\partial \bar{A}}{\partial x}{\bar{A}}^{-1}=  \frac{\partial A}{\partial x}A^{-1}, \hspace{5mm}
     \frac{\partial \bar{B}}{\partial x} - \frac{\partial \bar{A}}{\partial x}{\bar{A}}^{-1}\bar{B}  = \frac{\partial B}{\partial x} - \frac{\partial A}{\partial x}A^{-1}B   \]
ad thus:  \\
\[  \frac{\partial}{\partial x}(A^{-1}\bar{A})=0  \Rightarrow \bar{A}=AC, C=cst \Rightarrow \bar{B}=AD+B, D=cst  \]
This is a modern version of the so-called Darboux vectors which are replaced by sections of jet bundles.  \\

Coming back to the fundamental theorem of the differential Galois theory ([3,4,20,22]), the reciprocal distribution $\Delta$ commuting with the invariant distribution $\Theta$ generated by the ${\rho}_2({\theta}_{\tau}), \, \tau=1,2,3$, is generated by the $4=dim(J^0_2(E))=dim(J_2(E)-m$: \\
\[   {\delta}_1= y^1_x\frac{\partial}{\partial y^1_x}+ y^2_x\frac{\partial}{\partial y^2_x} , 
\hspace{4mm} {\delta}_2= y^1_{xx}\frac{\partial}{\partial y^1_x}+y^2_{xx}\frac{\partial}{\partial y^2_x}, 
\hspace{4mm}{\delta}_3= y^1_x\frac{\partial}{\partial y^1_{xx}}+y^2_x\frac{\partial}{\partial y^2_{xx}}, \]
\[  {\delta}_4= y^1_{xx}\frac{\partial}{\partial y^1_{xx}}+y^2_{xx}\frac{\partial}{\partial y^2_{xx}}   \]
First of all, it is easy to check that all the differential extensions already considered are stable by $\Delta$ and that the distribution $\delta \otimes \delta$ or $\delta + \bar{\delta}$ where $\bar{\delta}$ is obtained from any $\delta \in \Delta$ by adding a bar over the jet coordinates, satisfies:  \\
\[ \bar{\delta}.\bar{M}=A\delta.M, \, \bar{M}=AM , \, \bar{y}=Ay+B  \Rightarrow (\delta+\bar{\delta}).A=0, \, (\delta+\bar{\delta}).B=0  \]
and thus $(A,B)$ is killed by $ (\delta + \bar{\delta}), \forall \delta \in \Delta$. Moreover, if any vector field $\delta$ commutes with ${\delta}_1,...,{\delta}_4$, then, applying any $\theta \in \Theta$ to the linear combination $\delta= a_1{\delta}_1 + ... +a_4{\delta}_4$ we discover that the $a$ are kiled by $\theta$ and are thus only functions of the differential invariants. For example, if we consider $\delta=y^1_x\frac{\partial}{y^2_x} - y^2_x \frac{\partial}{\partial y^1_x}$ we obtain $\delta = - (\Gamma/\Sigma){\delta}_1 + (\Omega / \Sigma){\delta}_2$ a result not evident at all at first sight. We obtain therefore $\delta \Gamma= \Sigma$ and thus $\delta K\not\subset K$ because $\Sigma \notin K$ but $\delta K_0\subset K_0$.

We end this example with a few comments on variational calculus with constraints presented at the end of ([23]). To start with, we shall even suppose that the chain $C$ is elastic while hanging between two fixed points, with initial mass $m_0(x)$ per unit lentgh $x$, final mass $m(s)$ per unit legth $s$ in such a way that $m_0dx=m(s)ds$ and elastic coefficient $E$ in such a way that, under a tension $T$ we have $T=E(ds-dx)/dx=E(ds/dx)-1)$. Using previous notations, the total lagrangian including the gravitational potential and the energy of deformation stored in the chain, we have to study the following variational problem where $y^k=f^k(x), y^k_x=(df^k_x)/dx\ \forall k=1,2$ :  \\
\[  \delta {\int}_C (m_0gy^2 + \frac{1}{2}E(\frac{ds}{dx} - 1)^2 )dx={\int}_C (m_0g\delta y^2 + T (y^1_x\delta y^1_x+y^2_x\delta y^2_x)/(\frac{ds}{dx}))dx=0  \]
Integrating by part while using the relations $\frac{d y^1}{dx}/\frac{ds}{dx}= cos(\theta),\frac{d y^2}{dx}/\frac{ds}{dx}= sin(\theta) $, we obtain successively the two Euler-Lagrange equations:   \\
\[ \begin{array}{lcl}
\delta y^1 \rightarrow  \frac{d}{dx}(Tcos(\theta))=0    & \Rightarrow & Tcos(\theta)=T_0  \\
\delta y^2 \rightarrow \frac{d}{dx}(Tsin(\theta))=m_0g  &  \Rightarrow &  T sin(\theta)=m_0gx  
\end{array}   \]
Setting $T_0=am_0g$ and integrating, we finally get, up to a shift of the axes:  \\
\[ y^1= a \, arcsh\, (x/a) + (T_0/E)x , \hspace{1cm} y^2= \sqrt{a^2+x^2}+ (T_0/2Ea)x^2 \]
At the limit $T\rightarrow \infty$, we get $y^2=a\, ch(y^1/a)$ as expected. However, we may also find the same result by looking at the extremum of the potential energy under the differential constraint $(ds/dx)-1 =0$, that is by studying the variational problem:  \\
\[     \delta {\int}_C(m_0gy^2 + \lambda (\frac{ds}{dx}-1))dx=0    \]
with {\it Lagrange multiplier} $\lambda$  and getting the usual identification $\lambda=T$, though we notice that the founding principles are quite different. \\

\newpage 

\noindent
{\bf 4) FRENET-SERRET FORMULAS REVISITED }  \\

Though the Frenet (1847) and Serret (1851) formulas are well known by anybody studying curves in ${\mathbb{R}}^3$, we shall revisit them in the light of the differential Galois theory as they will provide one of the best examples of the criteria for automorphic systems already presented. For such a purpose, let $X=\mathbb{R}$ with local coordinates $(x)$ and $Y={\mathbb{R}}^3$ with cartesian coordinates $(y^1,y^2,y^3)$ as usual. We consider the Lie group of rigid motions $G=(A,B)$ of dimension $3+3=6$ where $A$ is a {\it constant} orthogonal matrix and $B$ is a {\it constant} vector. The action and its successive prolongations are defined by the formulas:  \\
\[ \bar{y}=A\, y + B \Rightarrow {\bar{y}}_x=A \, y_x, {\bar{y}}_{xx}= A\, y_{xx}, {\bar{y}}_{xxx}=A \, y_{xxx}, ...  \]
Compared to the previous examples, the rather striking fact is that now the only differential invariant of order one is:  \\
\[   \Omega  \equiv {\vec{y}}_x.{\vec{y}}_x\equiv (y^1_x)^2 +(y^2_x)^2 + (y^3_x)^2=\omega  \]
and there are thus {\it only two} differential invarints at order two, namely:  \\
\[   \Gamma\equiv {\vec{y}}_x.{\vec{y}}_{xx}\equiv y^1_xy^1_{xx}+y^2_xy^2_{xx}+y^3_xy^3_{xx}=\gamma= \frac{1}{2}{\partial}_x\omega  , \]
\[   \Sigma \equiv   {\vec{y}}_{xx}.{\vec{y}}_{xx}\equiv (y^1_{xx})^2+(y^2_{xx})^2+(y^3_{xx})^2=\sigma  \]
Indeed, another possibility could be:  \\
\[  R \equiv \mid {\vec{y}}_x\wedge {\vec{y}}_{xx} {\mid }^2 \equiv (y^2_xy^3_{xx}-y^3_xy^2_{xx})^2 + (y^3_xy^1_{xx}-y^1_xy^3_{xx})^2 + 
(y^1_xy^2_{xx}-y^2_xy^1_{xx})^2 = \rho  \]
but it is easy to check the identity:  \\
\[   R= \Omega \Sigma - (\Gamma)^2 \Rightarrow \rho=\omega \sigma -{\gamma}^2   \]
but the resulting equations for the symbol of order two should be:  \\
\[  ((y^3_xy^1_{xx}-y^1_xy^3_{xx})y^3_x -(y^1_xy^2_{xx} - y^2_xy^1_{xx})y^2_x ) v^1_{xx} + ... =0  \]
where the sum is done on the permutations of $(1,2,3)$ and these equations are linear combunations of the two other ones, namely: \\
\[  {\delta}_2\Gamma \equiv y^1_xv^1_{xx} +y^2_xv^2_{xx}+y^3_xv^3_{xx}=0, \, \, {\delta}_2\Sigma \equiv y^1_{xx}v^1_{xx} +y^2_{xx}v^2_{xx}+y^3_{xx}v^3_{xx}=0  \]
because we have $ {\delta}_2 R =\Omega \, {\delta}_2\Sigma - 2 \Gamma \,{\delta}_2\Gamma $ in $V(J_2(X\times Y))$.  \\
In the present situation, ${\cal{A}}_1$ is clearly {\it not} a PHS for $G$ as it has fiber dimension equal to $3+2=5 < 6$ while ${\cal{A}}_2$ is a PHS for $G$ with fiber dimension $3+2+1=6$ but is {\it not} an automorphic system. Indeed, using one prolongation, {\it we only get two third order OD equations} and the correct automorphic system ${\cal{A}}_3$ is defined by adding the following {\it three} third order equations: \\
\[  \Phi \equiv {\vec{y}}_x.{\vec{y}}_{xxx}=y^1_xy^1_{xxx}+y^2_xy^2_{xxx}+y^3_xy^3_{xxx}=\varphi={\partial}_x\gamma- \sigma,\]
\[\Psi\equiv {\vec{y}}_{xx}.{\vec{y}}_{xxx}=y^1_{xx}y^1_{xxx}+y^2_{xx}y^2_{xxx}+y^3_{xx}y^3_{xxx}=\psi= \frac{1}{2}{\partial}_x\sigma, \]
\[ \Upsilon \equiv  ({\vec{y}}_x \wedge {\vec{y}}_{xx}).{\vec{y}}_{xxx} \equiv ({\vec{y}}_x,{\vec{y}}_{xx}, {\vec{y}}_{xxx})= \upsilon  \]
with no CC involving $\upsilon$ , that is only $(\omega,\sigma, \upsilon)$ can be given arbitrarily.

In order to recover the classical formulas, we need to introduce the curvilinear abcissa $s$ defined by $ds^2= \omega \,dx^2$ and we get successively:  \\
\[  \frac{d\vec{y}}{ds}= \vec{t}, \,\, \frac{d\vec{t}}{ds}= \kappa \vec{n}, \,\, \frac{d\vec{n}}{ds}= - \kappa \vec{t} + \tau \vec{b}, 
\,\, \frac{d\vec{b}}{ds}= - \tau \vec{n}  \]
where $\kappa$ is called the {\it curvature}, $\tau$ is called the {\it torsion} and $\vec{b}=\vec{t} \wedge \vec{n}$. As $\mid \vec{t}{\mid}^2=1$, it follows that $\vec{n} $ is orthogonal to $\vec{t}$ and we shall exhibit the relations existing between $(\kappa,\tau)$ and the $(\omega, \sigma, \nu)$ already obtained. \\

Now, with $y=f(x),y_x={\partial}_xf(x), y_{xx}={\partial}_{xx}f(x)$, we obtain successively:  \\
\[  y_x=y_s\frac{ds}{dx}, y_{xx}=y_{ss}(\frac{ds}{dx})^2 + y_s \frac{d^2s}{dx^2}, y_{xxx}=y_{sss}(\frac{ds}{dx})^3+ 3y_{ss}\frac{ds}{dx}\frac{d^2s}{dx^2}+y_s\frac{d^3s}{dx^3}, ...\]
\[ {\vec{y}}_x=\frac{ds}{dx}\vec{t} , {\vec{y}}_{xx}=(\frac{ds}{dx})^2\frac{d\vec{t}}{ds} + \frac{d^2s}{dx^2}\vec{t}, 
{\vec{y}}_{xxx}= (\frac{ds}{dx})^3\frac{d^2 \vec{t}}{ds^2}+ 3\frac{ds}{dx}\frac{d^2s}{dx^2}\frac{d\vec{t}}{ds}+ \frac{d^3s}{dx^3}\vec{t}, ... \]

\[  \omega ={\vec{y}}_x.{\vec{y}}_x=(\frac{ds}{dx})^2 \mid \vec{t}\,{\mid}^2=(\frac{ds}{dx})^2 \]
\[ \gamma= {\vec{y}}_x.{\vec{y}}_{xx}= \frac{ds}{dx}\frac{d^2s}{dx^2}=\frac{1}{2}\frac{d}{dx}(\frac{ds}{dx})^2=\frac{1}{2}{\partial}_x\omega  \]
\[ \rho = \mid {\vec{y}}_x\wedge {\vec{y}}_{xx}{\mid}^2=(\frac{ds}{dx})^6 \mid \vec{t}\wedge \frac{d\vec{t}}{ds}{\mid}^2=
{\omega}^3{\kappa}^2\mid \vec{b}\,{\mid}^2={\omega}^3{\kappa}^2 \]

\[\sigma  ={\vec{y}}_{xx}.{\vec{y}}_{xx}=(\frac{ds}{dx})^4 \frac{d\vec{t}}{ds}.\frac{d\vec{t}}{ds}+ (\frac{d^2s}{dx^2})^2=
{\omega}^2{\kappa}^2 +    ( {\gamma}^2/\omega  )                \]

\[ \varphi = {\vec{y}}_x.{\vec{y}}_{xxx}= {\omega}^2\vec{t}.\frac{d^2\vec{t}}{ds^2} + \frac{ds}{dx}\frac{d^3s}{ds^3}= - {\omega}^2{\kappa}^2 +  \frac{ds}{dx}\frac{d^3s}{ds^3}= - {\omega}^2{\kappa}^2+{\partial}_x\gamma -({\gamma}^2/{\omega})      \]

\[  \begin{array}{rcl}
\upsilon  & = & ({\vec{y}}_x\wedge {\vec{y}}_{xx}). {\vec{y}}_{xxx}                        \\
      & = & (\frac{ds}{dx})^3(\frac{d\vec{y}}{ds}\wedge \frac{d^2\vec{y}}{ds^2}).((\frac{ds}{dx})^3 \frac{d^3\vec{y}}{ds^3}  \\
      & = & {\omega}^3(\vec{t}\wedge \frac{d\vec{t}}{ds}). \frac{d^2\vec{t}}{d^2s} \\
      & = &{\omega}^3 (\vec{t},\frac{d\vec{t}}{ds},\frac{d^2\vec{t}}{ds^2})  \\
      & = & {\omega}^3(\vec{t}\wedge \kappa \vec{n}).(\frac{d\kappa}{ds} \vec{n} + \kappa \frac{d\vec{n}}{ds}) \\
      & = & {\omega}^3\kappa \vec{b}.(-{\kappa}^2\vec{t} + \kappa \tau \vec{b})  \\
      & = & {\omega}^3{\kappa}^2\tau     \\
      & = & \rho \tau 
\end{array}   \]

We now provide the computation for the {\it helix} (draw a picture):  \\
\[ y^1=r \,cos(\theta), \hspace{5mm} y^2=r \,sin(\theta), \hspace{5mm}  y^3=h \, \theta, \hspace{5mm}  r=cst  \]
\[  dy^1=-r\, sin(\theta) d\theta, \hspace{3mm} dy^2=r \, cos(\theta)d\theta, \hspace{3mm} dy^3=h \, d\theta \Rightarrow (ds)^2=(r^2+h^2)(d\theta)^2 \]
\[t^1= - \frac{r\, sin(\theta)}{\sqrt{r^2+h^2}},\hspace{3mm}t^2=\frac{r \, cos(\theta)}{\sqrt{r^2+h^2}},\hspace{3mm}t^3= \frac{h}{\sqrt{r^2+h^2}} \]
\[ dt^1= - \frac{r \, cos(\theta)}{\sqrt{r^2+h^2}}d\theta, \hspace{3mm} dt^2= - \frac{r \, sin(\theta)}{\sqrt{r^2+h^2}}d\theta, \hspace{3mm}dt^3=0 \]
We obtain therefore $\kappa=\mid\frac{d\vec{t}}{ds}\mid =\frac{r}{r^2+h^2}$ and: \\
\[ n^1= - cos(\theta), n^2= - sin(\theta), n^3=0 \Rightarrow \frac{dn^1}{ds}=\frac{sin(\theta)}{\sqrt{r^2+h^2}},  \frac{dn^2}{ds}=\frac{cos(\theta)}{\sqrt{r^2+h^2}},  \frac{dn^3}{ds}=0 \]
A result leading to:\\
\[ \frac{dn^1}{ds}+\kappa t^1=\frac{1}{\sqrt{r^2+h^2}} \frac{h^2}{r^2+h^2}sin(\theta),  \frac{dn^2}{ds}+\kappa t^2= - \frac{1}{\sqrt{r^2+h^2}} \frac{h^2}{r^2+h^2} cos(\theta),\]
\[ \frac{dn^3}{ds}+\kappa t^3=\frac{1}{\sqrt{r^2+h^2}}\frac{hr}{r^2+h^2}  \]
and finally to $\tau=\mid \frac{d\vec{n}}{ds}+\kappa \vec{t} \mid=\frac{h}{r^2+h^2}$.  \\

Studying finally the gauging procedure existing in the differential Galois theory, we may introduce the $3 \times 3$ matrix:  \\
\[M=\left(  \begin{array}{lll}
y^1_x & y^1_{xx} & y^1_{xxx}  \\
y^2_x & y^2_{xx} & y^2_{xxx}  \\
y^3_x & y^3_{xx} & y^3_{xxx} 
\end{array}  \right)  \]
For two sections $f_3$ and ${\bar{f}}_3$ of ${\cal{A}}_3$, we have $\bar{M}=A\, M$ and thus $A(x)=\bar{M}M^{-1}$ whenever $M$ is invertible, 
that is $det(M)=\upsilon\neq 0$. Then $B(x)={\bar{f}}(x) - A(x) f(x)$ and, as before, $A$ and $B$ are constant whenever $y=f(x)$ and $\bar{y}={\bar{f}}(x)$ are two solutions  of ${\cal{A}}_3 $, a result not evident at first sight.  \\

\newpage

\noindent
{ \bf 5) LIE THEORY OF OD EQUATIONS REVISITED}  \\

As a very general algebraic problem, if $K$ and $L$ are two abstract extensions of a field $k$, the idea is to look for the rings and field that can be constructed from this only knowledge, in particular we should look for a bigger fields containing $K$ and $L$ in a nice way in order to define the smallest subfield $(K,L)$ containing $K$ and $L$. Coming back to ordinary Galois theory, the basic purpose is to know about all roots of a given polynomial or, in a more abstract way, to split a given polynomial into as many irreducible factors as possible through the knowledge of the roots of another polynomial. For example, if $K/k$ is a Galois extension and $L$ is an arbitrary extension of $k$ such that $K$ and $L$ are contained in a bigger field, then $(K,L)$ is a Galois extension of $L$. Moreover, if $K \cap L=k$, then $K$ and $L$ are linearly disjoint over $k$ in $(K,L)$ 
([20,22], p 131 or any basic textbook on Galois theory like [1,28]).\\  

As an elementary but tricky counter-example with $k=\mathbb{Q}$, let us consider the polynomials $P\equiv y^3-2\in k[y]$ with residue $y\rightarrow \eta$ and $Q\equiv z^2+z+1\in k[z]$ with residue $z\rightarrow \zeta$. As $P$ and $Q$ are irreducible over $k$, we may consider the fields $K=k(\eta)$ and $L=k(\eta\zeta)$. As $(P,Q)$ is a prime ideal in $k[y,z]$, then we may choose $(K,L)=k(\eta,\zeta)=Q(k[y,z]/(P,Q))=k[y,z]/(P,Q)$. In actual practice, we may choose for $\eta$ the real root $\sqrt[3]{2}$ of $P$ and for $\zeta$ the imaginery cubic root of unity $(-1+i\sqrt{3})/2$ in such a way that 
$\eta\zeta$ is an imaginary root of $P$. We have of course $K\cap L =k$ but $K$ is {\it not} a Galois extension of $k$. Accordingly, $K$ and $L$ are not linearly disjoint over $k$ in $(K,L)$ because we have $(\eta\zeta)^2 \times 1+ (\eta\zeta)\eta + 1 \times {\eta}^2=0$ while $\{1, \eta, {\eta}^2\}$ is a basis of $K$ over $k$. Our aim will be to extend these ideas to the differential framework. \\

While studying the integration of the nonlinear OD equation $\frac{dy}{dx}-F(x,y)=0$, namely looking for functions $y=u(x)$ such that ${\partial}_xu(x)-F(x,u(x))\equiv 0$, Lie discovered that the knowledge of a vector field $\theta= \xi (x,y){\partial}_x+\eta (x,y){\partial}_y$ preserving the OD equation brings the integration to a simple {\it quadrature}. Indeed, introducing the Lie derivative ${\cal{L}}(\theta)=i(\theta)+di(\theta) $ where now $i( )$ is the interior multiplication by a vector and $d$ is the exterior derivative, the invariant property amounts to the equation ${\cal{L}}(\theta)(dy-Fdx)=A(dy-Fdx)$ where $A$ is a multiplicative factor. Eliminating $A$ among the various factors of $dx$ and $dy$, we get the only condition:  \\
\[ \frac{\partial \eta}{\partial x} + F \frac{\partial \eta}{\partial y} - F \frac{\partial \xi}{\partial x} - F^2 \frac{\partial \xi}{\partial y} - \frac{\partial F}{\partial x}\xi - \frac{\partial F}{\partial y}\eta =0     \]
Bringing the terms together in a different way, Lie discovered that this condition can also be written:  \\
\[  \frac{\partial}{\partial x} (\eta - F \xi) +  F\frac{\partial}{\partial y}(\eta - F \xi) - \frac{\partial F}{\partial y} (\eta - F \xi) =0  \]
Finally, setting $\chi = 1/(\eta - F \xi)$ and $\omega = -F$, we obtain equivalently:  \\
\[ \frac{\partial \chi}{\partial x} - \frac{\partial (\omega \chi)}{\partial y}=0  \]
a result proving that $\chi$ is an integrating factor for the $1$-form $dy + \omega dx$, that is to say the $1$-form $\chi(dy+\omega dx)$ is closed. \\

It is not evident at all to establish a link between such a problem and the differential Galois theory by extending the previous purely algebraic comment to a differential framework and using the Spencer operator. First of all, {\it we must start with an automorphic system}. For this, changing slightly the notations while introducing a manifold $X$, say ${\mathbb{R}}^2$ with local coordinates $(x^1,x^2)$ instead of $(x,y)$, and a manifold $Y$, say $\mathbb{R}$ with local coordinate $(y)$, we may consider the automoprphic system ${\cal{A}}_1$ defined over a differential field $K$:  \\
\[     \Phi \equiv \frac{y_1}{y_2}=\omega \in K, \hspace{5mm} y_2\neq 0  \]
for the Lie pseudogroup $\Gamma = aut(Y)$ of invertible transformations $\bar{y}=g(y)$ with $ {\partial}_yg(y)\neq 0$. Then, let us notice that any solution $y=f(x^1,x^2)$ of the preceding system  is such that:  \\
\[      {\partial}_1f + F(x^1,x^2){\partial}_2 f=0, \hspace{5mm}  {\partial}_2f\neq 0  \]
while $F\in K$ where the independent variables $(x^1,x^2)$ are not explicitly appearing. If we have any {\it first integral} $f(x^1,x^2)=c=cst$, we may use the implicit funxction theorem in order to obtain $x^2=u(x^1,c)$ with an {\it identity} $f(x^1,u(x^1,c))\equiv c, \forall x^1$. We obtain therefore:  \\
\[ {\partial}_1f(x^1,u(x^1,c)) + {\partial}_2f(x^1, u(x^1,c)){\partial}_1u(x^1,c)\equiv 0 \Rightarrow {\partial}_2f ( {\partial}_1u-F(x^1,u))=0   \]
and thus ${\partial}_xu(x,c)-F(x,u(x,c))=0, \forall x,\forall c$. As for the Spencer operator, setting $x^1=x$ and introducing any function $u(x)$, we may define a section $u_1=(u(x),u_x(x))$ of the first jet bundle by choosing $u_x(x)=F(x,u(x))$ and the initial system amounts to $Du_1={\partial}_xu(x) - u_x(x)=0$.  \\

From the differential algebraic point of view, with $\omega= - F\in K$ when $K$ is a differential field with derivations $({\partial}_1,{\partial}_2)$, we get the differential automorphic extension $L/K=Q(K\{y\})/{\mathfrak{p}})$ where $\mathfrak{p}\subset K\{y\}$ is the prime linear differential ideal generated by the differential polynomial $P=y_1-\omega y_2\in K\{y\}$. The idea is to exhibit another differential extension $M/K$ with $M=Q(K\{z^1,z^2\}/\mathfrak{q})$ and $\mathfrak{q}\subset K\{z^1,z^2\}$ is the differential ideal generated by the differential polynomial:   \\
\[  Q \equiv z^2_1 - \omega z^2_2 + \omega z^1_1 - {\omega}^2z^1_2 + {\partial}_1 \omega z^1 + {\partial}_2 \omega z^2 \in K\{z^1, z^2\}  \]
in order to take into account the condition for $\theta = {\xi}^1{\partial}_1 + {\xi}^2 {\partial}_2$ already obtained with different notations. However, contrary to the algebraic case, {\it the intersection} $K'=L\cap M$ {\it has no meaning at all}, unless we could define such an intersection in a bigger differential field $N$ containing both $L$ and $M$, according to the following commutative diagram of field inclusions:  \\
\[ \begin{array}{ccl}
L  & \rightarrow & N  \\
\uparrow &   & \uparrow   \\
  \, K' &  \rightarrow & M  \\
\uparrow & \nearrow  &   \\
K &  &  
\end{array}  \]
For this, let us consider the differential extension $N/K$ with 
$N=Q(K\{y,z^1,z^2\}/\mathfrak{r})$ and $\mathfrak{r}\subset K\{y,z^1,z^2\}$ is generated by the two differential polynomials:   \\
\[  y_1 - \omega y_2, \hspace{5mm} y_1z^1 +y_2z^2 - 1   \]
We have of course $\mathfrak{r}\cap K\{y\}=\mathfrak{p} \Rightarrow L \subset N$ and it just remains to prove that $\mathfrak{r} \cap K\{z^1,z^2 \}= \mathfrak{q} \Rightarrow M \subset N$. For this, we obtain by substitution the two differential polynomials :  \\
\[    (z^2 + \omega z^1)y_1 - \omega,\hspace{1cm} (z^2 + \omega z^1)y_2 - 1   \]
The elimination of $y$ can be done by crossed derivatives and we just find for $(z^1,z^2)$ the integrating factor condition for $\chi=1/(z^2 + \omega z^1)$ already obtained. The last relations prove that $K'$ is indeed generated by $(y_1,y_2)$ modulo $\mathfrak{p}$ and we have the new automorphic extension $L/K'$ for the subpseudogroup ${\Gamma}' \subset \Gamma$ made by the translations $\bar{y}=y+ a$ with $a=cst$. The corresponding automorphic system ${\cal{A}}'_1 \subset {\cal{A}}_1$:  \\
\[  {\Psi}^1 \equiv y_1 = \chi \omega, {\Psi}^2 \equiv y_2 = \chi ,   \hspace{5mm} \chi \in K'  \]
is nothing else than a new description of the quadrature concept where we have now $\chi = 1/({\xi}^2 + \omega {\xi}^1)$ and $({\xi}^1,{\xi}^2)$ is the image of $(z^1,z^2)$ under the residue with respect to $\mathfrak{q}$. We finally notice that $L=K'(y)$ ({\it care}) while $M=K'<z^1>=K'<z^2>$. It follows that $L$ is {\it regular} over $K'$, that is $K'$ is algebraically closed in $L$, while $M$ is differentially transcendental over $K'$. It also follows that $L$ and $M$ are linearly disjoint over $K'$ in $N$. Up to our knowledge, such a Galois type approach has never been provided elsewhere.  \\

\newpage

\noindent
{\bf 6) DRACH-VESSIOT THEORY REVISITED}  \\

Roughly, in the preceding example we have used $n=2,m=1$ in order to study the {\it first order} OD equation $y_x -F(x,y)=0$ and our purpose is now to use $n=3,m=2$ in order to study the {\it second order} OD equation $y_{xx} - F(x,y,y_x)=0$. Accordingly, we shall introduce functions $u^k(x,y,y_x)$ for $k=1,2$ and consider the linear homogeneous system:  \\
\[  \frac{\partial u^k}{\partial x}+ y_x \frac{\partial u^k}{\partial y}+ F(x,y,y_x)\frac{\partial u^k}{\partial y_x}=0,  \hspace{7mm}
   \frac{\partial (u^1,u^2)}{\partial (y,y_x)}\neq 0  \]
For any couple of first integrals $u^k(x,y,y_x)=c^k$ with $k=1,2$, setting $c=(c^1,c^2)$, we can use the Jacobian condition and the implicit function theorem in order to solve these  two equations with respect to $(y,y_x)$ in order to get:    \\
\[  y=f( x;c), \hspace{1cm}   y_x=f_x(x;c)  \]
However, differentiating with respect to $x$ the identities $u^k(x,f(x;c),f_x(x;c))\equiv c^k$, we obtain the relations:   \\
\[ \frac{\partial u^k}{\partial x} + \frac{\partial u^k}{\partial y}\frac{\partial f}{\partial x} + 
\frac{\partial u^k}{\partial y_x}\frac{\partial f_x}{\partial x} =0  \]
Substracting the previous equations in order to eliminate the $\partial u^k/\partial x$ while taking into account the Jacobian condition:  \\
\[   \frac{\partial (u^1,u^2)}{\partial (y,y_x)}\neq 0 \Leftrightarrow \frac{\partial(f,f_x)}{\partial (c^1,c^2)}\neq 0  \]
we obtain the {\it Spencer operator} through this {\it vertical procedure}, namely:  \\
 \[      \frac{\partial f}{\partial x} - f_x=0, \hspace{1cm} \frac{\partial f_x}{\partial x} - F=0  \]
      
Changing slightly the notations as before, we may use a manifold $X$ with local coordinates $(x)=(x^1, x^2, x^3)$ instead of $(x,y,y_x)$ with $x^1=x$ and a manifold $Y$ with local coordinates $(y)=(y^1,y^2)$ in order to look for solutions $y^k=f^k(x^1,x^2,x^3)=f^k(x)$ for $k=1,2$ of the linear system:  \\
\[  y^k_1 + x^3 y^k_2 + F(x) y^k_3=0,   \hspace{7mm}  \frac{\partial (y^1,y^2)}{\partial (x^2,x^3)}\equiv y^1_2y^2_3 - y^1_3y^2_2 \neq 0  \]
It is not evident at all and it has been the discovery of Drach ([9]) and Vessiot ([29]), that this is indeed an automorphic system ${\cal{A}}_1\subset J_1(X\times Y)$ for the Lie pseudogroup $\Gamma =aut(Y)$ made by invertible transformations of the form $\bar{y}=g(y)$ with nonzero jacobian, that 
can be written:  \\
\[ {\Phi}^1 \equiv \frac{ \frac{\partial (y^1,y^2)}{\partial (x^1,x^2)}}{\frac{\partial (y^1,y^2)}{\partial (x^2,x^3)}}={\omega}^1(x)=F(x),\hspace{5mm}{\Phi}^1 \equiv \frac{ \frac{\partial (y^1,y^2)}{\partial (x^3,x^1)}}{\frac{\partial (y^1,y^2)}{\partial (x^2,x^3)}}=
{\omega}^2(x)=x^3 \] 
We have thus obtained:    \\

\noindent
{\bf THEOREM 6.1}: The search for a family of solutions of the given second order OD equations depending on $2$ parameters is equivalent to the knowlege of one solution of this automorphic system.  \\

As the two generating first order differential invariants are rational functions of the first jets, we may therefore use the differential Galois theory by introducing a differential field $K$ with derivations $({\partial}_1,{\partial}_2, {\partial}_3)$ in such a way that ${\omega}^1, {\omega}^2 \in K$. Also, the fiber dimension of the system is $2+(3\times 2)-2=6$ while the fiber dimension of the system of Lie equations (made by no equation !) is 
$2+(2\times 2)=6$ and the automorphic property follows from the fact that ${\cal{A}}_1$ is involutive with no CC.  \\

Now, following Jacobi, we shall call a function $M(x)$ (Jacobi) {\it multiplier} for $\theta={\theta}^i{\partial}_i$ if we have the relation ${\partial}_i(M{\theta}^i)=0$. In particular, if $\bar{x}=\varphi (x)$ is a change of independent variables with jacobian 
$\Delta (x)=det ({\partial}_i{\varphi}^j(x)\neq 0$, then it is well known that we have the {\it identity} ([21,23,26]):  \\
\[   \frac{\partial}{\partial {\bar{x}}^j}( \frac{1}{\Delta}\frac{\partial {\varphi}^j}{\partial x^i})\equiv 0, \hspace {1cm} \forall x\in X  \]
Setting ${\bar{\theta}}^j={\partial}_i{\varphi}^j{\theta}^i$, we obtain easily:  \\
\[ \frac{\partial}{\partial {\bar{x}}^j}(\frac{M}{\Delta}{\bar{\theta}}^j)=
 \frac{\partial}{\partial {\bar{x}}^j}    (\frac{M}{\Delta}\frac{\partial {\varphi}^j}{\partial x^i}{\theta}^i)=
 \frac{1}{\Delta}\frac{\partial}{\partial x^i}(M{\theta}^i)=0   \]
and it follows that $M/\Delta$ is a multiplier for $\bar{\theta}$. In particular, $M=1$ is a multiplier if and only if $\theta$ is divergence free, that is 
${\partial}_i{\theta}^i=0$. In the present situation, we should have ${\partial}_1(1) + {\partial}_2 (x^3) + {\partial}_3(F)=0$, that is 
${\partial}_3F=0$ and thus $F=F(x^1,x^2)$.  \\

Coming back to the initial notations, let us consider the second order OD equation $y_{xx} - F(x,y)=0$ and the corresponding system:  \\
\[    \theta.\phi\equiv \frac{\partial \phi}{\partial x} +y_x\frac{\partial \phi}{\partial y}+ F(x,y) \frac{\partial \phi}{\partial y_x}=0  \]
If $\phi(x,y,y_x)=c=cst$ is a first integral containing explicily $y_x$, we can use locally the implicit function theorem and find $y_x=\psi(x,y;c)$. 
Let us prove that the $1$-form $dy - \psi(x,y;c)dx$ has the integrating factor $1/\frac{\partial \phi}{\partial y_x}(x,y,\psi(x,y;c))$ which is also a Jacobi 
multiplier for the differential system $dx/1=dy/\psi(x,y,c)$, that is let us prove:   \\
\[ \frac{\partial}{\partial x}((1/\frac{\partial \phi}{\partial y_x}) + \frac{\partial}{\partial y}(\psi/\frac{\partial \phi}{\partial y_x})=0  \]
whenever $y_x=\psi(x,y;c)$ that is to say:   \\
\[  \frac{{\partial}^2 \phi}{\partial x\partial y_x} +(\frac{\partial \psi}{\partial x} + \psi \frac{\partial \psi}{\partial y})\frac{{\partial}^2\phi}{\partial y_x\partial y_x} + \psi \frac{{\partial}^2 \phi}{\partial y\partial y_x} - \frac{\partial \psi}{\partial y} \frac{\partial \phi}{\partial y_x} =0   \]
Indeed, differentiating the identity $\phi(x,y,\psi(x,y;c))\equiv c$ with respect to $c$, $x$ and $y$ successively, we get:  \\
\[   \frac{\partial \phi}{\partial y_x}\frac{\partial \psi}{\partial c}=1  \Rightarrow  \frac{\partial \phi}{\partial y_x}\neq 0   \]
 \[    \frac{\partial \phi}{\partial x} + \frac{\partial \phi}{\partial y_x}\frac{\partial \psi}{\partial x} =0, \hspace{5mm}  
     \frac{\partial \phi}{\partial y} + \frac{\partial \phi}{\partial y_x}\frac{\partial \psi}{\partial y} =0  
   \, \,  \Rightarrow  \,\, \frac{\partial \psi}{\partial x}+ \psi \frac{\partial \psi}{\partial y}=F(x,y)   \]
and it just remains to differentiate the PD equation satisfied by $\phi$ with respect to $y_x$.  \\

Among the possible reductions of the Galois pseudogroup $\Gamma$, we may consider the Lie sub-pseudogroup ${\Gamma}'=\{ \bar{y}=g(y) \mid \partial({\bar{y}}^1,{\bar{y}}^2)/\partial (y^1,y^2)=1\}$ leading to the automorphic sub-system ${\cal{A}}'_1\subset {\cal{A}}_1$:  \\
\[ {\Psi}^1\equiv \frac{\partial (y^1,y^2)}{\partial (x^2,x^3)}={\psi}^1, \hspace{5mm} {\Psi}^2\equiv \frac{\partial (y^1,y^2)}{\partial (x^3,x^1)}={\psi}^2, \hspace{5mm}{\Psi}^3\equiv \frac{\partial (y^1,y^2)}{\partial (x^1,x^2)}={\psi}^3   \]
where ${\psi}^1,{\psi}^2,{\psi}^3 \in K'$ with $K \subset K' \subset L$ and:  \\
\[  {\partial}_1{\psi}^1 + {\partial}_2 {\psi}^2 + {\partial}_3 {\psi}^3=0  \]
As we must have ${\psi}^2/{\psi}^1= x^3, \, \,  {\psi}^3/{\psi}^1=F(x^1,x^2,x^3) $, {\it this reduction amounts to the explicit knowledge of a Jacobi multiplier}. In the particular case $F=F(x^1,x^2)$, we may choose ${\psi}^1=1$ and {\it this situation is exactly the one obtained by passing from the Lagrangian formalism to the Hamiltonian formalism in analytical mechanics}. We provide the details of this striking result which does not seem to be known. \\

With $n=3,m=2$ in this case, let us consider a Lagrangian $L(t,x,\dot{x})$ and the corresponding Euler-Lagrange equations:  \\
\[ \frac{d}{dt}(\frac{\partial L}{\partial \dot{x}}) - \frac{\partial L}{\partial x}=0 \Leftrightarrow \frac{{\partial}^2L}{\partial t \partial \dot{x}}+\dot{x} \frac{{\partial}^2 L}{\partial x\partial \dot{x}}+\ddot{x} \frac{{\partial}^2 L}{\partial \dot{x}\partial \dot{x}}-
\frac{\partial L}{\partial x}=0 \]
When the {\it Hessian condition} ${\partial}^2L/\partial \dot{x}\partial \dot{x}\neq 0$ is satisfied, we get a second order OD equation of the form 
$\ddot{x} - F(t,x,\dot{x})=0$ and we may thus introduce the automorphic system:  \\
\[  \frac{\partial y^k}{\partial t} + \dot{x} \frac{\partial y^k}{\partial x} + F(t,x, \dot{x}) \frac{\partial y^k}{\partial \dot{x}}=0, \hspace{5mm}
    \forall k=1,2  \]
with: \\
\[ \frac{\partial}{\partial t} (1) +\frac{ \partial}{\partial x} (\dot{x}) +\frac{\partial}{\partial \dot{x}}(F(t,x,\dot{x}))=
\frac{\partial F}{\partial \dot{x}}\neq 0  \]
in general. However, if we now consider the corresponding Hamiltonian formalism obtained by setting $p=\frac{\partial L}{\partial \dot{x}}$ and 
$H=\dot{x} \frac{\partial L}{\partial \dot{x}} - L= H(t,x,p)$, we obtain at once:  \\
\[   dH= \dot{x} dp - \frac{\partial L}{\partial t} dt -\frac{\partial L}{\partial x} dx \]
\[   \frac{dp}{dt}= \frac{{\partial}^2 L}{\partial t \partial \dot{x}} + \dot{x} \frac{{\partial}^2 L}{\partial x \partial \dot{x}}+ 
                            \ddot{x}\frac{{\partial}^2 L}{\partial \dot{x} \partial \dot{x}}=\frac{\partial L}{\partial x}  \]
and thus the well known OD Hamiltonian equations: \\
\[  \frac{dx}{dt} = \frac{\partial H}{\partial p}, \hspace{1cm}  \frac{dp}{dt}= - \frac{\partial H}{\partial x}  \]
a result leading to the automorphic system: \\
\[    \frac{\partial y^k}{\partial t} + \frac{\partial H}{\partial p} \frac{\partial y^k}{\partial x} - 
\frac{\partial H}{\partial x}\frac{\partial y^k}{\partial p}=0, \hspace{1cm}  \forall k=1,2  \]
on which we check:  \\
\[ \frac{\partial}{\partial t} (1) +\frac{\partial }{\partial x} (\frac{\partial H}{\partial p}) + \frac{\partial }{\partial p}(- \frac{\partial H}{\partial x})=0\]
The previous Lagrangian automorphic system had independent variables $(t,x,\dot{x})$ while the new Hamiltonian automorphic system has independent variables $(t,x,p)$. As this latter system admits the Jacobi multiplier $1$, it follows from the general theory explained at the beginning of this example that the Lagrangian system {\it must} admit the Jacobi multiplier:  \\
\[      \frac{\partial (t,x,p)}{\partial (t,x,\dot{x})}=\frac{\partial p}{\partial \dot{x}}=\frac{{\partial}^2L}{\partial \dot{x} \partial \dot{x}}   \]
Indeed, multiplying respectively $(1,\dot{x}, F(t,x,\dot{x}))$ by the hessian and noticing that:  \\
\[    F(t,x,\dot{x}) \frac{{\partial}^2L}{\partial \dot{x} \partial \dot{x}}=\ddot{x}\frac{{\partial}^2L}{\partial \dot{x} \partial \dot{x}}=
\frac{\partial L}{\partial x} - \frac{{\partial}^2 L}{\partial t \partial \dot{x}}- \dot{x} \frac{{\partial}^2 L}{\partial x \partial \dot{x}}  \]
we finally check, after an easy computation, the identity:    \\
\[  \frac{\partial}{\partial t}(\frac{{\partial}^2L}{\partial \dot{x} \partial \dot{x}}) + \frac{\partial}{\partial x} (\dot{x}\frac{{\partial}^2L}{\partial \dot{x} \partial \dot{x}}) + \frac{\partial}{\partial \dot{x}}(\frac{\partial L}{\partial x} - \frac{{\partial}^2 L}{\partial t \partial \dot{x}}- \dot{x} \frac{{\partial}^2 L}{\partial x \partial \dot{x}})\equiv 0   \]
The reduction to ${\Gamma}'\subset \Gamma$ becomes respectively:  \\
\[ \frac{\partial (y^1,y^2)}{\partial (x, \dot{x})}=\frac{{\partial}^2L}{\partial \dot{x} \partial \dot{x}}  ,\hspace{5mm} 
   \frac{\partial (y^1,y^2)}{\partial ( \dot{x}, t)}= \dot{x} \frac{{\partial}^2L}{\partial \dot{x} \partial \dot{x}},\hspace{5mm}
            \frac{\partial (y^1,y^2)}{\partial ( \dot{x}, t)}=\frac{\partial L}{\partial x} - \frac{{\partial}^2 L}{\partial t \partial \dot{x}}- \dot{x} \frac{{\partial}^2 L}{\partial x \partial \dot{x}}   \]
\[     \frac{\partial (y^1,y^2)}{\partial (x,p)}=1, \hspace{5mm} \frac{\partial (y^1,y^2)}{\partial (p,t)}=\frac{\partial H}{\partial p}, \hspace{5mm}
       \frac{\partial (y^1,y^2)}{\partial (t,x)}= - \frac{\partial H}{\partial x}  \]
The knowledge of one first integral, say $y^2=\phi$, reduces the Galois pseudogroup to ${\Gamma}"=\{ {\bar{y}}^1=y^1 + h(y^2), {\bar{y}}^2=y^2 \}$ .  \\

\newpage

\noindent
{\bf 7) HAMILTON-JACOBI EQUATIONS REVISITED}  \\

We end this list of examples by revisiting the Hamilton-Jacobi equation. This is {\it by far} the most difficult example in the sense that no classical approach using exterior calculus can be used in order to exhibit the corresponding automorphic systems involved. At the same time, it uses the first criterion for automorphic systems and thus, in particular, formal integrability or involution become crucial tools that cannot be avoided. This is the reason for which the results we present have not been found during the last century.  \\
Let $z=f(t,x)$ be a solution of the non-linear PD equation $z_t+H(t,x,z,z_x)=0$ written with jet notations for the single unknown $z$. When dealing with applications, $t$ will be {\it time}, $x$ will be {\it space}, $z$ will be the {\it action} and, as usual, we shall set $p=z_x$ for the {\it momentum}. It is important to notice that, in this general setting, $H(t,x,z,p)$ {\it cannot be called Hamiltonian as it involves} $z$. By analogy with the preceding example, we shall set ([9,20-22,29]):    \\

\noindent
{\bf DEFINITION 7.1}: A {\it complete integral} $z=f(t,x;a,b)$ is a family of  solutions depending on two constant parameters $(a,b)$ in such a way that the Jacobian condition $\partial (z,p)/\partial (a,b)\neq 0$ whenever $p={\partial}_xf(t,x;a,b)$. Using the implicit function theorem, we may set\\

\noindent
{\bf THEOREM 7.2}: The search for a complete integral of the PD equation:  \\
 \[z_t+H(t,x,z,z_x)=0 \] 
is equivalent to the search for a {\it single} solution of the automorphic system ${\cal{A}}_1$ with $n=4,m=3$, obtained by eliminating $\rho(t,x,z,p)$ in the Pfaffian system:  \\
\[    dz-pdx+H(t,x,z,p) dt = \rho (dZ-PdX )  \]  
The corresponding Lie pseudogroup is the pseudogroup $\Gamma$ of {\it contact transformations} of $(X,Z,P)$ that reproduces the contact $1$-form $dZ-PdX$ up to a function factor.  \\

\noindent
{\it Proof}: If $z=f(t,x; a,b)$ is a complete integral, we have:  \\
\[   dz-pdx+H(t,x,z,p)dt= \frac{\partial f}{\partial a}da +\frac{\partial f}{\partial b} db   \]
Using the implicit function theorem and the Jacobian condition, we may set:  \\
\[  a =X(t,x,z,p), \, b=Z(t,x,z,p)  \Rightarrow \rho(t,x,z,p)=\frac{\partial f}{\partial b}, \,
 P(t,x,z,p)=\frac{\partial f}{\partial a}/\frac{\partial f}{\partial b}   \]
The converse is left to the reader.  \\

For another solution denoted wit a "bar", we have:  \\
\[ dz-pdx+H(t,x,z,p) dt = \bar{\rho} (d\bar{Z}-\bar{P} d\bar{X} )\,\,\Rightarrow \,\,d\bar{Z}-\bar{P} d\bar{X} = \frac{\rho}{\bar{\rho}}(dZ-PdX) \]
Closing this system, we obtain at once:  \\
\[    d\bar{X}\wedge d\bar{Z}\wedge d\bar{P}= (\frac{\rho}{\bar{\rho}})^2 dX\wedge dZ \wedge dP    \]
Closing again, we discover that $\rho/\bar{\rho}$ is in fact a function of $(X,Z,P)$, a result bringing the Lie pseudogroup of contact transformations and showing that no restriction must be imposed to $H(t,x,z,p)$.  \\
\hspace*{12cm}     Q.E.D.  \\

It is quite more dificult to exhibit the equations of the above automorphic sytem and the corresponding equations of the Lie pseudogroup $\Gamma$ in Lie form or even as involutive systems of PD equations. From what has been said, we obtain {\it at least}:   \\
\[  \frac{\frac{\partial \bar{Z}}{\partial X} - \bar{P}\frac{\partial \bar{X}}{\partial X}}{ \frac{\partial \bar{Z}}{\partial Z} - \bar{P}\frac{\partial \bar{X}}{\partial Z}}= - P , 
    \frac{\frac{\partial \bar{Z}}{\partial P} - \bar{P}\frac{\partial \bar{X}}{\partial P}}{ \frac{\partial \bar{Z}}{\partial Z} - \bar{P}\frac{\partial \bar{X}}{\partial Z}}= 0     \Rightarrow  \frac{\partial \bar{Z}}{\partial P} - \bar{P}\frac{\partial \bar{X}}{\partial P}=0  \]
 for defining ${\cal{R}}_1$, that is to say:  \\
 \[    \frac{ \partial \bar{Z}}{\partial X} - \bar{P} \frac{\partial \bar{X}}{\partial X} + P  (\frac{\partial \bar{Z}}{\partial Z} - \bar{P}\frac{\partial \bar{X}}{\partial Z})=0 , \hspace{1cm} \frac{\partial \bar{Z}}{\partial P} - \bar{P}\frac{\partial \bar{X}}{\partial P}=0   \]
Using now letters $(x,z,p)$ instead of the capital letters $(X,Z,P)$ and $(\xi, \eta, \zeta)$ for the corresponding vertical bundles, we obtain by linearization the system of first order infinitesimal Lie equations:  \\
\[  \frac{\partial \xi}{\partial x}- p \frac{ \partial \eta}{\partial  x} - \zeta + p ( \frac{\partial \xi}{\partial z}- p\frac{ \partial \eta}{\partial  z})=0, \,\, 
\frac{\partial \xi}{ \partial p }- p \frac{\partial \eta}{\partial p}=0  \]
This system is not involutive as it is not even formally integrable. Using crossed derivatives in $x/p$, we obtain the {\it only new first order} equation: \\
\[    \frac{\partial \eta}{\partial x} - \frac{\partial \xi}{\partial z}+ \frac{\partial \zeta }{\partial p} + 2p \frac{\partial \eta}{\partial z}=0  \]
and the resulting system ${\cal{R}}^{(1)}_1$ is involutive with two equations of class $x$ solved with respect to $(\frac{\partial \xi}{\partial x}, \frac{\partial \eta}{\partial x})$ and one equation of class $p$ solved with respecto $\frac{\partial \xi}{\partial p}$, that is $dim_Y({\cal{R}}^{(1)}_1)= (3+ 3 \times 3) - 3=9$. Accordingly, the non-linear system of Lie equations {\it must} become involutive by adding {\it only one equation in Lie form}, namely:   \\
\[   \frac{  \frac{\partial (\bar{Z},\bar{X}, \bar{P})}{\partial (Z,X,P)}}{ ( \frac{\partial \bar{Z}}{\partial Z} - \bar{P} \frac{\partial \bar{X}}{\partial X})^2} = 1   \]
and its linearization jus provides:  \\
\[  \frac{\partial \eta}{\partial x} + \frac{\partial \xi}{\partial z}+ \frac{\partial \zeta }{\partial p}= 2  (\frac{\partial \xi}{\partial z} - p \frac{\partial \eta}{\partial z})  \]
that is exactly the previous equation. The following {\it Janet board} provides the structure of an involutive {\it solved form}:  \\ 
\[  \left\{   \begin{array}{lcl}
  X  & \longrightarrow & \left\{  \begin{array}{c} \bar{Z} \\ \bar{X} \end{array}\right.  \\
   P  &  \longrightarrow &  \, \left\{  \begin{array}{c} \bar{Z}  \end{array}  \right.
 \end{array}  \right.                         
\fbox{ $ \begin{array}{lll}
 Z & P & X \\
 Z & P & X \\
 Z & P & \bullet \\
\end{array}$ }   \]
                            
Coming back to the original system and notations, we may suppose $\frac{\partial Z}{\partial z} -P\frac{\partial X}{\partial z}\neq 0$ and introduce the $7=3+4$ equations: \\
\[ \frac{\partial Z}{\partial x} - P\frac{\partial X}{\partial x} + p ( \frac{\partial Z}{\partial z} - P\frac{\partial X}{\partial z})=0, 
   \frac{\partial Z}{\partial t} - P\frac{\partial X}{\partial t} - H ( \frac{\partial Z}{\partial z} - P\frac{\partial X}{\partial z})=0,
   \frac{\partial Z}{\partial p} -  P \frac{\partial X}{\partial p}=0  \]
\[   \frac{\partial (Z,X,P)}{\partial (z,x,p)} - (\frac{\partial Z}{\partial z} - P\frac{\partial X}{\partial z})^2=0,
     \frac{\partial (Z,X,P)}{\partial (z,p,t)} - \frac{\partial H}{\partial p} (\frac{\partial Z}{\partial z} - P\frac{\partial X}{\partial z})^2=0, ...  \]  
                                
Starting now, the next results {\it canot} be obtained by exterior calculus and are therefore not known. Indeed, developping the $ 3 \times 3$ Jacobian determinant, the fourth equation provided can be written as:  \\
\[  \frac{\partial Z}{\partial x}. \frac{\partial (X,P)}{\partial (x,p)} - \frac{\partial Z}{\partial x} . \frac{\partial ( (X,P)}{\partial (z,p)} + \frac{\partial Z}{\partial p}. \frac{\partial (X,P)}{\partial (z,x)} - (\frac{\partial Z}{\partial z} - P\frac{\partial X}{\partial z})^2=0  \]
Using the previous equations in order to eliminate $\frac{\partial Z}{\partial x}$ and $\frac{\partial Z}{\partial p}$, we obtain:  \\
\[  \frac{\partial Z}{\partial x}. \frac{\partial (X,P)}{\partial (x,p)} + p(\frac{\partial Z}{\partial z}-P \frac{\partial X}{\partial z}) . \frac{\partial ( (X,P)}{\partial (z,p)}  - P\frac{\partial X}{\partial x}. \frac{\partial (X,P)}{\partial (z,p)}       + P \frac{\partial X}{\partial p}. \frac{\partial (X,P)}{\partial (z,x)}  =\] 
\[    (\frac{\partial Z}{\partial z} - P \frac{\partial X}{\partial z})( \frac{\partial (X,P)}{\partial (x,p)} + p \frac{\partial (X,P)}{\partial (z,p)})   = (\frac{\partial Z}{\partial z} - P\frac{\partial X}{\partial z})^2  \]
and thus:  \\
\[ \frac{\partial (X,P)}{\partial (x,p)} + p \frac{\partial (X,P)}{\partial (z,p)})  -  (\frac{\partial Z}{\partial z} - P\frac{\partial X}{\partial z}) =0   \]
which is nothing else than the first order equation that can be obtained from the first and third among the previous $7$ equations by using crossed derivatives in $x/p$. It follows that ${\cal{A}}^{(1)}_1$ may be defined by $6$ equations {\it only} and we have thus 
$dim_X ({\cal{A}}^{(1)}_1)= (3 + 4 \times 3)- 6=9$. The following {\it Janet board} provides the structure of an involutive {\it solved form}:  \\ 

\[  \left\{   \begin{array}{lcl}
  x  & \longrightarrow & \left\{  \begin{array}{c} Z \\ X \\ P \end{array}\right.  \\
  t  & \longrightarrow  & \left\{  \begin{array}{c} Z \\ X  \end{array} \right. \\
  p  &  \longrightarrow &  \, \left\{  \begin{array}{c} Z   \end{array}  \right.
   \end{array} \right.
\fbox{ $ \begin{array}{llll}
z & p & t & x \\
z & p & t & x \\
z & p & t & x \\
z & p & t & \bullet \\
z & p & t & \bullet \\
z & p  &\bullet &  \bullet
\end{array}$ }   \]
 
This result proves that the involutive system ${\cal{A}}^{(1)}_1$ is an automorphic system for the involutive Lie groupoid ${\cal{R}}^{(1)}_1$. \\

If $H=H(t,x,p)$ is an Hamiltonian function, then we may look for  a complete integral of the form $z=f(t,x;a)+b$ and we have:  \\

\noindent
{\bf COROLLARY 7.3}: The search for such a complete integral of the PD equation:  \\
 \[z_t+H(t,x,z_x)=0 \] 
is equivalent to the search for a {\it single} solution of the automorphic system ${\cal{A}}'_1$ with $n=4,m=3$ described by the Pfaffian system:  \\
\[    dz-pdx+H(t,x,p) dt = dZ-PdX   \]  
The corresponding Lie pseudogroup is the Lie pseudogroup ${\Gamma}'\subset \Gamma$ of {\it unimodular contact transformations} of $(X,Z,P)$ that preserve the contact $1$-form $dZ-PdX$ and we have thus $\partial(\bar{X},\bar{Z},\bar{P})/\partial (X,Z,P)=1$. \\

\noindent
{\it Proof}: Now, we have:   \\
 \[p=\frac{\partial f}{\partial x}(t,x;a)  \Rightarrow a=X(t,x,p)\Rightarrow b=Z(t,x,z,p)=z-\varphi (t,x,p), \rho(t,x,z,p)=1\]
and thus 
\[\frac{\partial Z}{\partial z}=1, \frac{\partial X}{\partial z}=0, \frac{\partial f}{\partial b}=1 \Rightarrow P=\frac{\partial f}{\partial a}(t,x;X(t,x,p))\Rightarrow \frac{\partial P}{\partial z}=0.  \]
We shall just prove that the Pfaffian system:  \\
\[    dz-pdx+H(t,x,z,p) dt = dZ-PdX   \]
is compatible if and only if $\partial H/\partial z=0$ as a new group theoretical justification for revisiting the mathematical foundations of analytical mechanics.\\
Indeed, closing the system, we get:  \\
\[  dx\wedge dp +dH\wedge dt=dX\wedge dP   \]
Using the exterior multiplication among the coresponding left and right members, we get:  \\
\[   dz\wedge dx\wedge dp + dz\wedge dH \wedge  dt- pdx\wedge dH\wedge dt +H dt\wedge dx\wedge dp = dZ\wedge dX \wedge dP   \]
Closing again, we finally obtain:  \\
\[   - dp\wedge dx \wedge dH \wedge dt + dH \wedge dt \wedge dx \wedge dp= 2 \frac{\partial H}{\partial z} dt \wedge dx \wedge dz \wedge dp=0  \]
and the desired condition on $H$.  \\
We invite the reader to discover this condition just using CC for the second member of the desired automorphic system and notice that ${\cal{A}}'_1$ 
{\it is neither involutive nor even formally integrable}. \\
\hspace*{12cm}   Q.E.D.  \\

It is again dificult to exhibit the equations of the above automorphic sytem and the corresponding equations of the Lie pseudogroup ${\Gamma}'$ in Lie form or even as involutives systems of PD equations. From what has been said, we obtain {\it at least}:   \\
\[  \frac{\partial \bar{Z}}{\partial Z} - \bar{P}\frac{\partial \bar{X}}{\partial Z}= 1, \frac{\partial \bar{Z}}{\partial X} - \bar{P}\frac{\partial \bar{X}}{\partial X}= - P , 
     \frac{\partial \bar{Z}}{\partial P} - \bar{P}\frac{\partial \bar{X}}{\partial P}=0  \]
 for defining ${\cal{R}}'_1$, that is to say:  \\
 \[    \frac{ \partial \bar{Z}}{\partial X} - \bar{P} \frac{\partial \bar{X}}{\partial X} + P  (\frac{\partial \bar{Z}}{\partial Z} - \bar{P}\frac{\partial \bar{X}}{\partial Z})=0   \]
 Now, contrary to the preceding situation, we have the Pfaffian system:  \\
 \[   dX\wedge dP =d\bar{X}\wedge d\bar{P}   \]
 and we may add the $3$ new first order equations:  \\
 \[\frac{\partial (\bar{X},\bar{P})}{\partial (X,P)}=1,\frac{\partial (\bar{X},\bar{P})}{\partial (P,Z)}=0,
 \frac{\partial (\bar{X},\bar{P})}{\partial (Z,X)}=0  \] 
 an obtain therefore the $6$ equations:  \\
\[  \frac{\partial \bar{Z}}{\partial Z}=1, \frac{\partial \bar{X}}{\partial Z}=0, \frac{\partial \bar{P}}{\partial Z}=0, 
\frac{\partial \bar{Z}}{\partial X} - \bar{P}\frac{\partial \bar{X}}{\partial X}= - P ,\frac{\partial (\bar{X},\bar{P})}{\partial (X,P)}=1,
     \frac{\partial \bar{Z}}{\partial P} - \bar{P}\frac{\partial \bar{X}}{\partial P}=0     \]
The following {\it Janet board} provides the structure of an involutive {\it solved form}: \\
 \[  \left\{   \begin{array}{lcl}
  Z  & \longrightarrow & \left\{  \begin{array}{c} \bar{Z} \\ \bar{X}\\ \bar{P} \end{array}\right.  \\
  X  & \longrightarrow  & \left\{  \begin{array}{c} \bar{Z} \\ \bar{X} \\ \end{array} \right. \\
  P  &  \longrightarrow &  \, \left\{  \begin{array}{c} \bar{Z}  \end{array}  \right.
 \end{array}  \right.                         
\fbox{ $ \begin{array}{lll}
 P & X & Z \\
 P & X & Z \\
 P & X & Z \\
 P & X & \bullet \\
 P & X  &\bullet \\
P & \bullet & \bullet
\end{array}$ }   \]

\noindent                            
The resulting system ${\cal{R}}'^{(1)}_1$ is thus involutive with $3$ equations of class $Z$ solved with respect to $(\frac{\partial \bar{Z}}{\partial Z}, \frac{\partial \bar{X}}{\partial Z},\frac{\partial \bar{P})}{\partial Z}$, $2$ equations of class $X$ solved with respect to 
$(\frac{\partial \bar{Z}}{\partial X}, \frac{\partial \bar{X}}{\partial X})$ and $1$ equation of class $P$ solved with respect to $(\frac{\partial \bar{Z}}{\partial P})$. The characters are $(2,1,0)$ and we have $dim_Y({\cal{R}}'^{(1)}_1)= (3+ 9) - 6=6$. \\

Taking into account the initial Pfaffian system $dz-pdx+H(t,x,p)dt=dZ-PdX$ and its exterior closure $dx\wedge dp +dH\wedge dt=dX\wedge dP$, we may proceed as before and find the first order system combining the $6$ equations:  \\
\[  \frac{\partial Z}{\partial z}=1, \frac{\partial X}{\partial z}=0, \frac{\partial P}{\partial z}=0, 
\frac{\partial Z}{\partial x} - P \frac{\partial  X}{\partial x}= - p , 
     \frac{\partial Z}{\partial p} - P \frac{\partial X}{\partial p}=0 , \frac{\partial Z}{\partial t} - P \frac{\partial X}{\partial t}=H  \]
with the $3$ equations:  \\
\[  \frac{\partial (X,P)}{\partial (x,p)}=1,\, \, \, \frac{\partial (X,P)}{\partial (x,t)}=\frac{\partial H}{\partial x}, \,  \, \, \frac{\partial (X,P)}{\partial (p,t)}=\frac{\partial H}{\partial p} \] 
The following {\it Janet board} provides the structure of an involutive {\it solved form}:   \\
\[  \left\{   \begin{array}{lcl}
  z &  \longrightarrow  & \left\{ \begin{array}{c} Z \\ X \\ P \end{array}\right. \\
  x  & \longrightarrow & \left\{  \begin{array}{c} Z \\ X \\ P \end{array}\right.  \\
  p  & \longrightarrow  & \left\{  \begin{array}{c} Z \\ X  \end{array} \right. \\
  t  &  \longrightarrow &  \, \left\{  \begin{array}{c} Z   \end{array}  \right.
   \end{array} \right.
\fbox{ $ \begin{array}{llll}
t & p & x & z \\
t & p & x & z \\
t & p & x & z \\
t & p & x & \bullet \\
t & p & x &  \bullet  \\
t & p  & x & \bullet \\
t & p & \bullet & \bullet  \\
t & p & \bullet & \bullet \\
t & \bullet & \bullet & \bullet
\end{array}$ }   \]
with $3$ equations of class $z$, $3$ equations of class $x$, $2$ equations of class $p$ and $1$ equation of class $t$ giving characters $(2,1,0,0)$ where the non-zero ones coincide with the non-zero ones of ${\cal{R}}^{'(1)}_1$. One must therefore use the involutive automorphic system 
${\cal{A}}'^{(1)}_1={\pi}^2_1 ({\cal{A}}'_2)$ where ${\cal{A}}'_2$ is the first prolongation of ${\cal{A}}'_1$. A similar difficulty has been found for the nonlinear Lie equations ${\cal{R}}'_1$ defining the Lie pseudogroup ${\Gamma}'\subset \Gamma$.  \\

\noindent
{\bf COROLLARY 7.4}: The search for a complete integral $z=u(t;a) + v(x;a) +b$ by {\it separation of variables} is equivalent to the search for a single solution of the automorphic system obtained by adding $\partial X/\partial t=0$ to the system of the preceding Corollary, provided that:  \\
\[      \frac{\partial H}{\partial z}=0 ,\hspace{1cm}  \frac{\partial}{\partial t} ( \frac{\partial H}{\partial x} / \frac{\partial H}{\partial p})=0   \] 
The corresponding Lie pseudogroup is:  \\
\[ {\Gamma}"= \{\bar{X}=g(X), \bar{Z}=Z+ h(X), \bar{P}=(P + \partial h/\partial X)/ (\partial g/ \partial X)\} \subset {\Gamma}' \subset \Gamma  \]

\noindent
{\it Proof}: We have $p=\partial v/\partial x \Rightarrow a=X(x,p) \Rightarrow \partial X / \partial t=0 $ and:  \\
\[b=z-u(t;X(x,p))- v(x;X(x,p))=Z(t,x,p)  \Rightarrow \frac{\partial Z}{\partial t}= \frac{\partial u}{\partial t}(t;X(x,p))=H(t,x,p)     \]
\[ \frac{\partial f}{\partial b}=1\Rightarrow P= \frac{\partial f}{\partial a}= \frac{\partial v}{\partial a}(x; X(x,p)) \Rightarrow \frac{\partial P}{\partial t}=0  \]
Using the last three equations of the preceding Corollary, we get:  \\
\[\frac{\partial X}{\partial t}=\frac{\partial H}{\partial p} / \frac{\partial X}{\partial x} - \frac{\partial H}{\partial x} / \frac{\partial X}{\partial p}=0\]        
We obtain therefore the additional differential invariant in Lie form:  \\
\[   \frac{\partial X}{\partial x} / \frac{\partial X}{\partial p} = \frac{\partial H}{\partial x} / \frac{\partial H}{\partial p}  \]
and the desired CC by differentiating with respect to $t$. This {\it necessary condition} for separating the variables in the integration of the Hamilton-Jacobi equation has been found by T. Levi-Civita in $1904$ ([17]) and integrated byA. Huaux in $1976$ ([12]).  \\
It also follows that the corresponding Lie pseudogroup ${\Gamma}"$ must contain transformations $\bar{X}=g(X)$. Then $\partial \bar{Z}/\partial Z=1, \partial \bar{Z}/\partial P=0 \Rightarrow \bar{Z}=Z + h(X)$ and finally:  \\
\[\frac{\partial \bar{P}}{\partial P}\frac{\partial \bar{X}}{\partial X}=1, \frac{\partial \bar{Z}}{\partial X} - \bar{P}\frac{\partial \bar{X}}{\partial X}= - P \Rightarrow  \bar{P}= (P + \partial h/\partial X)/ (\partial g/ \partial X) \].
\hspace*{12cm}   Q.E.D.   \\

\noindent
{\bf COROLLARY 7.5}: The search for a complete integral $z=u(t) + v(x;a) + at + b$ is equivalent to the search for $1$ solution of the automorphic system obtained by adding $\partial P/\partial t=1$ to the system of the preceding Corollary, provided that:  \\
\[\frac{\partial H}{\partial z}=0, \hspace{1cm} \frac{{\partial}^2 H}{\partial t\partial x}=0,\hspace{1cm} \frac{{\partial}^2 H}{\partial t \partial p}=0\]
The corresponding Lie pseudogroup is:  \\
\[   {\Gamma}'''= \{ \bar{X}=X+c, \bar{Z}=Z + h(X), \bar{P}= P + \partial h/\partial X \}   \]

{\it Proof}: We have:  \\
\[  p=\partial v / \partial x (x;a) \Rightarrow a=X(x,p), \hspace{5mm}  Z=b=z- u(t) - v(x;X(x,p)) -X(x,p)t \]
and:    \\
\[  \frac{\partial f}{\partial b}=1 \Rightarrow P=\frac{\partial v}{\partial a }(x; X(x,p)) + t \Rightarrow \frac{\partial P}{\partial t}=1  \]
\[ \frac{\partial Z}{\partial t}= - \frac{\partial u}{\partial t}(t) - X(x,p)=H(t,x,p) \Rightarrow    \frac{\partial H}{\partial t}= -\frac{{\partial}^2 u}{\partial t^2}(t)  \]
Collecting all these results, the system is defied by the $11$ equations;  \\
\[   \frac{\partial Z}{\partial z}=1, \frac{\partial X}{\partial z}=0, \frac{\partial P}{\partial z}=0, \frac{\partial X}{\partial t}=0,  
\frac{\partial Z}{\partial t}=H , \frac{\partial P}{\partial t}=1 , \frac{\partial X}{\partial x}=\frac{\partial H}{\partial x}, \frac{\partial X}{\partial p}=\frac{\partial H}{\partial p}, \]
\[   \frac{\partial Z}{\partial x} -  \frac{\partial  H}{\partial x}P= - p , \frac{\partial Z}{\partial p} -  \frac{\partial H}{\partial p} P=0 , 
\frac{\partial H}{\partial x}\frac{\partial P}{\partial p} - \frac{\partial H}{\partial p}\frac{\partial P}{\partial x}= 1  \]
 The following {\it Janet board} provides the structure of an involutive solved form:  \\ 
\[  \left\{   \begin{array}{lcl}
  z &  \longrightarrow  & \left\{ \begin{array}{c} Z \\ X \\ P \end{array}\right. \\
  x  & \longrightarrow & \left\{  \begin{array}{c} Z \\ X \\ P \end{array}\right.  \\
  t & \longrightarrow  & \left\{  \begin{array}{c} Z \\ X \\ P \end{array} \right. \\
  p  &  \longrightarrow &  \, \left\{  \begin{array}{c} Z \\ X  \end{array}  \right.
   \end{array} \right.
\fbox{ $ \begin{array}{llll}
p & t & x & z \\
p & t & x & z \\ 
p & t & x & z \\
p & t & x & \bullet \\
p & t & x &  \bullet  \\
p & t  & x & \bullet \\
p & t & \bullet & \bullet  \\
p & t & \bullet & \bullet \\
p & t & \bullet & \bullet \\
p & \bullet & \bullet & \bullet  \\
p & \bullet & \bullet & \bullet
\end{array}$ }   \]
and the fiber dimension of this non-linear involutive first order system is $(3+4\times 3)-11= 4$ with characters equal to $(1,0,0,0)$.  \\
As for the non-linear system of finite Lie equations, we have at once from the explicit transformations:  \\
\[  \frac{\partial \bar{Z}}{\partial Z}=1, \frac{\partial \bar{X}}{\partial Z}=0, \frac{\partial \bar{P}}{\partial Z}=0, 
\frac{\partial \bar{Z}}{\partial P}=0,  \frac{\partial \bar{X}}{\partial P}= 0 ,\frac{\partial \bar{P}}{\partial P}=1,
     \frac{\partial \bar{X}}{\partial X}=1, \frac{\partial \bar{Z}}{\partial X}-\bar{P}= - P   \]                                                                                                                                                                                                                                                                                                                                                                                                                                                                                                                                                                                                                                                                                                                                                                                                                                                                                                                                                                                                                                                                                                                                                                                                                                                                                                                                                                                                                                                                                                                                                                                                                                                                                                                                                                                                                                                                                                                                                                                                                                                                                                                                                                                                                                                                                                                                                                                                                 
This system is involutive and the following {\it Janet board} provides the structure of an involutive solved form:  \\
 \[  \left\{   \begin{array}{lcl}
  Z  & \longrightarrow & \left\{  \begin{array}{c} \bar{Z} \\ \bar{X}\\ \bar{P} \end{array}\right.  \\
  P  & \longrightarrow  & \left\{  \begin{array}{c} \bar{Z} \\ \bar{X} \\ \bar{P} \\ \end{array} \right. \\
  X  &  \longrightarrow &  \, \left\{  \begin{array}{c} \bar{Z}\\ \bar{X}  \end{array}  \right.
 \end{array}  \right.                         
\fbox{ $ \begin{array}{lll}
 X & P & Z \\
 X & P & Z \\
 X & P & Z \\
 X & P & \bullet \\
 X & P & \bullet \\
 X & P & \bullet  \\
X & \bullet & \bullet \\
X & \bullet & \bullet
\end{array}$ }   \]
The fiber dimension at order $1$ is $(3+3\times 3)-8=4$ and the characters are $(1,0,0)$. It follows that we have again an automorphic system. It is worthwhile to notice that the Janet boards of a system and a subsystem may be quite different.  \\

We end this list of examples with a situation providing an intransitive groupoid in the case of a cyclic variable and let the reader treat it as an exercise.  \\

\noindent
{\bf COROLLARY 7.6}: The search for a complete integral $z=v(x;a)+at +b$ is equivalent to the search for $1$ solution of the automorphic system obtained by adding $X=H $ to the system of the preceding Corollary, provided that:  \\
\[  \frac{\partial H}{\partial z}=0, \hspace{1cm} \frac{\partial H}{\partial t}=0   \]
The corresponding intransitive Lie pseudogroup is:  \\
\[   {\Gamma}''''= \{ \bar{X}=X, \bar{Z}=Z + h(X), \bar{P}= P + \partial h/\partial X \} \subset {\Gamma }'''  \subset \Gamma \]

For more details on these topics of analytical mechanics, the interested reader may look at the recapitulating board in ([21], p 506).  As a striking conclusion, there are as many specific situations reflected by the hamiltonian as the number of Lie subpseudogroups of the Lie pseudogroup of contact transformations.   \\

\newpage

\noindent
{\bf REFERENCES}  \\

\noindent
[1] Artin, E.: Galois Theory, Notre Dame Mathematical Lectures, 2, 1942, 1997.  \\
\noindent
[2] Bonnet,O.: M\'emoire sur la Th\'{e}orie des Surfaces Applicables sur une Surface Donn\'{e}e, Journal de l'Ecole Polytechnique, 25 (1867) 31-151. \\
\noindent
[3] Byalinicki-Birula,A.: On the field of rational functions of Algebraic Groups, Pacific Journal of Mathematics, 11 (1961) 1205-1210.  \\
\noindent
[4] Bialynicki-Birula, A.: On Galois Theory of Fields with Operators, Amer. J. Math., 84 (1962) 89-109.  \\
\noindent
[5] Cartan,E.: Sur la Possibilit\'{e} de Plonger une Surface Riemannienne dans un Espace Euclidien, Ann. Soc. Pol. Math., 6 (1927)  1-7.  \\
\noindent
[6] Chevalley, C., Samuel, P.: Twoo Proofs of a Theorem on Algebraic groups, Proc. Amer. Math. Soc., 2 (1951) 126-134;  \\
\noindent
[7] Ciarlet, P.G., Larsonneur, F.P.: On the Ricovery of a Surface with Prescribed First and Second Forms, J. Math. Pures Appl., 81 (2002) 167-185. \\
\noindent
[8] Codazzi, D.: Sulle Coordinate Curvilinee d'una Superficiedello Spazio, Ann. Math. Pura Applicata, 2 (1868-1869) 101-19  \\
\noindent
[9] Drach, J.: Th\`{e}se de Doctorat, Ann. Ec. Normale Sup. (3) 15 (1898) 243-384.   \\
\noindent
[10] Gauss, K.F.: Disquitiones Generales circa Superficies Curvas, Comm. Soc. Gott., 6 (1828).  \\
\noindent
[11] Han, Q.: Global Isometric Embedding of Surfaces in ${\mathbb{R}}^3$, Chapter 2 in " Differential geometry and Coninuum mechanics", Springer, 2015.  \\
\noindent
[12] Huaux, A.: Sur la S\'eparation des Variables dans l'Equation aux d\'eriv\'ees partielles de Hamilton-Jacobi, Annali di Matematica Pura ed Applicata, 108 (1976) 251-282.  \\
\noindent
[13] Janet, M: Sur la Possibilit\'{e} de Plonger une Surface Riemannienne dans un Espace Euclidien, Ann. Soc. Pol. Math., 5 (1926)  38-43. \\
\noindent
[14] Kaplanski, I.: An Introduction to Differential Algebra, Hermann, 1957, 1976.  \\
\noindent
[15] Kolchin,, E.R. KOLCHIN: Differential Algebra and Algebraic groups, Academic Press, 1973.  \\
\noindent
[16] Kolchin, E.R.: Differential Algebraic Groups, Pure and Applied Mathematics, 114, Academic Press, 1985.  \\
\noindent
[17] Levi-Civitta; : Sulla Integrazione delle Equazione di Hamilton-Jacobi per Separazione di Variabili, Math. Annalen, 59,(1904) 383 (Also "Opere Matematische", 2 (1901-1907) 395-410).  \\
\noindent
[18] Mainardi, G.: Su la Theoria Generale delle Superficie, Giornale dell' Instituto Lombardo, 9 (1856) 385-404.   \\
\noindent
[19] Pommaret, J.F.: Systems of Partial Differential Equations and Lie Pseudogroups, Gordon and Breach, New York, 1978 
(Russian translation by MIR, Moscow, 1983). \\
\noindent
[20] : Pommaret, J.F.Differential Galois Theory, Gordon and Breach, New York, 1983 (750 pp).\\
\noindent
[21] Pommaret, : Lie Pseudogroups and Mechanics, Gordon and Breach, New York, 1988 (590 pp).  \\
\noindent
[22] Pommaret, J.F.: Partial Differential Equations and Group Theory: New Perspectives for Applications, Kluwer, 1994.\\
http://dx.doi.org/10.1007/978-94-017-2539-2  \\
\noindent 
[23] Pommaret, J.F.: Partial Differential Control Theory, Kluwer, 2001 (957 pp).\\
\noindent
[24] Pommaret, J.F.: Algebraic Analysis of Control Systems Defined by Partial Differential Equations, in Advanced Topics in Control Systems Theory, Lecture Notes in Control and Information Sciences LNCIS 311, Chapter 5, Springer, 2005, 155-223.\\
\noindent
[25] Pommaret, J.F.: Relative Parametrization of Linear Multidimensional Systems, Multidim. Syst. Sign. Process., 26 (2015) 405-437.  \\
DOI 10.1007/s11045-013-0265-0   \\
\noindent
[26] Pommaret, J.F.: Deformation Theory of Algebraic and Geometric Structures, Lambert Academic Publisher (LAP), Saarbrucken, Germany, 2016.  \\
\noindent
[27] Ritt, J.F.: Differential Algebra, Dover, 1950, 1966.  \\
\noindent
[28] Stewart, I.: Galois Theory, Chapman and Hall, 1973.  \\
\noindent
[29] Vessiot, E.: Sur la Th\'{e}orie de Galois et ses Diverses G\'{e}n\'{e}ralisations, Ann.Ec. Normale Sup., 21 (1904) 9-85 (Can be obtained from  http://numdam.org).  \\
\noindent
[30] Zariski, O., Samuel, P: Commutative Algebra, Van Nostrand, 1958.  \\

\end{document}